\documentclass[reqno,11pt]{amsart}
\UseRawInputEncoding
\pdfoutput=1
\usepackage{amsfonts}
\usepackage{graphicx}
\usepackage{amsfonts,amsmath, amssymb}
\usepackage{bbm,mathrsfs,mathscinet}
\usepackage{hyperref}
\usepackage{marginnote}
\usepackage{color}
\usepackage{mathrsfs}
\allowdisplaybreaks

\numberwithin{equation}{section}

\newcommand{\fb}{\mathfrak{b}}
\newcommand{\fa}{\mathfrak{a}}
\newcommand{\fw}{\mathfrak{w}}
\newcommand{\R}{\mathbb{R}}

\newcommand{\tw}{\tilde{w}}

\newtheorem{theorem}{Theorem}[section]

\newtheorem{lemma}[theorem]{Lemma}
\newtheorem{proposition}[theorem]{Proposition}
\newtheorem{remark}[theorem]{Remark}

\def\v{\varepsilon}

\def\t{\theta}
\def\k{\kappa}

\def\g{\gamma}

\def\r{\rho}

\def\z{\zeta}
\def\o{\omega}

\def\i{\infty}
\def\f{\frac}
\def\pa{\partial}

\def\dis{\displaystyle}
\def\la{\langle}
\def\ra{\rangle}
\def\intr{\int_{\mathbb{R}^{3}}}
\def\ints{\int_{\mathbb{S}^{2}}}

\def\sp{\shortparallel}

\begin{document}
	
	\title[Incompressible Euler limit of Boltzmann equation]{Hilbert expansion of the Boltzmann equation in the incompressible Euler level in a channel}

\author[F. M. Huang]{Feimin Huang}
\address[F.M. Huang]{Academy of Mathematics and Systems Science, Chinese Academy of Sciences, Beijing 100190, China; School of Mathematical Sciences, University of Chinese Academy of Sciences, Beijing 100049, China.}
\email{fhuang@amt.ac.cn}

   	\author[W. Q. Wang]{Weiqiang Wang}
   \address[W. Q. Wang]{Academy of Mathematics and Systems Science, Chinese Academy of Sciences, Beijing 100190, China.}
   \email{wangweiqiang@amss.ac.cn}

	\author[Y. Wang]{Yong Wang}
	\address[Y. Wang]{Academy of Mathematics and Systems Science, Chinese Academy of Sciences, Beijing 100190, China; School of Mathematical Sciences, University of Chinese Academy of Sciences, Beijing 100049, China.}
	\email{yongwang@amss.ac.cn}
	
	\author[F. Xiao]{Feng Xiao}
    \address[F. Xiao]{Key Laboratory of Computing and Stochastic Mathematics (Ministry of Education), School of Mathematics and Statistics, Hunan Normal University, Changsha, Hunan 410081, P. R. China}
    \email{xf@hunnu.edu.cn}

\begin{abstract}
The study of hydrodynamic limit of the Boltzmann equation with physical boundary is a challenging problem due to appearance of the viscous and Knudsen boundary layers. In this paper, the hydrodynamic limit from the Boltzmann equation with specular reflection boundary condition to the incompressible Euler in a channel is investigated. Based on the multi-scaled Hilbert expansion, the equations with boundary conditions and compatibility conditions for interior solutions, viscous and Knudsen boundary layers are derived under different scaling, respectively. Then some uniform estimates for the interior solutions, viscous and Knudsen boundary layers are established. With the help of $L^2-L^\infty$ framework and the uniform estimates obtained above, the solutions to the Boltzmann equation are constructed by the truncated Hilbert expansion with multi-scales, and hence the hydrodynamic limit in the incompressible Euler level is justified.

\end{abstract}

\date{\today}
\maketitle
	
\setcounter{tocdepth}{2}
\tableofcontents
	
\thispagestyle{empty}

	
	\section{Introduction and main results}
	
\subsection{Introduction}
In this paper, we consider the scaled Boltzmann equation
\begin{equation}\label{1.0}
	\mathcal{S}_{\mathfrak{t}}\, \partial_{t}\mathscr{F}+v\cdot\nabla_x \mathscr{F}=\frac1{\mathscr{K}_n}Q(\mathscr{F}, \mathscr{F}),
\end{equation}
where $\mathscr{F}(t,x,v)\geq 0$ is the density distribution function for the gas particles with position $x=(x_1,x_2,x_3)\in \mathbb{T}^2\times (0,1)$, where $\mathbb{T}$ is the periodic interval $[0,1]$, and velocity $v=(v_1,v_2,v_3)\in\mathbb{R}^3$ at  time $t>0$, and $\mathscr{K}_n>0$ is Knudsen number, which is proportional to the mean free path, $\mathcal{S}_{\mathfrak{t}}$ is Strouhal number. The Boltzmann collision term $Q(F_1,F_2)$ is defined in the following bilinear form
\begin{align}\label{1.2}
	Q(F_1,F_2)&\equiv\int_{\mathbb{R}^3}\int_{\mathbb{S}^2} B(v-u,\t)F_1(u')F_2(v')\,{d\omega du}
	-\int_{\mathbb{R}^3}\int_{\mathbb{S}^2} B(v-u,\t)F_1(u)F_2(v)\,{d\omega du}
\end{align}
where the relationship between the post-collision velocity $(v',u')$ of two particles with the pre-collision velocity $(v,u)$ is given by
\begin{equation*}
	u'=u+[(v-u)\cdot\omega]\omega,\quad v'=v-[(v-u)\cdot\omega]\omega,
\end{equation*}
for $\omega\in \mathbb{S}^2$, which can be determined by  the conservation laws of momentum and energy
\begin{equation*}
	u'+v'=u+v,\quad |u'|^2+|v'|^2=|u|^2+|v|^2.
\end{equation*}
The Boltzmann collision kernel $B=B(v-u,\theta)$ in \eqref{1.2} depends only on $|v-u|$ and $\theta$ with $\cos\theta=(v-u)\cdot \omega/|v-u|$. Throughout the present paper,  we  consider the hard sphere model, i.e.,
\begin{equation*}
	B(v-u,\t)=|(v-u)\cdot \omega|.
\end{equation*}

From the von Karman relation, we know that the Reynolds number $\mathscr{R}_{e}$ is given by
\begin{equation}\label{1.2-1}
	\frac{1}{\mathscr{R}_{e}}=\frac{\mathscr{K}_{n}}{\mathscr{M}_{a}},
\end{equation}
where $\mathscr{M}_{a}$ is the Mach number. It was shown in Maxwell \cite{Maxwell} and Boltzmann \cite{Boltzmann} that the Boltzmann equation is closely related to the fluid dynamical systems for both compressible and incompressible flows. In fact, what is called hydrodynamic limit is to derive fluid systems from Boltzmann equation by taking appropriate scaling limits of $\mathscr{M}_{a}$ and $\mathscr{K}_{n}$. For instance, formally,  by setting $\mathscr{M}_{a}=\mathcal{S}_{\mathfrak{t}}=\mathscr{K}_n \to 0+$, one can obtain the incompressible Navier-Stokes equations with $\mathscr{R}_{e}=1$. One can derive the compressible Euler equations by setting $\mathscr{M}_{a}=\mathcal{S}_{\mathfrak{t}}=1$ and $\mathscr{K}_n\to 0+$. If taking $\mathscr{M}_{a}=\mathcal{S}_{\mathfrak{t}} \to 0+,\, \mathscr{K}_n \to 0+$
and $\mathscr{R}_{e}\to +\infty$, then we shall derive the incompressible Euler equations.

In the study of hydrodynamic limit, two kinds of classical expansions are usually used. The one is Hilbert expansion, which is proposed by Hilbert in 1912. The another is Enskog-Chapman expansion, which is independently proposed by Enskog and Chapman in 1916 and 1917, respectively. Either Hilbert or the Chapman-Enskog expansion gives formal derivations of the compressible and incompressible fluid equations. How to justify rigorously these asymptotic expansions is a challenging problem, and is closely related to Hilbert's sixth problem \cite{Hilbert}.

For the case of compressible Euler limit, that is, $\mathscr{M}_{a}=\mathcal{S}_{\mathfrak{t}}\sim O(1)$ and $\mathscr{K}_{n}\to 0+$. When the solution of compressible Euler equations is smooth, Caflisch \cite{Caflisch} rigorously established the compressible Euler limit from the Boltzmann equation through the truncated Hilbert expansion, see also \cite{Lachowicz, Nishida, Ukai-Asano} and \cite{Guo Jang Jiang-1, Guo Jang Jiang} via a $L^2-L^{\infty}$ framework. Recently, when the boundary effects are considered, by analyzing systematically the effects of viscous boundary layer and Knudsen layer, Guo, Huang and Wang \cite{Guo-Huang-Wang} justified rigorously the validity of Hilbert expansion of the Boltzmann with specular boundary conditions in half space. On the other hand, it is well known that there are three basic wave patterns for compressible Euler equations, that is, shock wave, rarefaction wave and contact discontinuity, the hydrodynamic limit of Boltzmann to such wave patterns had been established \cite{Huang-Jiang-Wang,Huang-Wang-Yang,Huang-Wang-Yang-1,Huang-Wang-Wang-Yang,Xin-Zeng,Yu} in one spatial dimensional case. We also refer to \cite{WWZ} for the case of planar rarefaction wave in three spatial dimensional case.

For the case of incompressible Euler limit, that is, $\mathscr{M}_{a}=\mathcal{S}_{\mathfrak{t}}\to 0, \mathscr{K}_{n}\to 0+$ and $\mathscr{R}_{e}\to +\infty$. Under some assumptions, Bouchut, Golse and Pulvirenti \cite{Bouchut-Golse-Pulvirenti} established the limit from the renormalized solutions of Boltzmann equation \cite{Diperna-Lions} to the dissipative weak solution \cite{Lions-1996} of incompressible Euler equations without heat transfer. Then, one of the assumptions in \cite{Bouchut-Golse-Pulvirenti} has been removed by Lions and Masmoudi \cite{Lions-Masmoudi-1, Lions-Masmoudi-2}. Later, all assumptions are removed by Saint-Raymond \cite{Saint-Raymond2003, Saint-Raymond2009} through relative entropy method. Moreover, the heat transfer and specular reflection boundary conditions were also studied in \cite{Saint-Raymond2009}, see also \cite{BGP-2012} for the case of Maxwell reflection boundary conditions. Recently, the case of complete diffusive boundary condition is established by Cao, Jang and Kim \cite{Cao-Jang-Kim} for analytic data when $\mathscr{K}_{n}$ satisfies a special relation with $\mathscr{M}_{a}$, see also \cite{Jang-Kim}. On the other hand, when the solution of incompressible Euler equations is smooth, Masi, Esposito and Lebowitz \cite{Masi-Esposito-Lebowitz-1989} justified the incompressible Euler limit of the Boltzmann equation on tours by truncated Hilbert expansion, see also \cite{Wu-Zhou-Li-2019} via a $L^2-L^{\infty}$ framework.

For the hydrodynamic limit from Boltzmann equation to Navier-Stokes(-Fourier) equations, there are also many works, such as \cite{Bardos,Bardos-2,Bardos-Ukai,Golse-Saint-Raymond,Grad,Guo2006,Guo-Liu} without physical boundary, and \cite{ELM-1994,E-Guo-M,E-Guo-K-M,Jang-Kim,Jiang-Masmoudi,Masmoudi-Raymond,Wu} with physical boundary. We shall not go into details about the Navier-Stokes(-Fourier) limits since we will focus on the incompressible Euler limit in present paper.

The purpose of present paper is to establish the hydrodynamic limit from the Boltzmann equation with specular reflection boundary condition to the incompressible Euler equations through Hilbert expansion. Hereafter, let $\Omega:=\mathbb{T}^2\times (0,1)$ and we denote by $\vec{n}_0=(0,0,-1)$ and $\vec{n}_1=(0,0,1)$ the outward normal vector of $\mathbb{T}^2\times [0,1]$ on $x_3=0$ and $x_3=1$ respectively. The phase boundary of $\Omega\times\mathbb{R}^3$ is denoted as $\gamma:=\partial\Omega\times\mathbb{R}^3$, which can be split into outgoing boundary $\gamma_+^{i}$, incoming boundary $\gamma_-^{i}$, and grazing boundary $\gamma_0^{i}$:
\begin{align}
	\begin{split}\nonumber
		\gamma_+^{i}=\{(x,v) : x_{3}=i, v\cdot \vec{n}_{i}=(-1)^{i+1}v_3>0\},\\
		\gamma_-^{i}=\{(x,v) : x_{3}=i, v\cdot \vec{n}_{i}=(-1)^{i+1}v_3<0\},\\
		\gamma_0^{i}=\{(x,v) : x_{3}=i, v\cdot \vec{n}_{i}=(-1)^{i+1}v_3=0\},
	\end{split}\qquad \text{with } i=0,1.
\end{align}
In this paper, we consider the Boltzmann equation with the specular reflection boundary condition:
\begin{align}\label{1.3}
	\mathscr{F}^{\v}(t,x,v)|_{\gamma_{-}^{i}}=\mathscr{F}^{\v}(t,x_{\sp},i,R_xv),\qquad \text{with }i=0,1,
\end{align}
where we have denoted $R_xv=(v_1,v_2,-v_3)$.

\subsection{Asymptotic expansion}
Throughout the present paper, we denote by $\mu(v)$ the normalized global Maxwellian, i.e.,
\begin{equation}\nonumber
	\mu(v)=\frac{1}{(2\pi)^{\frac32}} \exp\left\{-\frac{|v|^2}{2} \right\}.
\end{equation}
It is clear that the global Maxwellian $\mu(v)$ satisfies the specular reflection boundary conditions \eqref{1.3}. Noting \eqref{1.2-1}, we consider the case $\mathscr{R}_{e}\to +\infty$ (i.e., the viscosity goes to zero) as $\mathscr{K}_{n}\to 0+$. On the other hand, in order to use Hilbert expansion, we need the Knudsen number $\mathscr{K}_{n}$ to be an integer power of the thickness of viscous boundary layer $\frac{1}{\sqrt{\mathscr{R}_{e}}}$. Hence we assume that
\begin{align*}
	\mathcal{S}_{\mathfrak{t}}=\mathscr{M}_a=\v^{n-2},\quad \mathscr{K}_n=\v^n,\quad n\geq 3.
\end{align*}
Then the Boltzmann equation \eqref{1.0} is rewritten as
\begin{equation}\label{1.1}
	\v^{n-2}\partial_t \mathscr{F}^\v+v\cdot \nabla_x \mathscr{F}^\v =\frac{1}{\v^n} Q(\mathscr{F}^\v, \mathscr{F}^\v).
\end{equation}

\subsubsection{Interior expansion:} We define the interior expansion as
\begin{equation}\label{1.8}
	F^\v(t,x,v)\sim \mu+\mathscr{M}_{a}\sum\limits_{i=0}^{\infty}\v^{i}F_{i}(t,x,v)=\mu + \v^{n-2}\sum_{i=0}^{\infty}\v^{i}F_{i}(t,x,v).
\end{equation}
Substituting \eqref{1.8} into  \eqref{1.1}, we get
\begin{align}\label{1.8-1}
	&\sum_{i=0}^{\infty}\v^{n-2+i} \pa_tF_{i} + \sum_{i=0}^{\infty}\v^{i}v\cdot\nabla_x F_{i}\nonumber\\
	&=\sum_{i=0}^{\infty}\v^{-n+i} [ Q(\mu, F_{i}) + Q(F_i,\mu) ]
	+\sum_{i,j=0}^{\infty} \v^{i+j-2} Q(F_i, F_j).
\end{align}
Comparing the order of $\v$ in \eqref{1.8-1}, one obtains
\begin{align}\label{1.7-1}
	\begin{split}
		&\dis \v^{-n+i} (0\leq i\leq n-3): \quad
		0=Q(\mu,F_{i})+Q(F_i,\mu),\\[2mm]
		&\dis \v^{-2}:\qquad\qquad\quad 0 =Q(\mu,F_{n-2})+Q(F_{n-2},\mu)+Q(F_0,F_0),\\[2mm]
		&\dis \v^{-1}:\qquad\qquad\quad 0 =[Q(\mu,F_{n-1})+Q(F_{n-1},\mu)]+[Q(F_0,F_1)+Q(F_1,F_0)],\\[2mm]
		&\dis \v^0: \qquad v\cdot\nabla_x F_{0}=[ Q(\mu,F_{n})+Q(F_{n},\mu) ] + [ Q(F_{0},F_{2})+Q(F_2,F_0)] +Q(F_1,F_1),\\[2mm]
		&\dis \qquad\qquad\qquad\qquad\quad\quad\quad\ \vdots\\
		&\dis \v^{k}(k\geq 1):\quad   \pa_t F_{k+2-n}+  v\cdot\nabla_xF_{k}=[Q(\mu,F_{k+n})+Q(F_{k+n},\mu)]\\
		&  \qquad\qquad\qquad\qquad\qquad\qquad\qquad\qquad+\sum_{\substack{i+j=k+2\\ i,j\geq 0}} \f12[Q(F_{i},F_{j}) + Q(F_{j},F_{i})],
	\end{split}
\end{align}
whereafter we use the notations $F_{i}\equiv0$ with $2-n\leq i\leq -1$ for simplicity of presentation.

For later use, we define the linearized collision operator $\mathbf{L} $ by
\begin{equation}
	\mathbf{L}  g:=-\frac{1}{\sqrt{\mu }}[Q(\mu ,\sqrt{\mu } g)+Q(\sqrt{\mu } g,\mu )],\nonumber
\end{equation}
and the nonlinear operator
\begin{equation}\nonumber
	\Gamma(g_1,g_2):=\frac{1}{\sqrt{\mu}} Q(\sqrt{\mu}g_1,\sqrt{\mu}g_2).
\end{equation}
The null space $\mathcal{N}$ of $\mathbf{L}$ is generated by the following standard orthogonal bases
$$
\begin{aligned}
&\chi_{0}(v)\equiv \sqrt{\mu},\quad\chi_{i}(v)\equiv v_i \sqrt{\mu}\,\,\,i=1,2,3,\quad\chi_{4}(v)\equiv\f{1}{\sqrt{6}}(|v|^2-3)\sqrt{\mu}.
\end{aligned}
$$
We also define the collision frequency $\nu$:
\begin{equation}\nonumber
	\nu(v):=\intr\ints B(v-u,\t)\mu(u)d\o du\cong 1+|v|.
\end{equation}
Let $\mathbf{P}g$ be the $L_{v}^{2}$ projection with respect to $[\chi_0,...,\chi_4]$. It is well-known that there exists a positive number $c_{0}>0$ such that for any function $g$
\begin{equation}\label{1.7-2}
	\la \mathbf{L}g,g\ra\geq c_{0}\|\{\mathbf{I-P}\}g\|_{\nu}^{2},
\end{equation}
where the weighted $L^2$-norm $\|\cdot\|_{\nu}$ is defined as
\begin{equation}\nonumber
	\|g\|_{\nu}^2:=\int_{\Omega\times\mathbb{R}^3} g^2(x,v) \nu(v)dxdv.
\end{equation}

For each $i\geq 0$, let $f_{i}:=\frac{F_{i}}{\sqrt{\mu}}$. We define the macroscopic and microscopic part of  $f_{i}$ as
\begin{align}\nonumber
	f_{i}=\mathbf{P} f_{i}+(\mathbf{I}-\mathbf{P}) f_{i}=\{\rho_{i}+u_{i} \cdot v+\frac{1}{2} \theta_{i}(|v|^{2}-3)\} \sqrt{\mu}+(\mathbf{I}-\mathbf{P}) f_{i},
\end{align}
where $\mathbf{P}f_{i}:=\{\rho_{i}+u_{i} \cdot v+\frac{1}{2} \theta_{i}\left(|v|^{2}-3\right)\}$. It follows from \eqref{1.7-1} that
\begin{equation}\label{2.1}
	\begin{aligned}
		&f_{i}\equiv \mathbf{P} f_{i}=\{\rho_{i}+v\cdot u_{i}+\frac{1}{2}\t_{i}(|v|^2-3)\}\sqrt{\mu} \quad \text{for } 0 \leq i \leq n-3,\\
		&(\mathbf{I}-\mathbf{P}) f_{n-2} =\mathbf{L}^{-1}\left(\Gamma(f_{0},f_{0})\right), \\
        &\dis \qquad\qquad\qquad\ \vdots\\
		&\begin{aligned}
	(\mathbf{I}-\mathbf{P}) f_{k+n-1} =\mathbf{L}^{-1}\Big\{&(\mathbf{I}-\mathbf{P})[-\partial_{t} f_{k+1-n}-v \cdot \nabla_{x} f_{k-1}]\\
	&+\sum_{i+j=k+1 \atop i, j \geq 0} \frac{1}{2}[\Gamma(f_{i},f_{j})+\Gamma(f_{j},f_{i})]\Big\}\quad k\geq 0.
	\end{aligned}
	\end{aligned}
\end{equation}

We define the Burnett functions
$\mathcal{A}_{ij}$ and $\mathcal{B}_{i}$ as
\begin{align}\label{2.4-1}
	\mathcal{A}_{ij}:=\{v_iv_j-\delta_{ij}\f{|v|^2}{3}\}\sqrt{\mu},\quad
	\mathcal{B}_{i}:=\f{v_i}{2}(|v|^2-5)\sqrt{\mu},
\end{align}
where $\delta_{ij}$ is the Kronecker symbol. It follows from a direct calculation that
\begin{equation}\label{2.2}
	\left\{\begin{aligned}
		&\int_{\mathbb{R}^{3}} F_{k} d v=\rho_{k},\quad  \int_{\mathbb{R}^{3}} v_{i} F_{k} d v=u_{k, i},\quad \int_{\mathbb{R}^{3}}|v|^{2} F_{k} d v=3\left(\rho_{k}+\theta_{k}\right) ,\\
		&\int_{\mathbb{R}^{3}} v_{i} v_{j} F_{k} d v=\delta_{i j}\left(\rho_{k}+\theta_{k}\right)+\left\langle\mathcal{A}_{i j},(\mathbf{I}-\mathbf{P}) f_{k}\right\rangle, \\
		&\int_{\mathbb{R}^{3}} v_{i}|v|^{2} F_{k} d v=5 u_{k, i}+2\left\langle\mathcal{B}_{i},(\mathbf{I}-\mathbf{P}) f_{k}\right\rangle.
	\end{aligned}\right.
\end{equation}

Multiplying $\eqref{1.7-1}_4$ by $1, v$ respectively, and integrating the resultant equation over $\R^3$ with respect to $v$, we can get
\begin{equation}\label{2.4}
	\operatorname{div}_{x}u_{0}=0,\qquad \nabla_{x}(\r_0+\t_0)=0.
\end{equation}
Thus, we may assume
\begin{equation}\label{2.5-1}
	\r_{0}+\t_{0}={p}_{0}(t),
\end{equation}
where ${p}_{0}(t)$ is a function independent of $x$.  In fact, we can prove $p_{0}(t)$ is a constant and $p_{0}(t)\equiv (\rho_{0}+\theta_{0})(0)$ in Lemma \ref{lem2.2} below due to the compatibility condition for $u_{n-2}$.

Moreover, taking $k=n-2$ in $\eqref{1.7-1}_5$, and multiplying it by $v, |v|^2$ respectively, then integrating the resultant equation over $\R^3$ with respect to $v$, one obtains that
\begin{equation}\label{2.8}
	\left\{\begin{aligned}
		&\partial_{t} u_{0}+\left(u_{0} \cdot \nabla_{x}\right) u_{0}+\nabla_{x}(\rho_{n-2}+\theta_{n-2}-\frac{1}{3}\left|u_{0}\right|^{2})=0, \\
		&\partial_{t} \theta_{0}+\left(u_{0} \cdot \nabla_{x}\right) \theta_{0}=\frac{2}{5}\pa_{t}p_{0}(t)\equiv 0,
	\end{aligned}\right.
\end{equation}
which, together with \eqref{2.4} and \eqref{2.5-1}, yields that $(\rho_{0},u_{0},\theta_{0})$ satisfies the following incompressible Euler system
\begin{equation}\label{2.8-1}
	\left\{\begin{aligned}
		&\mathrm{div}u_{0}=0, \quad \rho_0+\theta_{0}=p_{0}(t)\equiv \rho_{0}(0)+\theta_{0}(0),\\
		&\partial_{t} u_{0}+(u_{0} \cdot \nabla_{x}) u_{0}+\nabla_{x}p_{n-2}=0, \\
		&\partial_{t} \theta_{0}+(u_{0} \cdot \nabla_{x}) \theta_{0}=\frac{2}{5}\pa_{t}p_{0}(t)\equiv 0,
	\end{aligned}\right.
\end{equation}
where we have used the notation $p_{n-2}:=\rho_{n-2}+\theta_{n-2}-\frac{1}{3}\left|u_{0}\right|^{2}$ and the fact $p_{0}(t)$ is a constant. The detailed calculations of \eqref{2.8-1} are given in Appendix \ref{AppendixA}.

For the incompressible Euler equations \eqref{2.8-1}, we impose the slip boundary condition
\begin{equation}\label{2.8-2}
	u_{0}\cdot \vec{n}\vert_{\partial \Omega}=0, \quad \text{with } \vec{n}=\vec{n}_{0}\text{ at }x_{3}=0,\,\,\,\vec{n}=\vec{n}_{1}\text{ at }x_{3}=1,
\end{equation}
and the initial data
\begin{equation}\label{2.8-3}
	(\rho_0,u_{0},\theta_0)(0,x)=(\rho_{0}(0),u_{0}(0),\theta_{0}(0))\in H^{s_{0}}(\Omega).
\end{equation}
Similarly, to guarantee the compatibility condition for the equation of velocity higher order, we supplement the restriction:
\begin{equation}\label{2.8-4}
\int_{\Omega}p_{n-2}(t,x)dx=\mathfrak{c}_{n-2}(t)\in L^{\infty}(0,\tau),
\end{equation}
where $\mathfrak{c}_{n-2}(t)$ will be determined in Lemma \ref{lem2.2} below.

It is well known that there exists a unique local smooth solution $(\rho_{0},u_{0},\theta_{0},p_{n-2})$ of \eqref{2.8-1}-\eqref{2.8-4} with $s_{0}\geq 3$ such that
\begin{equation*}
	(\rho_{0},u_{0},\theta_{0})\in L^{\infty}(0,\tau; H^{s_0}(\Omega)),\quad p_{n-2}\in L^{\infty}(0,\tau; H^{s_0+1}(\Omega)),\quad \int_{\Omega}p_{n-2}dx=\mathfrak{c}_{n-2}(t).
\end{equation*}
see Lemma \ref{lem1.1} for details.

\subsubsection{Viscous boundary layer expansion:}
For viscous boundary layer expansion near $x_{3}=0$ and $x_{3}=1$, we define the scaled normal coordinates:
\begin{equation}\nonumber
	y_{-}:=\frac{x_3}{\v}>0,\qquad y_{+}:=\frac{1-x_{3}}{\v}>0.
\end{equation}
For simplicity of presentation, we denote
\begin{equation}\nonumber
	x_{\shortparallel}=(x_1,x_2),\quad  \nabla_{\sp}=(\partial_{x_1},\partial_{x_2})\quad\mbox{and}\quad  v_{\sp}=(v_1,v_2).
\end{equation}

For specular boundary conditions, it follows from the formal analysis that the zero-order viscous boundary layer will not appear. Therefore, motivated by \cite{Sone-2002, Sone-2007}, we define the upper and lower viscous boundary layer expansion as
\begin{align*}
F^{\pm,\v}(t,x_\sp, y_{\pm},v)\sim \sum_{i=1}^\infty \v^{n-2+i} \bar{F}_{i}^{\pm}(t,x_\sp, y_{\pm},v).
\end{align*}
Since $y_{+}, y_{-}\in (0,+\infty)$, we shall always use $y$ to represent $y_{+}$ and $y_{-}$ for simplicity of notations in the following, and it will be clear from the context which one is being used.

Plugging $F^\v+\bar{F}^{\pm,\v}$  into the Boltzmann equation \eqref{1.1} and using \eqref{1.7-1}, one obtains that
\begin{align}\label{1.13}
	&\sum_{i=1}^{\infty}\v^{n-2+i} \pa_t \bar{F}_{i}^{\pm}+ \sum_{i=1}^{\infty}\v^{i}v_{\sp}\cdot\nabla_{\sp} \bar{F}_{i}^{\pm}\mp \sum_{i=1}^{\infty}\v^{i-1} \, v_3 \pa_{y} \bar{F}_{i}^{\pm}\nonumber\\
	&=\sum_{i=1}^{\infty}\v^{-n+i} [ Q(\mu, \bar{F}_{i}^{\pm}) + Q(\bar{F}_{i}^{\pm},\mu) ]
	+\sum_{i=0}^{\infty} \sum_{j=1}^{\infty} \v^{i+j-2} [Q(F_{i}, \bar{F}_{j}^{\pm}) + Q(\bar{F}_{j}^{\pm},F_{i})]\nonumber\\
	&\qquad + \sum_{i,j=1}^{\infty} \f12 \v^{i+j-2} [ Q(\bar{F}_{i}^{\pm}, \bar{F}_{j}^{\pm}) + Q(\bar{F}_{j}^{\pm}, \bar{F}_{i}^{\pm})].
\end{align}
Comparing the order of $\v$ in \eqref{1.13}, we can obtain
\begin{align}\label{1.14}
	\begin{split}
		&	\dis   \v^{-n+i} (1\leq i\leq n-2): \qquad\quad 0=Q(\mu ,  \bar{F}_{i}^{\pm}) + Q(\bar{F}_{i}^{\pm}, \mu), \\[2mm]
		&\dis  \v^{-1}:\quad 0=[Q(\mu, \bar{F}_{n-1}^{\pm})+Q(\bar{F}_{n-1}^{\pm},\mu)]+ [Q(F_{0}^{\pm}, \bar{F}_{1})+Q( \bar{F}_{1}^{\pm},F_0^{\pm})],\\[2mm]
		&\dis  \v^{0}:\quad \mp v_3 \frac{\partial\bar{F}_{1}^{\pm}}{\partial y }=[Q(\mu,\bar{F}_{n}^{\pm})+Q(\bar{F}_{n}^{\pm},\mu)]
		+ \sum_{\substack{i+j=2\\i\geq0,j\geq1}} [Q(F_i^{\pm}, \bar{F}_{j}^{\pm})+Q(\bar{F}_{j}^{\pm}, F_i^{\pm})]\\
		&\qquad\qquad\qquad+ \sum_{\substack{i+j+l=2\\i\geq0,j\geq1, 1\leq l\leq \fb}} \frac{(\mp y)^l}{l!} [Q(\partial_{x_3}^lF_i^{\pm},\bar{F}_{j}^{\pm})+Q(\bar{F}_{j}^{\pm}, \partial_{x_3}^l F_i^{\pm})] + Q(\bar{F}_{1}^{\pm}, \bar{F}_{1}^{\pm}), \\
		&\dis \quad\quad\quad\quad\quad\quad\quad\ \vdots\\
		\dis & \v^{k}(k\geq 1): \quad   \partial_t \bar{F}_{k+2-n}^{\pm}+v_\sp\cdot\nabla_\sp\bar{F}_{k}^{\pm}\mp v_3 \cdot \frac{\pa\bar{F}_{k+1}^{\pm}}{\partial y} \\
		&\qquad=[Q(\mu,\bar{F}_{k+n}^{\pm})+Q(\bar{F}_{k+n}^{\pm},\mu)] +\sum_{\substack{i+j=k+2\\ i\geq 0, j\geq 1}}  [Q(F_i^{\pm}, \bar{F}_{j}^{\pm})+Q( \bar{F}_{j}^{\pm},F_i^{\pm}) ]\\
		&\qquad\qquad+\sum_{\substack{i+j+l=k+2\\i\geq0,j\geq 1,    1\leq l\leq \fb,}}\frac{(\mp y)^l}{l!} [Q(\partial_{x_3}^l F_i^{\pm}, \bar{F}_{j}^{\pm})+Q( \bar{F}_{j}^{\pm},\partial_{x_3}^lF_i^{\pm}) ] \\ &\qquad\qquad+\sum_{\substack{i+j=k+2\\i,j\geq1}}\frac{1}{2} [ Q(\bar{F}_{i}^{\pm},\bar{F}_{j}^{\pm}) + Q(\bar{F}_{j}^{\pm},\bar{F}_{i}^{\pm}) ],
	\end{split}
\end{align}
where we have used the notations $\bar{F}_i\equiv0$ for $i\leq 0$. To avoid the complicate interplay of interior layers and viscous boundary layers,  we have used  the Taylor expansions of   $F_i$ at $x_3=0,1$  for $i\geq 0$, i.e.,
\begin{equation}\label{1.14-2}
\begin{aligned}
	F_i(t,x_{\sp},x_3,v)
	&=F_i^{-}+\sum_{l=1}^{\fb} \frac{1}{l!} \partial_{x_3}^l F_{i}^{-}\cdot x_3^l + \frac{x_3^{\fb+1}}{(\fb+1)!} \partial_{x_3}^{\fb+1}\mathfrak{F}_{i}^{-},\\
	F_i(t,x_{\sp}, x_3,v)
	&=F_i^{+}+\sum_{l=1}^{\fb} \frac{1}{l!} \partial_{x_3}^l F_i^{+}\cdot (x_{3}-1)^l + \frac{(x_{3}-1)^{\fb+1}}{(\fb+1)!} \partial_{x_3}^{\fb+1}\mathfrak{F}_{i}^{+},
\end{aligned}
\end{equation}
where we have used the simplified notations
\begin{align}\label{1.14-3}
	\begin{split}
		\partial_{x_3}^lF_i^{-}:&=(\partial_{x_3}^lF_i)(t,x_{\sp}, 0, v),
		\quad \partial_{x_3}^{\fb+1}\mathfrak{F}_{i}^{-}:=(\partial_{x_3}^{\fb+1}F_i)(t,x_{\sp}, \xi_{i}^{-},v),\\
		\partial_{x_3}^lF_i^{+}:&=(\partial_{x_3}^lF_i)(t,x_{\sp}, 1, v),
		\quad \partial_{x_3}^{\fb+1}\mathfrak{F}_{i}^{+}:=(\partial_{x_3}^{\fb+1}F_i)(t,x_{\sp}, \xi_{i}^{+},v),
	\end{split}
\end{align}
for $0\leq l\leq \fb$ and some $\xi_{i}^{-} \in (0,x_3)$ and $\xi_{i}^{+}\in (x_{3},1)$ with $i\geq0$.  The number $\fb\in \mathbb{N}_+$ will be chosen later. As in \cite{Guo-Huang-Wang}, the main reason to use \eqref{1.14-2} is to make the coefficients of viscous boundary layer systems \eqref{1.14} be independent of $\v$.

\subsubsection{ Knudsen boundary layer expansion:}
To construct the solution that satisfies the boundary condition at higher orders, we have to introduce the Knudsen boundary layer. Firstly, we define the new scaled normal coordinates:
\begin{equation}\nonumber%
	\eta_{-}:=\frac{x_3}{\v^n}\qquad \eta_{+}:=\frac{1-x_3}{\v^n}.
\end{equation}
The Knudsen boundary layer expansion is defined as
\begin{equation}\nonumber
	\hat{F}^{\pm,\v}(t,x_\sp, \eta_{\pm},v)\sim \v^{n-2}\sum_{i=1}^\infty \v^i \hat{F}_{i}^{\pm}(t,x_\sp, \eta_{\pm},v).
\end{equation}
Similarly, we always use $\eta$ to represent $\eta_{+}$ and $\eta_{-}$ for simplicity of presentations.

Plugging $F^\v+\bar{F}^{\pm,\v}+\hat{F}^{\pm,\v}$  into \eqref{1.1} and comparing the order of $\v$,  then using \eqref{1.7-1} and \eqref{1.14}, one obtains
\begin{equation}\label{1.19-1}
	\begin{split}
		\dis& \v^{-2+i} \, (1\leq i\leq n-2):  \quad \mp v_3 \frac{\partial\hat{F}_{i}^{\pm}}{\partial \eta}-[Q(\mu,\hat{F}_{i}^{\pm})+Q(\hat{F}_{i}^{\pm},\mu)]=0,\\[2mm]
		\dis& \v^{n-3}:  \quad \mp v_3 \frac{\partial\hat{F}_{n-1}^{\pm}}{\partial \eta }-[Q(\mu,\hat{F}_{n-1}^{\pm})+Q(\hat{F}_{n-1}^{\pm},\mu)]=Q(F_0^{\pm},\hat{F}_{1}^{\pm}) + Q(\hat{F}_{1}^{\pm}, F_0^{\pm}),\\[2mm]
		\dis & \v^{n-2}:\quad \mp v_3 \frac{\partial\hat{F}_{n}^{\pm}}{\partial \eta }-[ Q(\mu,\hat{F}_{n}^{\pm})+Q(\hat{F}_{n}^{\pm},\mu) ]\\
		&\qquad\quad= \sum_{\substack{i+j=2\\  i\geq 0, j\geq1}} [ Q(F_i^{\pm}, \hat{F}_{j}^{\pm}) + Q(\hat{F}_{j}^{\pm}, F_i^{\pm}) ] + [ Q(\bar{F}_1^{\pm,0}, \hat{F}_{1}^{\pm}) + Q(\hat{F}_{1,\pm}, \bar{F}_1^{\pm,0})],\\
		&\qquad\qquad + Q(\hat{F}_{1}^{\pm}, \hat{F}_{1}^{\pm})\\[2mm]
		&\dis \quad\quad\quad\quad\quad\quad\quad\ \vdots \\
		\dis & \v^{n-2+k}(k\geq 1):\quad   \mp v_3\frac{\partial\hat{F}_{k+n}^{\pm}}{\partial \eta }-[ Q(\mu,\hat{F}_{k+n}^{\pm})+Q(\hat{F}_{k+n}^{\pm},\mu) ] \\
		&\qquad\quad=-\partial_t \hat{F}_{k+2-n}^{\pm} -v_\sp\cdot\nabla_\sp\hat{F}_{k}^{\pm}
		+\sum_{\substack{i+j=k+2\\ i\geq 0,j\geq1}}  [Q(F_i^{\pm}, \hat{F}_{j}^{\pm})+Q( \hat{F}_{j}^{\pm}, F_i^{\pm}) ] \\
		&\qquad\qquad+\sum_{\substack{i+j=k+2\\ i,j\geq1}} \big\{[Q(\bar{F}_{i}^{\pm,0},\hat{F}_{j}^{\pm})+Q(\hat{F}_{j}^{\pm},\bar{F}_{i}^{\pm,0})]  +\f12[Q(\hat{F}_{i}^{\pm},\hat{F}_{j}^{\pm}) + Q(\hat{F}_{j}^{\pm}, \hat{F}_{i}^{\pm})] \big\}\\
		&\qquad\qquad+\sum_{\substack{i+j+nl=k+2\\ i\geq0,j\geq1, 1\leq l\leq \fb}} \frac{(\mp \eta)^l}{l!} [Q(\partial_{x_3}^l F_i^{\pm}, \hat{F}_{j}^{\pm})+Q( \hat{F}_{j}^{\pm},\partial_{x_3}^l F_i^{\pm}) ]   \\
		&\qquad\qquad+\sum_{\substack{i+j+(n-1)l=k+2\\ i,j\geq1, 1\leq l\leq \fb}}\frac{\eta^l}{l!} [Q(\partial_{y}^l \bar{F}_{i}^{\pm,0}, \hat{F}_{j}^{\pm})+Q( \hat{F}_{j}^{\pm},\partial_{y}^l \bar{F}_{i}^{\pm,0}) ] ,
	\end{split}
\end{equation}
where we have used \eqref{1.14-2} and  the Taylor expansion of $\bar{F}_{i}^{\pm}$
\begin{equation*}
	\bar{F}_{i}^{\pm}(t, x_{\sp}, y,v)
	=\bar{F}_{i}^{\pm,0}+\sum_{l=1}^{\fb} \frac{y^l}{l!} \partial_{y}^l \bar{F}_i^{\pm,0} + \frac{y^{\fb+1}}{(\fb+1)!} \partial_{y}^{\fb+1}\bar{\mathfrak{F}}_{i}^{\pm},
\end{equation*}
with
\begin{equation}\label{1.14-4}
	\partial_{y}^l\bar{F}_{i}^{\pm,0}:=(\partial_{y}^l\bar{F}_{i}^{\pm})(t,x_{\sp}, 0, v),
	\quad \partial_{y}^{\fb+1}\bar{\mathfrak{F}}_{i}^{\pm}:=(\partial_{y}^{\fb+1}\bar{F}_{i}^{\pm})(t,x_{\sp}, \bar{\xi}_{i}^{\pm},v),
	\, \mbox{for}\,\, 0\leq l\leq \fb,
\end{equation}
for some $\bar{\xi}_{i}^{\pm}\in [0,y_{\pm}]$ with $i\geq 1$. It is noted that the Knudsen boundary layer \eqref{1.19-1} is in fact a steady problem with $(t,x_\sp)$ as parameters, and the well-posedness has already been established in \cite{GPS-1988,Jiang-Wang} under some conditions on the source term and boundary conditions.

\subsection{Hilbert expansion}
Let $\phi(z)$ be a smooth cutoff function satisfying
\begin{equation}\label{1.22-1}
	\phi(z)=\left\{
	\begin{aligned}
		1,&\quad 0\leq z\leq \frac{1}{4},\\
		0,&\quad \frac{1}{2}\leq z\leq 1.
	\end{aligned}
	\right.
\end{equation}
Now, we consider the scaled Boltzmann solution with the following truncated Hilbert expansion with multi-scales
\begin{equation}\nonumber
	\begin{aligned}
		\mathscr{F}^\v=\mu+ \v^{n-2}\Big\{&\sum_{i=0}^{N} \v^i F_i(t,x,v)+\sum\limits_{\pm}\sum_{i=1}^{N} \v^i\phi(\v y)\bar{F}_{i}^{\pm}(t,x_\sp,y,v)\\
		&+\sum\limits_{\pm}\sum_{i=1}^{N} \v^i\phi(\v^{n}\eta)\hat{F}_{i}^{\pm}(t,x_\sp,\eta,v)\Big\} +\v^{n-2+k_0} F^\v_{R},
	\end{aligned}
\end{equation}
where $k_{0}, N$ are some positive parameters determined later.

It follows from \eqref{1.8-1}, \eqref{1.14} and \eqref{1.19-1} that the remainder $F_{R}^{\v}$ satisfies
\begin{equation}\label{1.23}
	\begin{aligned}
		&\partial_{t}F_{R}^{\v}+\frac{1}{\v^{n-2}}v\cdot \nabla_{x}F_{R}^{\v}-\frac{1}{\v^{2n-2}}\big\{[Q(\mu,F_{R}^{\v})+Q(F_{R}^{\v},\mu)]+\v^{n-2}[Q(F_{0},F_{R}^{\v})+Q(F_{R}^{\v},F_{0})]\big\}\\
		&=\sum\limits_{\pm}\sum\limits_{i=1}^N\frac{1}{\v^{n-i}}[Q(F_{R}^{\v},F_{i}+[\phi(\v y)\bar{F}_{i}^{\pm}+\phi(\v^{n}\eta)\hat{F}_{i}^{\pm}])+Q(F_{i}+[\phi(\v y)\bar{F}_{i}^{\pm}+\phi(\v^{n}\eta)\hat{F}_{i}^{\pm}],F_{R}^{\v})]\\
		&\quad +\v^{k_0-n}Q(F_{R}^{\v},F_{R}^{\v})+\frac{1}{\v^{2n-4+k_0}}(R^{\v}+\bar{R}^{+,\v}+\bar{R}^{-,\v}+\hat{R}^{+,\v}+\hat{R}^{-,\v}),
	\end{aligned}
\end{equation}
where
\begin{align}
	&\begin{aligned}
		R^{\v}=&-\sum\limits_{i=N-2n+3}^{N}\v^{i+2n-4}\partial_{t}F_{i}-\sum\limits_{i=N-n+1}^{N}\v^{i+n-2}v\cdot \nabla_{x}F_{i}\\
		&+\sum\limits_{i+j\geq N-n+3\atop 0\leq i,j\leq N}\frac{1}{2}\v^{i+j+n-4}[Q(F_{i},F_{j})+Q(F_{j},F_{i})],
	\end{aligned}\label{1.24}\\
&\begin{aligned}
		\bar{R}^{\pm,\v}=&\pm v_{3}\phi'(\v y)\sum\limits_{i=1}^{N}\v^{n-2+i}\bar{F}_{i}^{\pm}-\sum\limits_{i=N-2n+3}^{N}\v^{i+2n-4}\phi(\v y)\partial_{t}\bar{F}_{i}^{\pm}\\
		&-\sum\limits_{i=N-n+1}^{N}\v^{i+n-2}\phi(\v y)v_{\sp}\cdot \nabla_{\sp}\bar{F}_{i}^{\pm}\pm \sum\limits_{i=N-n+2}\v^{i+n-3}\phi(\v y)v_3\partial_{y}\bar{F}_{i}^{\pm}\\
		&+\sum\limits_{i+j\geq N-n+3\atop 0\leq i\leq N,1\leq j\leq N}\phi(\v y)\v^{i+j+n-4}[Q(F_{i}^{\pm},\bar{F}_{j}^{\pm})+Q(\bar{F}_{j}^{\pm},F_{i}^{\pm})]\\
		&+\sum\limits_{i+j+l\geq N-n+3\atop 0\leq i\leq N, 1\leq j\leq N,1\leq l\leq \fb}\v^{i+j+l+n-4}\frac{(\mp y)^{l}}{l!}\phi(\v y)[Q(\partial_{x_3}^{l}F_{i}^{\pm},\bar{F}_{j}^{\pm})+Q(\bar{F}_{j}^{\pm},\partial_{x_3}^{l}F_{i}^{\pm})]\\
		&+\sum\limits_{0\leq i\leq N,1\leq j\leq N}\v^{i+j+\fb+n-3}\frac{(\mp y)^{\fb+1}}{(\fb+1)!}\phi(\v y)[Q(\partial_{x_3}^{\fb+1}\mathfrak{F}_{i}^{\pm},\bar{F}_{j}^{\pm})+Q(\bar{F}_{j}^{\pm},\partial_{x_3}^{\fb+1}\mathfrak{F}_{i}^{\pm})],
	\end{aligned}\label{1.26}
\end{align}
and
\begin{equation}\label{1.28}
	\begin{aligned}
		\hat{R}^{\pm,\v}=
		&\pm v_3\phi'(\v^{n}\eta)\sum\limits_{i=1}^{N}\v^{i+n-2}\hat{F}_{i}^{\pm}\\
		&-\sum\limits_{i=N-2n+3}^{N}\v^{i+2n-4}\phi(\v^{n}\eta)\partial_{t}\hat{F}_{i}^{\pm}-\sum\limits_{i=N-n+1}^{N}\phi(\v^{n}\eta)\v^{n-2+i}v_{\sp}\cdot \nabla_{\sp}\hat{F}_{i}^{\pm}\\
		&+\sum\limits_{i+j\geq N-n+3\atop  0\leq i\leq N,1\leq j\leq N}\phi(\v^{n}\eta)\v^{i+j+n-4}[Q(F_{i}^{\pm},\hat{F}_{j}^{\pm})+Q(\hat{F}_{j}^{\pm},F_{i}^{\pm})]\\
		&+\sum\limits_{i+j\geq N-n+3\atop 1\leq i\leq N,1\leq j\leq N}\v^{i+j+n-4}\phi^2(\v^{n}\eta)[Q(\bar{F}_{i}^{\pm,0},\hat{F}_{j}^{\pm})+Q(\hat{F}_{j}^{\pm},\bar{F}_{i}^{\pm,0})]\\
		&+\sum\limits_{i+j+nl\geq N-n+3\atop 0\leq i\leq N,1\leq j\leq N,1\leq l\leq \fb}\frac{(\mp\eta)^{l}}{l!}\v^{i+j+nl+n-4}\phi(\v^{n}\eta)[Q(\partial_{x_3}^{l}F_{i}^{\pm},\hat{F}_{j}^{\pm})+Q(\hat{F}_{j}^{\pm},\partial_{x_3}^{l}F_{i}^{\pm})]\\
		&+\sum\limits_{i+j+(n-1)l\geq N-n+3\atop 0\leq i\leq N,1\leq j\leq N,1\leq l\leq \fb}\frac{\eta^{l}}{l!}\v^{i+j+(n-1)l+n-4}\phi^2(\v^{n}\eta)[Q(\partial_{y}^{l}\bar{F}_{i}^{\pm,0},\hat{F}_{j}^{\pm})+Q(\hat{F}_{j}^{\pm},\partial_{y}^{l}\bar{F}_{i}^{\pm,0})]\\
		&+\sum\limits_{0\leq i\leq N,1\leq j\leq N}\frac{(\mp\eta)^{\fb+1}}{(\fb+1)!}\v^{i+j+n(\fb+1)+n-4}\phi(\v^{n}\eta)[Q(\partial_{x_3}^{\fb+1}\mathfrak{F}_{i}^{\pm},\hat{F}_{j}^{\pm})+Q(\hat{F}_{j}^{\pm},\partial_{x_3}^{\fb+1}\mathfrak{F}_{i}^{\pm})]\\
		&+\sum\limits_{1\leq i,j\leq N}\frac{\eta^{\fb+1}}{(\fb+1)!}\v^{i+j+(n-1)(\fb+1)+n-4}\phi^2(\v^n\eta)[Q(\partial_{y}^{\fb+1}\bar{\mathfrak{F}}_{i}^{\pm},\hat{F}_{j}^{\pm})+Q(\hat{F}_{j}^{\pm},\partial_{y}^{\fb+1}\bar{\mathfrak{F}}_{i}^{\pm})],
	\end{aligned}
\end{equation}
where $\partial_{x_3}^{l}F_{i}^{\pm}, \partial_{y}^{l}\bar{F}_{i}^{\pm,0}, \partial_{x_3}^{\fb+1}\mathfrak{F}_{i}^{\pm}$ and $\partial_{y}^{\fb+1}\bar{\mathfrak{F}}_{i}^{\pm}$ are the ones defined in \eqref{1.14-3} and \eqref{1.14-4}.

To apply $L^2-L^{\infty}$ framework for \eqref{1.23}, we introduce following local Maxwellian
\begin{equation}\label{1.30}
	\mathcal{M}_{\v}:=\frac{1+\v^{n-2}\rho_{0}}{(2\pi(1+\v^{n-2}\theta_{0}))^{\frac{3}{2}}}\exp\left\{-\frac{|v-\v^{n-2}u_{0}|^2}{2(1+\v^{n-2}\theta_0)}\right\},
\end{equation}
where $(\rho_{0}, u_{0}, \theta_0)$ are the macro quantities of $F_{0}$, and satisfy the incompressible Euler system \eqref{2.8-1}. We point out that we have to use the local Maxwellian $\mathcal{M}_{\v}$ so that the low decay term $\displaystyle\frac{1}{\v^{n}}\big[Q(F_{0},F_{R}^{\v})+Q(F_{R}^{\v},F_{0})\big]$ can be absorbed into the linearized operator. Otherwise, the decay of remainder is not enough to close the estimate. Taking the Taylor expansion of $\mathcal{M}_{\v}$ with respect to $\v$, we get
\begin{equation}\label{1.31}
	\mathcal{M}_{\v}=\mu+\v^{n-2} F_{0}+\v^{2n-4}r_0,
\end{equation}
where
$$
\begin{aligned}
r_{0}=\frac{1}{2}\big[&\big(\frac{\rho_{0}}{1+\tilde{\v}^{n-2}\rho_{0}}+\frac{(v-\tilde{\v}^{n-2} u)\cdot u}{1+\tilde{\v}^{n-2}\theta_{0}}+\frac{|v-\tilde{\v}^{n-2}u_{0}|^2\theta_{0}}{2(1+\tilde{\v}^{n-2}\theta_{0})^2}-\frac{3\theta_{0}}{2(1+\tilde{\v}^{n-2}\theta_{0})}\big)^2-\frac{\rho_{0}^2}{(1+\tilde{\v}^{n-2}\rho_{0})^2}\\
&-\frac{|u_{0}|^2}{1+\tilde{\v}^{n-2}\theta_{0}}+\frac{3\theta_{0}^2}{2(1+\tilde{\v}^{n-2}\theta_{0})^2}-\frac{2(v-\tilde{\v}^{n-2}u_{0})\cdot u_{0}}{(1+\tilde{\v}^{n-2}\theta_{0})^2}\theta_{0}-\frac{|v-\tilde{\v}^{n-2}u_{0}|^2\theta_{0}^2}{(1+\tilde{\v}^{n-2}\theta_{0})^3}\big] \mathcal{M}_{\tilde{\v}},
\end{aligned}
$$
for some $\tilde{\v}\in (0,\v)$.

We rewrite the remainder as
\begin{equation}\label{1.32}
	F_{R}^{\v}=\sqrt{\mathcal{M}_{\v}}f_{R}^{\v}.
\end{equation}
On the other hand, we introduce a global Maxwellian
\begin{equation}\nonumber
	\mu_{M}:=\frac{1}{(2\pi \theta_{M})^{\frac{3}{2}}}\exp\left\{-\frac{|v|^2}{2\theta_{M}}\right\},
\end{equation}
where $\theta_{M}>0$ satisfies the condition
\begin{equation}\label{1.34}
	\theta_{M}<\min_{(t,x)\in [0,\tau]\times \Omega}(1+\v^{n-2}\theta_0(t,x))\leq \max_{(t,x)\in [0,\tau]\times \Omega}(1+\v^{n-2}\theta_0(t,x))<2\theta_{M}.
\end{equation}
By the assumption \eqref{1.34}, one can easily deduce that there exists positive constant $C>0$ such that
\begin{equation}\label{1.35}
	\frac{1}{C}\mu_{M}\leq \mathcal{M}_{\v}(t,x,v)\leq C\mu_{M}^{\alpha},\qquad \text{for some }\frac{1}{2}<\alpha<1.
\end{equation}
We further define
\begin{equation}\label{1.36}
	F_{R}^{\v}=(1+|v|^2)^{-\frac{\beta}{2}}\sqrt{\mu_{M}}h_{R}^{\v}=\frac{1}{\omega_{\beta}(v)}\sqrt{\mu_{M}}h_{R}^{\v},
\end{equation}
with the velocity weight function
\begin{equation}\nonumber
	\omega_{\beta}(v):=(1+|v|^2)^{\frac{\beta}{2}},\quad \text{for }\beta\geq 3.
\end{equation}

\begin{theorem}\label{theorem1.1}
	Let $\tau>0$ be the lifespan of smooth solution of incompressible Euler equations given in Lemma \ref{lem1.1}. For fixed $n\geq 3$, let $k_0>3n-2$, $\beta\geq \frac{9}{2}+\frac{k_0-1}{n-1}$, $N\geq 2n+k_0-3$ and $\fb\geq n+k_0-\frac{5}{2}$.  We assume the initial data
    \begin{equation}\label{1.37-2}
	\begin{aligned}
		F^{\v}(0,x,v)=\mu+ \v^{n-2} \Big\{&\sum_{i=0}^{N}\v^i  F_i(0,x,v)+\sum\limits_{\pm}\sum_{i=1}^{N}\v^i\phi(\v y) \bar{F}_{i}^{\pm}(0,x_{\sp},y,v)\\
		&+\sum\limits_{\pm}\sum_{i=1}^{N} \v^{i}\phi(\v^n \eta)\hat{F}_{i}^{\pm}(0,x_\sp,\eta,v) \Big\}+\v^{k_0+n-2}F^\v_R(0,x,v)\geq0,
	\end{aligned}
	\end{equation}
	and $F_i(0), \bar{F}_{j}^{\pm}(0), \hat{F}_{j}^{\pm}(0),  i=0,\cdots, N$, $j=1,\cdots, N$  satisfy the regularity and compatibility conditions described in Proposition \ref{prop5.1}, and
	\begin{equation}\label{1.37-1}
	\left\|\left(\frac{F_{R}^{\v}}{\sqrt{\mathcal{M}_{\v}}}\right)(0,x,v)\right\|_{L_{x,v}^2}+\v^{3n-3}\left\|\omega_{\beta}(v)\left(\frac{F_{R}^{\v}}{\sqrt{\mu_{M}}}\right)(0,x,v)\right\|_{L_{x,v}^{\infty}}<\infty.
	\end{equation}
	Then there exists a small positive constants $\v_0>0$ such that IBVP of Boltzmann equation \eqref{1.3}-\eqref{1.1}, \eqref{1.37-2} has a unique solution for $\v\in (0,\v_0]$ over the interval $t\in [0,\tau]$ in the following form of expansion
	\begin{equation}\label{1.22-2}
		\begin{aligned}
		\mathscr{F}^\v=\mu+ \v^{n-2}\Big\{&\sum_{i=0}^{N} \v^i F_i(t,x,v)+\sum\limits_{\pm}\sum_{i=1}^{N} \v^i\phi(\v y)\bar{F}_{i}^{\pm}(t,x_\sp,y,v)\\
		&+\sum\limits_{\pm}\sum_{i=1}^{N} \v^i\phi(\v^{n}\eta)\hat{F}_i(t,x_\sp,\eta,v)\Big\} +\v^{k_0+n-2} F^\v_{R}\geq 0
		\end{aligned}
	\end{equation}
	with
	\begin{equation}\label{1.22-3}
		\sup_{t\in [0,\tau]}\left\{\left\|\frac{F_{R}^{\v}(t)}{\sqrt{\mathcal{M}_{\v}}}\right\|_{L_{x,v}^2}+\v^{3n-3}\left\|\omega_{\beta}(v)\frac{F_{R}^{\v}(t)}{\sqrt{\mu_{M}}}\right\|_{L_{x,v}^{\infty}}\right\}\leq C(\tau)<\infty.
	\end{equation}
	Here the functions $F_{i}(t,x,v), \bar{F}_{i}^{\pm}(t,x_{\sp},y,v)$ and $\hat{F}_{i}^{\pm}(t,x_{\sp},\eta,v)$ are respectively the interior expansion, viscous and Knudsen boundary layers constructed in Proposition \ref{prop5.1}.
\end{theorem}

\begin{remark}
	From \eqref{1.22-2}-\eqref{1.22-3} and the uniform estimates in Proposition \ref{prop5.1}, it is clear that
	$$
	\sup_{t\in [0,\tau]}\Big\{\Big\|\frac{1}{\sqrt{\mathcal{M}_{\v}}}\Big[\v^{2-n}(\mathscr{F}^{\v}-\mu)-\big(\rho_0(t,x)+u_{0}(t,x)\cdot v+\frac{\theta_0(t,x)}{2}(|v|^2-3)\big)\mu\Big]\Big\|_{L_{x,v}^2}\Big\}\leq C\v\rightarrow 0,
	$$
	and
	$$
	\sup_{t\in [0,\tau]}\Big\{\Big\|\frac{\omega_{\beta}(v)}{\sqrt{\mu_{M}}}\Big[\v^{2-n}(\mathscr{F}^{\v}-\mu)-\big(\rho_0(t,x)+u_{0}(t,x)\cdot v+\frac{\theta_0(t,x)}{2}(|v|^2-3)\big)\mu\Big]\Big\|_{L_{x,v}^{\infty}}\Big\}\leq C\v\rightarrow 0,
	$$
	where $(\rho_{0}(t,x),u_{0}(t,x),\theta_0(t,x))$ is the unique local smooth solution of the incompressible Euler equation \eqref{2.8-1}-\eqref{2.8-3}. Therefore we have established the hydrodynamic limit from the Boltzmann equation to the incompressible Euler system in the channel $\mathbb{T}^2\times (0,1)$.
\end{remark}


\begin{remark}
Compared to the Hilbert expansion of compressible Euler scaling \cite{Guo-Huang-Wang}, due to the influence of $\mathscr{M}_{a}\to 0$, the order of expansion should be more higher, and the expansions are more complicated in the case of incompressible Euler scaling.
\end{remark}

Now, we briefly mention some key analysis of the present paper. Recalling \eqref{2.8-1}, the macroscopic quantities of $f_{0}$ satisfies the incompressible Euler equations, while the macroscopic quantities of $f_{k}(k\geq 1)$ satisfies a kind of linearized Euler system with a non-zero divergence:
\begin{equation}\label{1.41}
\operatorname{div}_{x}u_{k}=-\pa_{t}\rho_{k+2-n},
\end{equation}
see Lemma \ref{lem2.1} below for details.  Integrating \eqref{1.41} to see the boundary conditions should satisfy
\begin{equation}\label{cc}
\int_{\mathbb{T}^2}(u_{k,3}(t,x_{\sp},1)-u_{k,3}(t,x_{\sp},0))dx_{\sp}
=\int_{\Omega}(-\pa_{t}\rho_{k+2-n})dx,\quad 0\leq k\leq N.
\end{equation}
In fact, using \eqref{2.17-1} and \eqref{4.33}-\eqref{4.34} below, then \eqref{cc} leads to the restriction on the $\int_{\Omega}p_{k+2-n}(t,x)dx$. These restrictions
will help us to determine uniquely the pressure terms $p_k$ for $0\leq k\leq N$, see Lemma \ref{lem2.2} and Remarks \ref{rem2.3}-\ref{rem2.3-1} for details.

For the uniform estimates of the $u_{k}$ in \eqref{2.17}, since the compatibility conditions can be guaranteed by imposing some restrictions on $\int_{\Omega}p_{k+2-n}(t,x)$, we can find an auxiliary function $\hat{u}_{k}$ such that $\operatorname{div}_{x}(u_{k}-\hat{u}_{k})=0$ and $\hat{u}_{k,3}\vert_{\pa \Omega}=u_{k,3}\vert_{\pa \Omega}$. In fact, we consider $\hat{u}_{k}=\nabla_{x} q_{k}$, then $\Delta_{x} q_{k}=\operatorname{div}u_{k}$ with $\nabla_{x} q_{k}\cdot \vec{n}\vert_{\partial \Omega}=u_{k,3}\vert_{\pa \Omega}$, and the uniform estimates of $q_{k}$ can be obtained from classical theory of elliptic equation. Moreover, $\tilde{u}_{k}:=u_{k}-\hat{u}_{k}$ satisfies the usual linearized incompressible Euler equations and homogeneous boundary condition: $\tilde{u}_{k,3}\vert_{\pa_{\Omega}}=0$. Then applying the standard energy method to $\tilde{u}_{k}$, we succeed in obtaining $\|\tilde{u}_{k}\|_{H^{k}}$ and hence $\|u_{k}\|_{H^{k}}$.
We note that, in some known results \cite{Masi-Esposito-Lebowitz-1989,Wu-Zhou-Li-2019}, they usually skip the construction of the linear parts of Hilbert expansion in the case of Cauchy (or Torus)  problem. Here we point out that one should be very careful in the case of initial boundary value problem due to the boundary effect.

For the viscous boundary layer, we deduce the corresponding equations into a linear parabolic system with partial
viscosity and linear growth coefficients, and the solutions of such systems have already been established in \cite{Guo-Huang-Wang}. However, the source terms and boundary conditions are more complicated than those in \cite{Guo-Huang-Wang}, and we need to analyze more carefully. In fact, for the case of $n=3$, three terms $\la \mathcal{A}_{3i}, J_{k-1}^{\pm}\ra(i=1,2), \la \mathcal{B}_{3}, J_{k-1}^{\pm}\ra $ in the source terms and boundary conditions for the equations of the $k-$th viscous boundary layer, which cannot be generally calculated explicitly, and depend on the macroscopic velocity and temperature of $\bar{f}_{k}^{\pm}$ from a first glance, see Proposition \ref{prop3.1} and Lemma \ref{lem4.3}. Fortunately, due to $u_{0,3}\vert_{\pa \Omega}=\bar{u}_{1,3}^{\pm}\equiv 0$ and the property that $\mathbf{L}^{-1}$ preserves the odevity of velocity, we can prove that the highest order terms in those terms depend only on $\bar{u}_{k,3}^{\pm}$, which is already determined by $\bar{f}_{i}^{\pm}(1\leq i\leq k-1)$, see Remark \ref{rem2.1} and Appendix \ref{AppendixD}. This observation helps us to construct $\bar{f}_{k}^{\pm}$ inductively by using similar arguments as in \cite{Guo-Huang-Wang}.

For the weighted $L^{\infty}-$estimate of remainder term, since the characteristic lines of \eqref{1.1} depend on $\varepsilon$, the collision number $k$ for each backward bi-characteristic depends on $\v$, and $k\to \infty$ as $\v\to 0$. This is different from previous works, such as \cite{Cao-Jang-Kim,Guo2010,Guo-Huang-Wang,Jiang-Wang}. Thus we need to analyze carefully such that all estimates are independent of $\v$, especially after taking Vidav's iteration.

The rest of this paper is organized as follows. In \S 2, we reformulate the interior expansions \eqref{1.7-1}, viscous boundary layers \eqref{1.14} and Knudsen boundary layers \eqref{1.19-1}. The boundary conditions for interior expansions and viscous boundary layers, and criteria of uniqueness for the pressure terms will be given in \S 3.
In \S 4, the existence theory and uniform energy estimates for a kind of linearized Euler equations with nonvanishing divergence will be presented. Then  we briefly state an existence theory of IBVP for a kind of parabolic systems.
In \S 5, we present the complete procedures of constructions of interior expansions, the viscous and Knudsen boundary layers.
Finally, in \S 6,  by using $L^2-L^{\infty}$ framework, we prove Theorem \ref{theorem1.1}. Detailed calculations on the derivations of macroscopic equations for interiors expansions and viscous layers are given in Appendices \ref{AppendixA} and \ref{AppendixD}, respectively.

\vspace{1.5mm}
\noindent{\bf Notations.}  Throughout this paper, $C$ denotes a generic positive constant  and may vary from line to line. And $C(a),C(b),\cdots$ denote the generic positive constants depending on $a,~b,\cdots$, respectively, which also may vary from line to line. We use $\langle\cdot ,\cdot \rangle$ to denote the standard $L^2$ inner product in $\R^3_v$.
$\|\cdot\|_{L^2}$ denotes the standard $L^2(\Omega\times\mathbb{R}^3_v)$-norm, and $\|\cdot\|_{L^\infty}$ denotes the $L^\infty(\Omega\times\mathbb{R}^3_v)$-norm. \vspace{1.5mm}

\section{Reformulation of Hilbert expansions and boundary conditions}
\subsection{Interior expansion}
It follows from $\eqref{2.1}_1$ that
$$
f_{0}=\mathbf{P}f_{0}=\{\rho_{0}+v\cdot u_{0}+\frac{1}{2}\theta_{0}(|v|^2-3)\}\sqrt{\mu},
$$
where $(\rho_0,u_{0},\theta_0)$ satisfy the incompressible Euler equations \eqref{2.8-1}. For system \eqref{2.8-1}-\eqref{2.8-4}, we have the following existence result and uniform estimate.

\begin{lemma}\label{lem1.1}
	Assume the initial data \eqref{2.8-3} satisfies the compatibility conditions (the compatibility conditions means the initial data satisfies \eqref{2.4} and \eqref{2.8-3}, and the time-derivatives of initial data are defined through \eqref{2.17} inductively). Let $s_{0}\geq 3$, then there exist a local time $\tau>0$, which is proportional to the inverse of $\|u_{0}(0)\|_{H^{s_0}}$, and a unique local smooth solution $(\rho_{0},u_{0},\theta_{0},p_{n-2})$ of system \eqref{2.8-1}-\eqref{2.8-4} such that
	\begin{equation*}
		(\rho_{0},u_{0},\theta_{0})\in L^{\infty}(0,\tau; H^{s_0}(\Omega)),\quad \nabla p_{n-2}\in L^{\infty}(0,\tau; H^{s_0}(\Omega)),
	\end{equation*}
Moreover, it holds that
	\begin{equation*}
		\sup_{t\in [0,\tau]}\|(\rho_{0},u_{0},\theta_{0})(t)\|_{H^{s_{0}}(\Omega)}+\sup_{t\in [0,\tau]}\|\nabla_{x} p_{n-2}(t)\|_{H^{s_{0}}(\Omega)}\leq C\|(\rho_{0}, u_{0},\theta_{0})(0)\|_{H^{s_{0}}(\Omega)},
	\end{equation*}
In addition, by the classical Poincar\'{e} inequality, one has
$$
\begin{aligned}
\sup_{t\in [0,\tau]}\|p_{n-2}(t)\|_{L^2(\Omega)}&\leq \sup_{t\in [0,\tau]}\|\nabla_{x} p_{n-2}\|_{L^2(\Omega)}+\|\mathfrak{c}_{n-2}(t)\|_{L^{\infty}}\\
&\leq C\|(\rho_{0}, u_{0},\theta_{0})(0)\|_{H^{s_{0}}(\Omega)}+\|\mathfrak{c}_{n-2}(t)\|_{L^{\infty}}.
\end{aligned}
$$
\end{lemma}
The proof of Lemma \ref{lem1.1} is classical, for instance, one can see \cite[Theorem 4.1]{Petcu}, so we omit it here for brevity.

The other interior layers, i.e., $F_{1}(t,x,v),\cdots,F_{N}(t,x,v)$ are constructed inductively in following lemma, whose proof will be presented in Appendix \ref{AppendixA}.
\begin{lemma}[Reformulation of interior expansion]\label{lem2.1}
	Let $F_{i}=\sqrt{\mu}f_{i}$ ($i\geq 1$) be the solution of \eqref{1.7-1}, then the microscopic part of $f_{i}$ is given by
	\begin{equation}\label{2.0}
		\begin{aligned}
			&(\mathbf{I}-\mathbf{P}) f_{i} =0, \quad 0 \leq i \leq n-3 ; \\
			&\begin{aligned}
			(\mathbf{I}-\mathbf{P}) f_{k+n-2} =\mathbf{L}^{-1}\big\{&(\mathbf{I}-\mathbf{P})-\partial_{t} f_{k-n}-v \cdot \nabla_{x} f_{k-2})\\
			&+\sum_{i+j=k \atop i \geq 0, j \geq 0} \frac{1}{2} [\Gamma\left(f_{i}, f_{j}\right)+\Gamma\left(f_{j}, f_{i}\right)]\big\}\quad k\geq 0,
			\end{aligned}
		\end{aligned}
	\end{equation}
	where we have used the notation $f_{i}\equiv 0$ for $i\leq -1$, that is, the microscopic part of $f_{i}$ is determined by $f_{j}$ with $1\leq j\leq i-(n-2)$.
	
For the macroscopic part $(\r_{i},u_{i},\t_{i})$ of $f_{i}$, it holds
	\begin{equation}\label{2.17}
		\left\{\begin{aligned}
			&\operatorname{div}_{x} u_{m}=-\partial_{t}\rho_{m+2-n}, \\
			&\partial_{t} u_{m}+(u_{0} \cdot \nabla_{x}) u_{m}+(u_{m} \cdot \nabla_{x}) u_{0}+\nabla_{x} p_{n+m-2}=\mathfrak{h}_{m-1},\\
			&\partial_{t} \theta_{m}+(u_{0} \cdot \nabla_{x}) \theta_{m}+(u_{m} \cdot \nabla_{x}) \theta_{0}=\mathfrak{q}_{m-1},\\
		\end{aligned}\right.\qquad\text{for }m\geq 1,
	\end{equation}
	where we have denoted
	\begin{equation}\label{2.17-1}
		\begin{aligned}
		&p_{k}:=\rho_{k}+\theta_{k}\equiv p_{k}(t),\quad \text{for }0\leq k\leq n-3,\\
		 &p_{n-2}:=\rho_{n-2}+\theta_{n-2}-\frac{1}{3}|u_0|^2,\\ &p_{m+n-2}:=\rho_{m+n-2}+\theta_{m+n-2}-\frac{2}{3}u_{0}\cdot u_{m},\quad \text{for }m\geq 1,
		\end{aligned}
	\end{equation}
	and $\mathfrak{h}_{m-1}=(\mathfrak{h}_{m-1,1},\mathfrak{h}_{m-1,2},\mathfrak{h}_{m-1,3})$, $\mathfrak{q}_{m-1}$ are given by
	\begin{equation}\label{2.17-2}
		\begin{aligned}
			&\mathfrak{h}_{m-1,i}:=\partial_{t}\rho_{m+2-n} u_{0,i}-\sum_{j=1}^{3} \partial_{x_{i}}\left\langle\mathcal{A}_{i j}, G_{m-1}\right\rangle, \quad i=1,2,3,\\
			&\mathfrak{q}_{m-1}:=\frac{2}{5} \partial_{t}\left(\rho_{m}+\theta_{m}\right)+\partial_{t}\rho_{m+2-n} \theta_{0}-\frac{2}{5} \sum\limits_{i=1}^{3}\partial_{x_{i}}\left\langle\mathcal{B}_{i}, G_{m-1}\right\rangle\\
			&=\left\{
			\begin{aligned}
				&\frac{2}{5}\pa_{t}p_{m}-\frac{2}{5}\sum\limits_{i=1}^3\partial_{x_{i}}\la \mathcal{B}_{i}, G_{m-1}\ra,\quad \text{for }1\leq m\leq n-3,\\
				&\frac{2}{5}\pa_{t}p_{n-2}+\frac{2}{15}\pa_{t}|u_{0}|^2+\pa_{t}\rho_{0}\t_{0}-\frac{2}{5}\sum\limits_{i=1}^3\partial_{x_{i}}\la \mathcal{B}_{i}, G_{n-3}\ra,\quad \text{for }m=n-2,\\
				&\frac{2}{5}\pa_{t}p_{m}+\frac{4}{15}\pa_{t}(u_{0}\cdot u_{m+2-n})+\pa_{t}\rho_{m+2-n}\t_{0}-\frac{2}{5}\sum\limits_{i=1}^3\partial_{x_{i}}\la \mathcal{B}_{i}, G_{m-1}\ra,\quad \text{for }m\geq n-1,\\
			\end{aligned}
			\right.
		\end{aligned}
	\end{equation}
	with $\mathcal{A}_{ij}$ and $\mathcal{B}_{i}$ are the Burnett functions defined in \eqref{2.4-1}, and $G_{m-1}$ is defined by
	\begin{equation}\label{2.13}
		\begin{aligned}
			G_{m-1}&=
		 \mathbf{L}^{-1}\big\{(\mathbf{I}-\mathbf{P})[-\partial_{t} f_{m-n}-v \cdot \nabla_{x} f_{m-2}]+\sum_{i+j=m \atop i \geq 1, j \geq 1} \frac{1}{2}[\Gamma\left(f_{i}, f_{j}\right)+\Gamma\left(f_{j}, f_{i}\right)]\\
			&\qquad \qquad +[\Gamma\left(f_{0},(\mathbf{I}-\mathbf{P}) f_{m}\right)+\Gamma\left((\mathbf{I}-\mathbf{P}) f_{m}, f_{0}\right)]\big\}\\
		&=(\mathbf{I}-\mathbf{P})f_{m+n-2}-\mathbf{L}^{-1}\{\Gamma(f_{0},\mathbf{P}f_{m})+\Gamma(\mathbf{P}f_{m},f_{0})\}\quad \text{for }m\geq 1,
		\end{aligned}
	\end{equation}
which depends only on $f_{i}$ with $1\leq i\leq m-1$. Here $p_{k}(t)$ in \eqref{2.17-1} with $0\leq k\leq n-3$ means a function depending on $t$, and will be determined in Lemma \ref{lem2.2} below to guarantee the compatibility conditions for the equations of velocity.
\end{lemma}

\subsection{Viscous boundary layer}
To formulate $\bar{F}_{i}^{\pm}$, we define the macro-micro decomposition of $\bar{F}_{i}^{\pm}$ as
$$
\begin{aligned}
	\frac{\bar{F}_{k}^{\pm}}{\sqrt{\mu}}&=\bar{f}_{k}^{\pm}=\mathbf{P}\bar{f}_{k}^{\pm}+(\mathbf{I-P})\bar{f}_{k}^{\pm}\\
	&=(\bar{\rho}_{k}^{\pm}+v\cdot \bar{u}_{k}^{\pm}+\frac{\bar{\theta}_{k}^{\pm}}{2}(|v|^2-3))\sqrt{\mu}+(\mathbf{I-P})\bar{f}_{k}^{\pm}\quad (k \geq 1).
\end{aligned}
$$
It follows from $\eqref{1.14}$ that
\begin{equation}\label{3.1}
	\begin{aligned}
		&\bar{f}_{i}^{\pm}\equiv \mathbf{P} \bar{f}_{i}^{\pm}=\{\bar{\rho}_{i}^{\pm}+\bar{u}_{i}^{\pm} \cdot v+\frac{1}{2}\bar{\theta}_{i}^{\pm}(|v|^{2}-3)\} \sqrt{\mu}, \quad1 \leq i \leq n-2, \\
		&(\mathbf{I}-\mathbf{P}) \bar{f}_{n-1}^{\pm}=\mathbf{L}^{-1}\{\Gamma(f_{0}^{\pm},\bar{f}_{1}^{\pm})+\Gamma(\bar{f}_{1}^{\pm},f_{0}^{\pm})\},\\
		&\begin{aligned}
			(\mathbf{I}-\mathbf{P}) \bar{f}_{n}^{\pm}=&\mathbf{L}^{-1}\big\{(\mathbf{I}-\mathbf{P})(\pm v_{3} \partial_{y} \bar{f}_{1}^{\pm})+ \Gamma(\bar{f}_{1}^{\pm}, \bar{f}_{1}^{\pm})\\
			&\qquad \quad +\sum_{i+j=2\atop i\geq 0,j\geq 1}[\Gamma\left(f_{i}^{\pm}, \bar{f}_{j}^{\pm}\right)+\Gamma\left(\bar{f}_{j}^{\pm}, f_{i}^{\pm}\right)] \\	
			&\qquad\quad+\sum_{i+j+l=2 \atop i \geq 0, j \geq 1, 1\leq l\leq \fb}\frac{(\mp y)^{l}}{l !} [\Gamma\left(\partial_{x_3}^{l} f_{i}^{\pm}, \bar{f}_{j}^{\pm}\right)+\Gamma\left(\bar{f}_{j}^{\pm}, \partial_{x_3}^{l} f_{i}^{\pm}\right)] \big\},
		\end{aligned}\\
		&\begin{aligned}(\mathbf{I-P})\bar{f}_{k+n}^{\pm}=&\mathbf{L}^{-1}\big\{(\mathbf{I}-\mathbf{P})(-\partial_{t}\bar{f}_{k+2-n}^{\pm}-v_{\sp}\cdot \nabla_{\sp}\bar{f}_{k}^{\pm}\pm v_{3}\partial_{y}\bar{f}_{k+1}^{\pm})\\
			&\qquad \quad +\sum\limits_{i+j=k+2\atop i\geq 0,j\geq 1}[\Gamma(f_{i}^{\pm},\bar{f}_{j}^{\pm})+\Gamma(\bar{f}_{j}^{\pm},f_{i}^{\pm})]\\
			&\qquad \quad +\sum\limits_{i+j+l=k+2\atop i\geq 0,j\geq 1,1\leq l\leq \fb}\frac{(\mp y)^{l}}{l!}[\Gamma(\partial_{x_3}^{l}f_{i}^{\pm},\bar{f}_{j}^{\pm})+\Gamma(\bar{f}_{j}^{\pm},\partial_{x_3}^{l}f_{i}^{\pm})]\\
			&\qquad\quad+\sum\limits_{i+j=k+2\atop i,j\geq 1}\frac{1}{2}[\Gamma(\bar{f}_{i}^{\pm},\bar{f}_{j}^{\pm})+\Gamma(\bar{f}_{j}^{\pm},\bar{f}_{i}^{\pm})]\big\},\quad k\geq 1.
		\end{aligned}
	\end{aligned}
\end{equation}
where we have used the notations
$$
(\partial_{x_3}^{l})f_{i}^{\pm}:=\frac{1}{\sqrt{\mu}}(\partial_{x_3}^{l})F_{i}^{\pm}\quad \text{for } 0\leq l\leq \fb.
$$

Projecting $\eqref{1.14}_3$ onto $1, v$ and using \eqref{2.2}, we get
\begin{equation}\label{3.3}
	\left\{\begin{aligned}
		&\partial_{y}\bar{u}_{1,3}^{\pm}=0, \\
		&\partial_{y}\left(\bar{\rho}_{1}^{\pm}+\bar{\theta}_{1}^{\pm}\right)=0.
	\end{aligned}\right.
\end{equation}
Throughout the present paper, we always assume the far field condition
\begin{equation}\label{3.3-3}
	\lim _{y \rightarrow +\infty} \bar{f}_{k}^{\pm}\left(t, x_{\sp}, y\right)=0,\quad \text{for }k\geq 1.
\end{equation}
which, together with \eqref{3.3}, yields that
\begin{equation}\label{3.3-1}
	\bar{u}_{1,3}^{\pm} \equiv 0,\quad \bar{\rho}_{1}^{\pm}+\bar{\theta}_{1}^{\pm}\equiv 0.
\end{equation}
It is noted that $\eqref{3.3-1}_2$ is the Boussinesq relation in the diffusive limit of Boltzmann equation.

For later use, we define the viscosity and thermal conductivity coefficients $\lambda$, $\kappa$ of the viscous boundary layer
\begin{equation}\nonumber
		\lambda:= \la \mathcal{A}_{ij},\mathbf{L}^{-1}\mathcal{A}_{ij}\ra,\quad \forall i\neq j,\quad \kappa:=\la \mathcal{B}_{i},\mathbf{L}^{-1}\mathcal{B}_{i}\ra\quad \forall i=1,2,3.
\end{equation}
It follows from \cite[Lemma 4.4]{Bardos-2} that $\la \mathcal{A}_{ii},\mathbf{L}^{-1}\mathcal{A}_{ii}\ra=\frac{4}{3}\lambda$ with $i=1,2,3$.

The following Proposition \ref{prop3.1} is devoted to the reformulation of viscous boundary layer, whose proof is given in Appendix \ref{AppendixD} for simplicity of presentation.

\begin{proposition}[Reformulation of viscous boundary layer]\label{prop3.1}
	Let $\bar{f}_{k}^{\pm}(k \geq 1)$ be the solution of $\eqref{1.14}$. Then, $\left(\bar{\rho}_{k,\sp}^{\pm},\bar{u}_{k,\sp}^{\pm}, \theta_{k}^{\pm}\right) $ for $ k \geq 1,$ satisfies
	\begin{equation}\label{3.12}
		\begin{aligned}
			&\partial_{t} \bar{u}_{k,1}^{\pm}+( y \partial_{x_3} u_{0,3}^{\pm}\mp u_{1,3}^{\pm}) \partial_{y} \bar{u}_{k, 1}^{\pm}-\lambda \partial_{y}^{2} \bar{u}_{k, 1}^{\pm}+(u_{0,\sp}^{\pm} \cdot \nabla_{\sp}) \bar{u}_{k,1}^{\pm}+(\bar{u}_{k,\sp}^{\pm} \cdot \nabla_{\sp}) u_{0,1}^{\pm}\\
			&\qquad \quad =\pm \partial_{y}[(u_{1,1}^{\pm}\mp y \partial_{x_3} u_{0,1}^{\pm}+\bar{u}_{1,1}^{\pm}) \bar{u}_{k, 3}^{\pm}]-\sum_{j=1}^{2} \partial_{j}\la\mathcal{A}_{1j}, W_{k-1}^{\pm}\ra\pm \partial_{y}\la\mathcal{A}_{31}, J_{k-1}^{\pm}\ra\\
			&\qquad\quad\quad\mp \partial_{x_{1}} \int_{y}^{+\infty} H_{k-1}^{\pm}\left(t, x_{\sp}, s\right) d s+u_{0,1}^{\pm} \partial_{t} \bar{\rho}_{k+2-n}^{\pm}:=\mathfrak{f}_{k-1,1}^{\pm},\\
			&\partial_{t} \bar{u}_{k, 2}^{\pm}+(y \partial_{x_3} u_{0,3}^{\pm}\mp u_{1,3}^{\pm}) \partial_{y} \bar{u}_{k, 2}^{\pm}-\lambda \partial_{y}^{2} \bar{u}_{k, 2}^{\pm}+(u_{0,\sp}^{\pm} \cdot \nabla_{\sp}) \bar{u}_{k, 2}^{\pm}+(\bar{u}_{k,\sp}^{\pm} \cdot \nabla_{\sp}) u_{0,2}^{\pm}\\
			&\qquad\quad=\pm\partial_{y}[(u_{1,2}^{\pm}\pm y\partial_{x_3} u_{0,2}^{\pm}+\bar{u}_{1,2}^{\pm}) \bar{u}_{k, 3}^{\pm}]-\sum_{j=1}^{2} \partial_{j}\la \mathcal{A}_{2 j}, W_{k-1}^{\pm}\ra\pm \partial_{y}\la\mathcal{A}_{32}, J_{k-1}^{\pm}\ra\\
			&\qquad\quad\quad\mp \partial_{x_2} \int_{y}^{+\infty} H_{k-1}^{\pm}\left(t, x_{\sp}, s\right) d s+u_{0,2}^{\pm} \partial_{t} \bar{\rho}_{k+2-n}^{\pm}:=\mathfrak{f}_{k-1,2}^{\pm},\\
			&\partial_{t} \bar{\theta}_{k}^{\pm}+(y \partial_{x_3} u_{0,3}^{\pm}\mp u_{1,3}^{\pm}) \partial_{y} \bar{\theta}_{k}^{\pm}-\frac{2}{5} \kappa \partial_{y}^{2} \bar{\theta}_{k}^{\pm}+(u_{0,\sp}^{\pm} \cdot \nabla_{\sp}) \bar{\theta}_{k}^{\pm}+(\bar{u}_{k,\sp}^{\pm} \cdot \nabla_{\sp}) \theta_{0}^{\pm}\\
			&\qquad\quad=\frac{2}{5} \partial_{t}(\bar{\rho}_{k}^{\pm}+\bar{\theta}_{k}^{\pm})\pm \partial_{y}[(\mp y \partial_{x_3} \theta_{0}^{\pm}+\bar{\theta}_{1}^{\pm}+\theta_{1}^{\pm}) \bar{u}_{k,3}^{\pm}]\\
			&\qquad\quad\quad\pm \frac{2}{5} \partial_{y}\la \mathcal{B}_{3}, J_{k-1}^{\pm}\ra -\frac{2}{5} \sum_{j=1}^{2} \partial_{j}\la\mathcal{B}_{j}, W_{k-1}^{\pm}\ra+\theta_{0}^{\pm} \partial_{t} \bar{\rho}_{k+2-n}^{\pm}:=\mathfrak{g}_{k-1}^{\pm},\\
			&\bar{\rho}_{k}^{\pm}+\bar{\t}_{k}^{\pm}=\frac{2}{3}u_{0,\sp}^{\pm}\cdot \bar{u}_{k+2-n,\sp}^{\pm}\pm \int_{y}^{\infty}H_{k+2-n}^{\pm}(t,x_{\sp},s)ds:=\mathfrak{s}_{k-1}^{\pm},
		\end{aligned}
	\end{equation}
	where we have used the notations $\bar{u}_{k}^{\pm}=(\bar{u}_{k,1}^{\pm},\bar{u}_{k,2}^{\pm},\bar{u}_{k,3}^{\pm})=(\bar{u}_{k,\sp}^{\pm},\bar{u}_{k,3}^{\pm})$,
	$$
	(u_{k}^{+},\theta_{k}^{+}):=(u_{k},\theta_{k})(x_{1},x_{2},1),\quad (u_{k}^{-},\theta_{k}^{-}):=(u_{k},\theta_{k})(x_{1},x_{2},0).
	$$
	It holds that $\mathfrak{f}_{k-1,1}^{\pm},\mathfrak{f}_{k-1,2}^{\pm}, \mathfrak{g}_{k-1}^{\pm}$ and $\mathfrak{s}_{k-1}^{\pm}$ depend on $f_{i}(0\leq i\leq k)$ and $\bar{f}_{j}^{\pm}(1\leq j\leq k-1)$ with
	\begin{align}
		&\begin{aligned}
			W_{k-1}^{\pm}&:=\mathbf{L}^{-1} \Big\{-(\mathbf{I}-\mathbf{P})(\partial_{t} \bar{f}_{k-n}^{\pm}+v_{\sp} \cdot \nabla_{\sp} \bar{f}_{k-2}^{\pm})\pm(\mathbf{I}-\mathbf{P})[v_{3} \partial_{y} \bar{f}_{k-1}^{\pm}]\\
			&\quad +\sum_{i+j=k \atop i \geq 1, j \geq 1} [\Gamma\left(f_{i}^{\pm}, \bar{f}_{j}^{\pm}\right)+\Gamma\left(\bar{f}_{j}^{\pm}, f_{i}^{\pm}\right)]+\sum_{i+j=k\atop i,j\geq 1} \frac{1}{2} [\Gamma(\bar{f}_{i}^{\pm}, \bar{f}_{j}^{\pm})+\Gamma(\bar{f}_{j}^{\pm}, \bar{F}_{i}^{\pm})]\\
			&\quad+\sum_{i+j+l=k\atop i\geq 0, j \geq 1,1\leq l\leq \fb} \frac{(\mp y)^{l}}{l !} [\Gamma(\partial_{x_3}^{l} f_{i}^{\pm}, \bar{f}_{j}^{\pm})+Q(\bar{f}_{j}^{\pm}, \partial_{x_3}^{l} f_{i}^{\pm})]\\
			&\quad +[\Gamma(f_{0}^{\pm},(\mathbf{I}-\mathbf{P}) \bar{f}_{k}^{\pm})+\Gamma((\mathbf{I}-\mathbf{P}) \bar{f}_{k}^{\pm}, f_{0}^{\pm})]\Big\},
		\end{aligned}\label{3.15}\\
		&\begin{aligned}
			H_{k-1}^{\pm}:=-\partial_{t} \bar{u}_{k-1,3}^{\pm}-\sum\limits_{j=1}^{2}\pa_{j}\la \mathcal{A}_{3j},(\mathbf{I-P})\bar{f}_{k+n-3}^{\pm}\ra \pm\pa_{y}\la \mathcal{A}_{33},W_{k-1}^{\pm}\ra,
		\end{aligned}\label{3.17}
	\end{align}
and
\begin{equation}\label{3.13}
\begin{aligned}
	J_{k-1}^{\pm}&:=\mathbf{L}^{-1} \Big\{-(\mathbf{I}-\mathbf{P})(\partial_{t} \bar{f}_{k+1-n}^{\pm}+v_{\sp} \cdot \nabla_{\sp} \bar{f}_{k-1}^{\pm})\pm(\mathbf{I}-\mathbf{P})[v_{3} \partial_{y}(\mathbf{I}-\mathbf{P}) \bar{f}_{k}^{\pm}]\\
	&\quad+\sum_{i+j=k+1 \atop i \geq 2, j \geq 1} [\Gamma(f_{i}^{\pm}, \bar{f}_{j}^{\pm})+\Gamma(\bar{f}_{j}^{\pm}, f_{i}^{\pm})]+\sum_{i+j=k+1\atop i,j\geq 2}\frac{1}{2} [\Gamma(\bar{f}_{i}^{\pm}, \bar{f}_{j}^{\pm})+\Gamma(\bar{f}_{j}^{\pm}, \bar{f}_{i}^{\pm})] \\
	&\quad+\sum_{i+j+l=k+1 \atop i \geq 0,1 \leq j \leq k-1,1 \leq l \leq \fb}\frac{(\mp y)^{l}}{l !} [\Gamma(\partial_{x_3}^{l} f_{i}^{\pm}, \bar{f}_{j}^{\pm})+\Gamma(\bar{f}_{j}^{\pm}, \partial_{x_3}^{l} f_{i}^{\pm})]\\
	&\quad +[\Gamma(f_{0}^{\pm},(\mathbf{I}-\mathbf{P}) \bar{f}_{k+1}^{\pm})+\Gamma((\mathbf{I}-\mathbf{P}) \bar{f}_{k+1}^{\pm}, f_{0}^{\pm})]\\
	&\quad +[\Gamma(f_{1}^{\pm},(\mathbf{I}-\mathbf{P}) \bar{f}_{k}^{\pm})+\Gamma((\mathbf{I}-\mathbf{P}) \bar{f}_{k}^{\pm}, f_{1}^{\pm})]\\
	&\quad +[\Gamma((\mathbf{I-P})f_{1}^{\pm},\mathbf{P} \bar{f}_{k}^{\pm})+\Gamma(\mathbf{P} \bar{f}_{k}^{\pm}, (\mathbf{I-P})f_{1}^{\pm})]\\
	&\quad\mp y\cdot[\Gamma(\partial_{x_3} f_{0}^{\pm},(\mathbf{I}-\mathbf{P}) \bar{f}_{k}^{\pm})+\Gamma((\mathbf{I}-\mathbf{P}) \bar{f}_{k}^{\pm}, \partial_{x_3} f_{0}^{\pm})] \\
	&\quad+[\Gamma(\bar{f}_{1}^{\pm},(\mathbf{I}-\mathbf{P}) \bar{f}_{k}^{\pm})+\Gamma((\mathbf{I}-\mathbf{P}) \bar{f}_{k}^{\pm}, \bar{f}_{1}^{\pm})]\Big\}.
\end{aligned}
\end{equation}

Moreover, once one solves $(u_{k,\sp}^{\pm},\theta_{k}^{\pm})$, then $\bar{u}_{k+1,3}^{\pm}, \partial_{y}\left(\bar{\rho}_{n-2+k}^{\pm}+\bar{\theta}_{n-2+k}^{\pm}\right),(\mathbf{I}-\mathbf{P}) \bar{f}_{k+n-2}^{\pm}$ can be determined by the following equations
	\begin{align}
		\partial_{y}\bar{u}_{k+1,3}^{\pm}&=\pm\partial_{t}\bar{\rho}_{k+2-n}^{\pm}\pm\operatorname{div}_{\sp}\bar{u}_{k,\sp}^{\pm},\label{3.21}\\
		(\mathbf{I}-\mathbf{P}) \bar{f}_{k+n-2}^{\pm}&=\mathbf{L}^{-1}\{[\Gamma(f_{0}^{\pm}, \mathbf{P} \bar{f}_{k}^{\pm})+\Gamma( \mathbf{P} \bar{f}_{k}^{\pm}, f_{0}^{\pm})]\}+W_{k-1}^{\pm},\label{3.18}\\
		\partial_{y}(\bar{\rho}_{n-2+k}^{\pm}+\bar{\theta}_{n-2+k}^{\pm})&=\frac{2}{3} \partial_{y}(u_{0,\sp}^{\pm} \cdot \bar{u}_{k, \sp}^{\pm})\mp H_{k-1}^{\pm},\label{3.20}
	\end{align}
\end{proposition}

\begin{remark}\label{rem2.1}
	From \eqref{3.15}-\eqref{3.17}, it is clear  that $W_{k-1}^{\pm}, H_{k-1}^{\pm}$ and $\mathfrak{s}_{k-1}^{\pm}$ depend on $f_{i}^{\pm}, \bar{f}_{i}^{\pm}$ up to $(k-1)-$th order.
	
	On the other hand,   it is clear that $J_{k-1}^{\pm}$ depends on $\bar{f}_{k}^{\pm}$ due to the terms $(\mathbf{I}-\mathbf{P})\bar{f}_{k+1}^{\pm}$ and $\mathbf{P}\bar{f}_{k}^{\pm}$, see \eqref{3.13} for details. Fortunately, we note that the equations \eqref{3.12} for $(\bar{\rho}_{k,\sp}^{\pm},\bar{u}_{k,\sp}^{\pm}, \theta_{k}^{\pm}) $ depend only on  $\la \mathcal{A}_{3i}, J_{k-1}^{\pm}\ra, i=1,2$ and $\la \mathcal{B}_{3}, J_{k-1}^{\pm}\ra$ which are indeed independent of $\bar{f}_{k}^{\pm}$. For instance, one can show
	\begin{equation}\label{3.20-1}
		\la \mathcal{A}_{3i}, J_{0}^{\pm}\ra=\la \mathcal{B}_{3}, J_{0}^{\pm}\ra=0,\quad i=1,2.
	\end{equation}
	In fact, taking $k=1$ in \eqref{3.13}, one has
	$$
	\begin{aligned}
		J_{0}^{\pm}&=\mathbf{L}^{-1}\big\{[\Gamma(f_{0}^{\pm},(\mathbf{I}-\mathbf{P})\bar{f}_{2}^{\pm})+\Gamma((\mathbf{I}-\mathbf{P})\bar{f}_{2}^{\pm},f_{0}^{\pm})]\\
		&\qquad \quad +[\Gamma((\mathbf{I}-\mathbf{P})f_{1}^{\pm},\mathbf{P}\bar{f}_{1}^{\pm})+\Gamma(\mathbf{P}\bar{f}_{1}^{\pm},(\mathbf{I}-\mathbf{P})f_{1}^{\pm})]\big\}\\
		&:=I_{1}+I_{2}.
	\end{aligned}
	$$
	
	If $n\geq 4$, it follows from \eqref{2.0} and \eqref{3.1} that $(\mathbf{I}-\mathbf{P})f_{1}=(\mathbf{I}-\mathbf{P})\bar{f}_{2}^{\pm}=0$, and hence $J_{0}\equiv 0$. Therefore, \eqref{3.20-1} holds automatically.
	
	If $n=3$, for $I_{1}$, it follows from \eqref{3.1} and \cite[p. 649]{Guo2006} that
	$$
	\begin{aligned}
		(\mathbf{I}-\mathbf{P})\bar{f}_{2}^{\pm}&=\mathbf{L}^{-1}\big\{\Gamma(f_{0}^{\pm},\bar{f}_{1}^{\pm})+\Gamma(\bar{f}_{1}^{\pm},{f}_{0})\big\}\\
		&=\sum\limits_{l,s=1}^{3}u_{0,l}^{\pm}\bar{u}_{1,s}^{\pm}\mathcal{A}_{ls}+\sum\limits_{l=1}^{3}\t_{0}^{\pm}\bar{u}_{1,l}^{\pm}\mathcal{B}_{l}+\frac{1}{4}\t_{0}^{\pm}\bar{\t}_{1}^{\pm}(\mathbf{I}-\mathbf{P})\{(|v|^2-5)^2\sqrt{\mu}\}.
	\end{aligned}
	$$
	Noting $f_{0}^{\pm}=\mathbf{P}f_{0}^{\pm}$ and $\mathbf{L}^{-1}$ preserves the odevity of $v$, we have
	$$
	\begin{aligned}
		\la \mathcal{A}_{3i}, I_{1}\ra&=\rho_{0}^{\pm}(u_{0,3}^{\pm}\bar{u}_{1,i}+u_{0,i}^{\pm}\bar{u}_{1,3}^{\pm})\la \mathcal{A}_{3i},\mathbf{L}^{-1}\{\Gamma(\sqrt{\mu}, \mathcal{A}_{3i})+\Gamma(\mathcal{A}_{3i},\sqrt{\mu})\}\ra\\
		&\quad +\frac{\theta_{0}^{\pm}}{2}(u_{0,3}^{\pm}\bar{u}_{1,i}+u_{0,i}^{\pm}\bar{u}_{1,3}^{\pm})\la \mathcal{A}_{3i},\mathbf{L}^{-1}\{\Gamma((|v|^2-3)\sqrt{\mu}, \mathcal{A}_{3i})+\Gamma(\mathcal{A}_{3i},(|v|^2-3)\sqrt{\mu})\}\ra
		\\
		&\equiv 0,\quad i=1,2,
	\end{aligned}
	$$
	where we have used $u_{0,3}^{\pm}=\bar{u}_{1,3}^{\pm}=0$. Using \eqref{2.0}, a similar calculation shows that $\la \mathcal{A}_{3i}, I_{2}\ra=0$. Similarly, one can get $\la \mathcal{B}_{3}, J_{0}^{\pm}\ra=0$. Hence, we conclude \eqref{3.20-1}.
	
	For $k\geq 2$, by similar arguments, one has the highest order terms of $\la \mathcal{A}_{3i}, J_{k-1}^{\pm}\ra(i=1,2)$ and $\la \mathcal{B}_{3}, J_{k-1}^{\pm}\ra$ depend only on $\bar{u}_{k,3}^{\pm}$, which, together with \eqref{3.21}, implies that they depend on $\bar{f}_{i}^{\pm}(1\leq i\leq k-1)$, and so are $\mathfrak{f}_{k-1,1}^{\pm}, \mathfrak{f}_{k-1,2}^{\pm}$ and $\mathfrak{g}_{k-1}^{\pm}$, one can see the last two paragraphs in the proof of Proposition \ref{prop3.1} for more details.
\end{remark}

\subsection{Knudsen boundary layer}
Let $\hat{f}_{k}^{\pm}:=\frac{\hat{F}_{k}^{\pm}}{\sqrt{\mu}}$,  then we rewrite \eqref{1.19-1} as
\begin{equation}\label{4.1}
	\mp v_{3} \frac{\partial \hat{f}_{k}^{\pm}}{\partial \eta}+\mathbf{L} \hat{f}_{k}^{\pm}=\hat{S}_{k}^{\pm},\quad k \geq 1,
\end{equation}
where $\hat{S}_{k}^{\pm}=\hat{S}_{k,1}^{\pm}+\hat{S}_{k,2}^{\pm}$ $(k \geq 1)$ with
\begin{equation}\label{4.3}
	\begin{aligned}
		&\hat{S}_{i,1}^{\pm}=\hat{S}_{i, 2}^{\pm}=\hat{S}_{n-1,1}^{\pm}=\hat{S}_{n, 1}^{\pm}\equiv 0, \quad 1 \leq i \leq n-2,\\
		&\hat{S}_{n-1,2}^{\pm}=[\Gamma(f_{0}^{\pm}, \hat{f}_{1}^{\pm})+\Gamma(\hat{f}_{1}^{\pm}, f_{0}^{\pm})],\\
		&\hat{S}_{n,2}^{\pm}=[\Gamma(\bar{f}_{1}^{\pm,0},  \hat{f}_{1}^{\pm})+\Gamma(\hat{f}_{1}^{\pm}, \bar{f}_{1}^{\pm,0})]+\sum_{i+j=2 \atop i \geq 0, j \geq 1} [\Gamma(f_{i}^{\pm},  \hat{f}_{j}^{\pm})+\Gamma( \hat{f}_{j}^{\pm}, f_{i}^{\pm})]+\Gamma(\hat{f}_{1}^{\pm},\hat{f}_{1}^{\pm}),\\
		&\hat{S}_{n+k,1}^{\pm}=-\mathbf{P}\{\partial_{t} \hat{f}_{k+2-n}^{\pm}+v_{\sp} \cdot \nabla_{\sp} \hat{f}_{k}^{\pm}\},\quad k\geq 1,
	\end{aligned}
\end{equation}
and
\begin{equation}\label{4.3-1}
\begin{aligned}
	\hat{S}_{n+k, 2}^{\pm}&=-(\mathbf{I}-\mathbf{P})\{\partial_{t} \hat{f}_{k+2-n}^{\pm}+v_{\sp} \cdot \nabla_{\sp} \hat{f}_{k}^{\pm}\}+\sum_{i+j=k+2\atop i\geq 0,j \geq 1}[\Gamma(f_{i}^{\pm},  \hat{f}_{j}^{\pm})+\Gamma( \hat{f}_{j}^{\pm}, f_{i}^{\pm})]\\
	&\quad+\sum_{i+j=k+2\atop i,j \geq 1} \big\{[\Gamma(\bar{f}_{i}^{\pm,0},  \hat{f}_{j}^{\pm})+\Gamma(\hat{f}_{j}^{\pm}, \bar{f}_{i}^{\pm,0})] +\frac{1}{2}[\Gamma( \hat{f}_{i}^{\pm}, \hat{f}_{j}^{\pm})+\Gamma(\hat{f}_{i}^{\pm},  \hat{f}_{j}^{\pm})]\big\}\\
	&\quad+\sum_{i+j+nl=k+2 \atop i \geq 0, j \geq 1,1 \leq \ell \leq\fb} \frac{(\mp\eta)^{l}}{l!} \frac{1}{\sqrt{\mu}}[\Gamma(\partial_{x_3}^{l} f_{i}^{\pm},  \hat{f}_{j}^{\pm})+\Gamma(\hat{f}_{j}^{\pm}, \partial_{x_3}^{l} f_{i}^{\pm})]\\
	&\quad+\sum_{i+j+(n-1) l=k+2\atop i, j \geq 1,1 \leq \ell \leq \fb} \frac{\eta^{l}}{l !} [\Gamma(\partial_{y}^{l} \bar{f}_{i}^{\pm,0}, \hat{f}_{j}^{\pm})+\Gamma(\hat{f}_{j}^{\pm}, \partial_{y}^{l} \bar{f}_{i}^{\pm,0})] \quad k \geq 1.
\end{aligned}
\end{equation}
It is easy to notice that $\hat{S}_{k,1}^{\pm} \in \mathcal{N}, \hat{S}_{k, 2}^{\pm} \in \mathcal{N}^{\perp}$. We assume that
\begin{equation}\label{4.4}
	\hat{S}_{k,1}^{\pm}=\{\hat{a}_{k}^{\pm}+\hat{b}_{k}^{\pm} \cdot v+\hat{c}_{k}^{\pm}|v|^{2}\} \sqrt{\mu},
\end{equation}
where $(\hat{a}_{k}^{\pm}, \hat{b}_{k}^{\pm}, \hat{c}_{k}^{\pm})=(\hat{a}_{k}^{\pm}, \hat{b}_{k}^{\pm}, \hat{c}_{k}^{\pm})(t, x_{\sp}, \eta)$. By similar arguments as in \cite{Bardos}, we have the following lemma. The details of proof are omitted for simplicity of presentation.

\begin{lemma}\label{lem4.1}
	For $(\hat{a}_{k}^{\pm},\hat{b}_{k}^{\pm},\hat{c}_{k}^{\pm})$ defined in \eqref{4.4}, we assume that
	$$
	\lim\limits_{\eta\to \infty}e^{\zeta\eta}|(\hat{a}_{k}^{\pm},\hat{b}_{k}^{\pm},\hat{c}_{k}^{\pm})(t,x_{\sp},\eta)|=0,
	$$
	for some positive constant $\zeta>0$. Then there exists a function
	\begin{equation}\label{4.5}
		\hat{f}_{k,1}^{\pm}=\{\hat{A}_{k}^{\pm}v_{3}+\hat{B}_{k, 1}^{\pm} v_{3}v_{1}+\hat{B}_{k,2}^{\pm} v_{3}v_{2}+\hat{B}_{k,3}^{\pm}+\hat{C}_{k}^{\pm} v_{3}|v|^{2}\} \sqrt{\mu},
	\end{equation}
	such that
	$$
	\mp v_{3} \partial_{\eta} \hat{f}_{k,1}^{\pm}-\hat{S}_{k, 1}^{\pm} \in \mathcal{N}^{\perp},
	$$
	where
	\begin{equation}\label{4.6}
		\begin{aligned}
			\hat{A}_{k}^{\pm}\left(t, x_{\sp}, \eta\right) &=-\int_{\eta}^{\infty}\left(2 \hat{a}_{k}^{\pm}+3 \hat{c}_{k}^{\pm}\right)\left(t, x_{\sp}, s\right) d s, \\
			\hat{B}_{k, i}^{\pm}\left(t, x_{\sp}, \eta\right) &=-\int_{\eta}^{\infty} \hat{b}_{k,i}^{\pm}\left(t, x_{\sp}, s\right) d s \qquad i=1,2, \\
			\hat{B}_{k,3}^{\pm}\left(t, x_{\sp}, \eta\right) &=-\int_{\eta}^{\infty} \hat{b}_{k,3}^{\pm}\left(t, x_{\sp}, s\right) d s, \\
			\hat{C}_{k}^{\pm}\left(t, x_{\sp}, \eta\right) &=\frac{1}{5} \int_{\eta}^{\infty} \hat{a}_{k}^{\pm}\left(t, x_{\sp}, s\right) d s.
		\end{aligned}
	\end{equation}
	Moreover it holds that
	$$
	|\mp v_{3} \partial_{\eta} \hat{f}_{k, 1}^{\pm}-\hat{S}_{k,1}^{\pm}| \leq C|(\hat{a}_{k,\pm}, \hat{b}_{k}^{\pm}, \hat{c}_{k}^{\pm})(t, x_{\sp}, \eta)|(1+|v|)^{4} \sqrt{\mu}
	$$
	and
	$$
	|\hat{f}_{k,1}^{\pm}(t, x_{\sp}, \eta, v)| \leq C(1+|v|)^{3} \sqrt{\mu} \int_{\eta}^{\infty}|(\hat{a}_{k}^{\pm}, \hat{b}_{k}^{\pm}, \hat{c}_{k}^{\pm})| \rightarrow 0 \text { as } \eta \rightarrow \infty.
	$$
\end{lemma}
\begin{remark}\label{rem4.1}
	It is very important to note that $\hat{S}_{k,1}^{\pm}$ depends only on $\hat{f}_{k+2-2n}^{\pm}$ and $\hat{f}_{k-n}^{\pm}$, which are already known functions when we consider the existence of $\hat{f}_{k}^{\pm}$. Thus $\hat{f}_{k,1}^{\pm}$ (or equivalent $(\hat{A}_{k}^{\pm}, \hat{B}_{k}^{\pm}, \hat{C}_{k}^{\pm})$) is determined by $\hat{f}_{k+2-2 n}^{\pm}$ and $ \hat{f}_{k-n}^{\pm} .$ On the other hand, by \eqref{4.3}, one has
	\begin{equation}\nonumber
		\hat{f}_{i,1}^{\pm}=\hat{f}_{n-1,1}^{\pm}=\hat{f}_{n,1}^{\pm} \equiv 0, \quad1 \leq i \leq n-2,
	\end{equation}
	which yields that
	\begin{equation}\label{4.7-1}
	\begin{aligned}
		&(\hat{A}_{i}^{\pm}, \hat{B}_{i}^{\pm}, \hat{C}_{i}^{\pm})\equiv(0,0,0),\quad \text{for }1\leq i\leq n.\\
	\end{aligned}
	\end{equation}
\end{remark}

Now we consider $\hat{f}_{k,2}^{\pm}$, which satisfies
\begin{equation}\label{4.8}
	\mp v_{3} \frac{\partial \hat{f}_{k,2}^{\pm}}{\partial \eta}+\mathbf{L} \hat{f}_{k,2}^{\pm}=\hat{S}_{k,2}^{\pm}-\mathbf{L} \hat{f}_{k,1}^{\pm}\pm(v_{3} \partial_{\eta} \hat{f}_{k, 1}^{\pm}-\hat{S}_{k,1}^{\pm}) \in \mathcal{N}^{\perp},
\end{equation}
Then it is easy to check that
$$
\hat{f}_{k}^{\pm}:=\hat{f}_{k,1}^{\pm}+\hat{f}_{k,2}^{\pm},
$$
is a solution of $\eqref{4.1}$.

In the following, to obtain the existence of solution to the Knudsen boundary layers  problem with perturbed boundary conditions, we consider the following half-space linear problem
\begin{equation}\label{4.9}
	\left\{\begin{aligned}
		&v_{3} \partial_{\eta} f+\mathbf{L} f=S\left(t, x_{\sp}, \eta, v\right), \\
		&\left.f\left(t, x_{\sp}, 0, v_{\sp}, v_{3}\right)\right|_{v_{3}>0}=f\left(t, x_{\sp}, 0, v_{\sp},-v_{3}\right)+f_{b}\left(t, x_{\sp}, v_{\sp},-v_{3}\right), \\
		&\lim _{\eta \rightarrow \infty} f\left(t, x_{\sp}, \eta, v\right)=0,
	\end{aligned}\right.
\end{equation}
where $\eta \in \mathbb{R}_{+}$ and $\left(t, x_{\sp}\right) \in[0, \tau] \times \mathbb{T}^{2}$. In fact, we regard $\left(t, x_{\sp}\right) \in[0, \tau] \times \mathbb{T}^{2}$ as parameters in \eqref{4.9}. The function $f_{b}\left(t, x_{\sp}, v\right)$ is defined only for $v_{3}<0$, and we always assume that it is extended to be 0 for $v_{3}>0$.

\begin{lemma}[{\cite{Huang-Jiang-Wang}}]\label{lem4.2}
	Let $0\leq \fa<\frac{1}{2}$ and $\beta\geq 3$. For each $(t,x_{\sp})\in [0,\tau]\times \mathbb{T}^2$, we assume that
	\begin{equation*}
		S\in \mathcal{N}^{\perp}\quad \text{ and }\quad ||\omega_{\beta}\mu^{-\fa}f_{b}(t,x_{\sp},0,\cdot)||_{L_{v}^{\infty}}+||\nu^{-1}\omega_{\beta}\mu^{-\fa}e^{\zeta_{0}\eta}S(t,x_{\sp},\cdot,\cdot)||_{L_{\eta,v}^{\infty}}<\infty,
	\end{equation*}
	for some positive constant $\zeta_{0}>0$, and
	\begin{equation}\label{4.11}
		\left\{\begin{aligned}
			&\int_{\mathbb{R}^{3}} v_{3} f_{b}\left(t, x_{\sp}, v\right) \sqrt{\mu} d v  \equiv 0,\\
			&\int_{\mathbb{R}^{3}}v_{i} v_{3} f_{b}\left(t, x_{\sp}, v\right) \sqrt{\mu} d v \equiv 0\quad i=1,2, \\
			&\int_{\mathbb{R}^{3}}|v|^{2} v_{3} f_{b}\left(t, x_{\sp}, v\right) \sqrt{\mu} d v  \equiv 0.
		\end{aligned}\right.
	\end{equation}
	Then the boundary value problem \eqref{4.9}-\eqref{4.11} has a unique solution $f$ satisfying
	\begin{equation*}
		\begin{aligned}
			&\|w_{\beta} \mu^{-\fa} e^{\zeta\eta} f\left(t, x_{\sp}, \cdot, \cdot \right)\|_{L_{\eta, v}^{\infty}}+\|w_{\beta} \mu^{-\mathfrak{a}} f\left(t, x_{\sp}, 0, \cdot\right)\|_{L_{v}^{\infty}}\\
			&\leq \frac{C}{\zeta_{0}-\zeta}\Big(\big\|w_{\beta} \mu^{-\mathfrak{a}} f_{b}\left(t, x_{\sp}, \cdot\right)\big\|_{L_{v}^{\infty}}+\big\|\nu^{-1} w_{\beta} \mu^{-\fa} e^{\zeta_{0}\eta} S\left(t, x_{\sp}, \cdot, \cdot\right)\big\|_{L_{\eta, v}^{\infty}}\Big),
		\end{aligned}
	\end{equation*}
	for all $\left(t, x_{\sp}\right) \in[0, \tau] \times \mathbb{T}^{2}$, where $C>0$ is a positive constant independent of $\left(t, x_{\sp}\right)$, and $\zeta$ is any positive constant such that $0<\zeta<\zeta_{0} .$ Moreover, if $S$ is continuous in $\left(t, x_{\sp}, \eta, v\right) \in$ $[0, \tau] \times \mathbb{T}^{2} \times \mathbb{R}_{+} \times \mathbb{R}^{3}$ and $f_{b}$ is continuous in $\left(t, x_{\sp}, v_{\sp},-v_{3}\right) \in[0, \tau] \times \mathbb{T}^{2} \times \mathbb{R}^{2} \times \mathbb{R}_{+}$, then the solution $f$ is continuous away from the grazing set $[0, \tau] \times \gamma_{0}$.
\end{lemma}

\begin{remark}\label{rem4.4}
	As indicated in \cite[Remark 2.7]{Guo-Huang-Wang},
	it is hard to obtain  the normal derivatives estimates for the boundary value problem \eqref{4.9}. However, we can obtain the tangential and time derivatives estimates for the solution of \eqref{4.9}. In fact, by the same arguments in \cite[Remark 2.7]{Guo-Huang-Wang}, we have
	\begin{align}
		& \sum_{i+j\leq r} \sup_{t\in[0,\tau]} \{\|w_{\beta} \mu^{-\fa}e^{\zeta\eta}\partial_t^{i}\nabla_\sp^{j}f(t)\|_{L^\infty_{x_\sp,\eta,v}}
		+\|w_{\beta} \mu^{-\fa} \partial_t^{i}\nabla_\sp^{j}f(t,\cdot,0,\cdot)\|_{L^\infty_{x_\sp,v}}\}\nonumber\\
		&\leq \frac{C}{\zeta_0-\zeta}  \sup_{t\in[0,\tau]}\Big\{ \sum_{i+j\leq r} \|w_{\beta} \mu^{-\fa} \partial_t^{i}\nabla_\sp^{j}f_b(t)\|_{L^\infty_{x_\sp,v}}+
		\sum_{i+j\leq r} \|\nu^{-1}w_{\beta} \mu^{-\fa} e^{\zeta_0 \eta} \partial_t^{i}\nabla_\sp^{j}S(t)\|_{L^\infty_{x_\sp,\eta,v}}\Big\}\nonumber.
	\end{align}
\end{remark}

\

\section{Boundary conditions and criteria of uniqueness for  pressure terms}

\subsection{Boundary conditions}\label{sec2.4}
In this subsection, we derive the boundary conditions for interior expansions and viscous boundary layers from the solvability conditions \eqref{4.11} of Knudsen boundary layers. Noting \eqref{1.3}, 
and then comparing the order of $\varepsilon$, we have for $k\geq 1$ that
\begin{align}
	 &\hat{f}_{k,2}^{-}\left(t, x_{\sp}, 0,v_{\sp}, v_{3}\right)\vert_{ v_{3}>0}=\hat{f}_{k,2}^{-}\left(t, x_{\sp},0, v_{\sp},-v_{3}\right)\nonumber\\
	&\quad +f_{k}(t,x_{\sp},0,v_{\sp},-v_{3})+(\bar{f}_{k}^{-}+\hat{f}_{k,1}^{-})(t,x_{\sp},0,v_{\sp},-v_{3})\nonumber\\
	&\quad-\big[f_{k}(t,x_{\sp},0,v_{\sp},v_{3})+(\bar{f}_{k}^{-}+\hat{f}_{k,1}^{-})(t,x_{\sp},0,v_{\sp},v_{3})\big],\label{4.14}\\
	&\hat{f}_{k,2}^{+}\left(t, x_{\sp}, 1,v_{\sp}, v_{3}\right)\vert_{ v_{3}<0}=\hat{f}_{k,2}^{+}\left(t, x_{\sp},1, v_{\sp},-v_{3}\right)\nonumber\\
	&\quad +f_{k}(t,x_{\sp},1,v_{\sp},-v_{3})+(\bar{f}_{k}^{+}+\hat{f}_{k,1}^{+})(t,x_{\sp},1,v_{\sp},-v_{3})\nonumber\\
	&\quad-\big[f_{k}(t,x_{\sp},1,v_{\sp},v_{3})+(\bar{f}_{k}^{+}+\hat{f}_{k,1}^{+})(t,x_{\sp},1,v_{\sp},v_{3})\big],\label{4.14-1}
\end{align}
For simplicity of presentation, we introduce the following notations for interior expansions
$$
f_{k}^{+}(v):=f_{k}(t,x_{\sp},1,v_{\sp},v_{3}),\quad f_{k}^{-}(v):=f_{k}(t,x_{\sp},0,v_{\sp},v_{3}).
$$
Let
\begin{equation}\label{4.15}
	\hat{g}_{k}^{\pm}(t, x_{\sp}, v_{\sp}, v_{3})=\left\{\begin{aligned}
		&0, \quad \mp v_{3}>0, \\
		& [f_{k}^{\pm}(v)+(\bar{f}_{k}^{\pm}+\hat{f}_{k,1}^{\pm})(t, x_{\sp},0, v_{\sp}, v_{3})]\\
		&\,\,-[f_{k}^{\pm}(R_{x}v)+(\bar{f}_{k}^{\pm}+\hat{f}_{k,1}^{\pm})(t, x_{\sp},0, v_{\sp},-v_{3})], \quad \mp v_{3}<0,
	\end{aligned}\right.
\end{equation}
then
$$
\hat{f}_{k,2}^{\pm}(t, x_{\sp}, 0, v_{\sp}, v_{3})\vert_{\mp v_{3}>0}=\hat{f}_{k,2}^{\pm}(t, x_{\sp}, 0, v_{\sp},-v_{3})+\hat{g}_{k}^{\pm}(t, x_{\sp}, v_{\sp},-v_{3})\quad \text{for }k\geq 1.
$$

On the other hand, we impose the far field boundary condition
\begin{equation}\label{4.16}
	\lim _{\eta \rightarrow \infty} \hat{f}_{k, 2}^{\pm}(t, x_{\sp}, \eta, v)=0.
\end{equation}

From Lemma \ref{lem4.2}, in order to solve \eqref{4.9}, \eqref{4.14}-\eqref{4.16}, we need $\hat{g}_{k}$ to satisfy \eqref{4.11}, i.e.,
\begin{equation}\label{4.20}
		\int_{\mathbb{R}^{3}} v_{3} \hat{g}_{k}^{\pm} \sqrt{\mu} d v =\int_{\mathbb{R}^{3}}v_{1} v_{3} \hat{g}_{k}^{\pm} \sqrt{\mu} d v=\int_{\mathbb{R}^{3}}v_{2} v_{3} \hat{g}_{k}^{\pm} \sqrt{\mu} d v =\int_{\mathbb{R}^{3}} v_{3}|v|^{2} \hat{g}_{k}^{\pm} \sqrt{\mu} d v=0,
\end{equation}
which are equivalent the followings:
\begin{equation}\label{4.21}
	\begin{aligned}
		&\int_{\mathbb{R}^{3}} v_{3} \sqrt{\mu} \mathbf{P}\{f_{k}^{\pm}+\bar{f}_{k}^{\pm}(t, x_{\sp},0, v)\}d v=-(\hat{A}_{k}^{\pm}+5 \hat{C}_{k}^{\pm})(t, x_{\sp},0), \\
		&\int_{\mathbb{R}^{3}} v_{3}v_{i} \sqrt{\mu}(\mathbf{I}-\mathbf{P})\{f_{k}^{\pm}+\bar{f}_{k}^{\pm}(t, x_{\sp},0, v)\}d v=-\hat{B}_{k, i}^{\pm}(t, x_{\sp},0) \quad i=1,2, \\
		&\int_{\mathbb{R}^{3}} v_{3}(|v|^{2}-5) \sqrt{\mu}(\mathbf{I}-\mathbf{P})\{f_{k}^{\pm}+\bar{f}_{k}^{\pm}(t, x_{\sp},0, v) \}d v=-10 \hat{C}_{k}^{\pm}(t, x_{\sp},0),
	\end{aligned}\quad \text{for }k\geq 1,
\end{equation}
where we have used \eqref{4.5}.

\begin{lemma}\label{lem4.3}
	Let $f_{k}^{\pm},\,\bar{f}_{k}^{\pm}$ and $\hat{f}_{k}^{\pm}$ satisfy \eqref{4.21}. Then it holds that
	\begin{align}
		&u_{k,3}^{\pm}+\bar{u}_{k,3}^{\pm,0}=-(\hat{A}_{k}^{\pm}+5 \hat{C}_{k}^{\pm})(t, x_{\sp}, 0) \equiv 0,\qquad\quad\,\,\,1\leq k\leq n,\label{4.33}\\
		&u_{n+m,3}^{\pm}+\bar{u}_{n+m,3}^{\pm,0}=-(\hat{A}_{n+m}^{\pm}+5 \hat{C}_{n+m}^{\pm})(t, x_{\sp}, 0) \quad 1 \leq m,\label{4.34}
	\end{align}
and
\begin{align}
	\left\{\begin{aligned}
		&\begin{aligned}
			\mp \partial_{y} \bar{u}_{m, i}^{\pm}(t, x_{\sp}, 0)
			&=\frac{1}{\lambda}\Big\{[(u_{1 . i}^{\pm}+\bar{u}_{1, i}^{\pm,0}) \bar{u}_{m,3}^{\pm,0}]+\hat{B}_{n-1+m, i}^{\pm,0}-u_{0, i}^{\pm}(\hat{A}_{m+1}^{\pm,0}+\hat{C}_{m+1}^{\pm,0})\\
			&\qquad\quad +\la\mathcal{A}_{3 i}, G_{m}^{\pm}\ra+\la\mathcal{A}_{3 i}, J_{m-1}^{\pm,0}\ra\Big\} \quad i=1,2,\quad \text{for }m \geq 1,
		\end{aligned}\\
	&\begin{aligned}
		\mp \partial_{y} \bar{\theta}_{m}^{\pm}\left(t, x_{\sp}, 0\right)
			&=\frac{1}{2\kappa}\Big\{5[\theta_{1}^{\pm}+\bar{\theta}_{1}^{\pm,0}] \bar{u}_{m,3}^{\pm,0}+5 \hat{C}_{n-1+m}^{\pm,0}-5\theta_{0}^{\pm}(\hat{A}_{m+1}^{\pm,0}+\hat{C}_{m+1}^{\pm,0})\\
			&\qquad\quad +2\la\mathcal{B}_{3}, G_{m}^{\pm}\ra+2\la\mathcal{B}_{3}, J_{m-1}^{\pm,0}\ra\Big\} \quad \text{for }m \geq 1,
		\end{aligned}
		\end{aligned}\right.\label{4.37}
	\end{align}
	where we have denoted $(\bar{u}_{k}^{\pm,0},\bar{\t}_{k}^{\pm,0}):=(\bar{u}_{k}^{\pm},\bar{\t}_{k}^{\pm})(t,x_{\sp},0)$, $(A_{k}^{\pm,0},B_{k}^{\pm,0},C_{k}^{\pm,0}):=(A_{k}^{\pm},B_{k}^{\pm},C_{k}^{\pm})(t,x_{\sp},0)$, and $G_{m}, J_{m-1}^{\pm,0}$ are the ones defined in \eqref{2.13} and \eqref{3.13}, respectively.
\end{lemma}

\noindent\textbf{Proof.} Firstly, we consider $\eqref{4.21}_{1}$. By using \eqref{4.7-1} and a direct calculation, we get \eqref{4.33}-\eqref{4.34}. 

Next, we consider $\eqref{4.21}_{2}$ and $\eqref{4.21}_{3}$. If $1 \leq k \leq n-3$, it follows from \eqref{2.1} and \eqref{3.1} that $f_{k}=\mathbf{P} f_{k}, \bar{f}_{k}^{\pm}=\mathbf{P} \bar{f}_{k}^{\pm}$. Thus $\eqref{4.21}_{2}$ and $\eqref{4.21}_{3}$ are satisfied automatically.

If $k=n-2$ or $n-1$, by using \eqref{2.0} and \eqref{3.1}, then $\eqref{4.21}_2-\eqref{4.21}_3$ can be reduced as
\begin{equation}\nonumber
	\begin{aligned}
	&u_{0, i}^{\pm} u_{0,3}^{\pm}=-\hat{B}_{n-2,i}^{\pm}(t, x_{\sp}, 0)\equiv 0,\quad i=1,2, \quad 5\theta_{0}^{\pm} u_{0,3}^{\pm}=-10 \hat{C}_{n-2}^{\pm}(t, x_{\sp}, 0)\equiv 0,\\
	&u_{0, i}^{\pm} u_{1,3}^{\pm}=-\hat{B}_{n-1, i}^{\pm}(t, x_{\sp}, 0)\equiv0, \quad i=1,2, \quad 5\theta_{0}^{\pm} u_{1,3}^{\pm}=-10 \hat{C}_{n-1}^{\pm}(t, x_{\sp}, 0)\equiv 0,
	\end{aligned}
\end{equation}
which holds automatically by \eqref{4.7-1}, \eqref{4.33} and the fact $u_{0,3}^{\pm}=\bar{u}_{1,3}^{\pm}\equiv 0$.

If $k \geq n$, denoting $k=n-1+m$ with $m \geq 1$, it follows from \eqref{2.7-1}-\eqref{2.7-4}, \eqref{3.16-2}-\eqref{3.16-3}, \eqref{4.33}-\eqref{4.34} and direct calculations that
\begin{equation*}
	\begin{aligned}
		\mp\partial_{y} \bar{u}_{m, i}^{\pm}\left(t, x_{\sp}, 0\right)&=\frac{1}{\lambda}\Big\{[(u_{1,i}^{\pm}+\bar{u}_{1, i}^{\pm,0}) \bar{u}_{m, 3}^{\pm,0}]+u_{0, i}^{\pm} \bar{u}_{1+m, 3}^{\pm,0}+\la\mathcal{A}_{3 i}, J_{m-1}^{\pm,0}\ra\\
		&\quad+\la\mathcal{A}_{3 i},(\mathbf{I}-\mathbf{P}) f_{n-1+m}^{\pm}\ra+\hat{B}_{n-1+m, i}^{\pm,0}\Big\},\\
		&=\frac{1}{\lambda}\Big\{[(u_{1,i}^{\pm}+\bar{u}_{1, i}^{\pm,0}) \bar{u}_{m, 3}^{\pm,0}]+\la\mathcal{A}_{3 i}, J_{m-1}^{\pm,0}\ra-u_{0, i}^{\pm}(\hat{A}_{m+1}^{\pm,0}+\hat{C}_{m+1}^{\pm,0})\\
		&\quad+\la\mathcal{A}_{3 i}, G_{m}^{\pm}\ra+\hat{B}_{n-1+m, i}^{\pm,0}\Big\} \quad i=1,2, \text{for }m \geq 0,
	\end{aligned}
\end{equation*}
\begin{equation*}
	\begin{aligned}
		\mp \partial_{y} \bar{\theta}_{m}^{\pm}(t, x_{\sp}, 0)&=\frac{1}{2\kappa}\Big\{5[\theta_{1}^{\pm}+\bar{\theta}_{1}^{\pm,0}] \bar{u}_{m, 3}^{\pm,0}+5 \theta_{0}^{\pm} \bar{u}_{m+1, 3}^{\pm,0}+2\la\mathcal{B}_{3}, J_{m-1}^{\pm,0}\ra\\
		&\qquad +2\la\mathcal{B}_{3}, (\mathbf{I}-\mathbf{P})f_{n-1+m}^{\pm}\ra+10\hat{C}_{n-1+m}^{\pm,0}\Big\}\\
		&=\frac{1}{2\kappa}\Big\{5[\theta_{1}^{\pm}+\bar{\theta}_{1}^{\pm,0}] \bar{u}_{m, 3}^{\pm,0}+2\la\mathcal{B}_{3}, J_{m-1}^{\pm,0}\ra-5\theta_{0}^{\pm}(\hat{A}_{m+1}^{\pm,0}+\hat{C}_{m+1}^{\pm,0})\\
		&\qquad +2\la\mathcal{B}_{3}, G_{m}^{\pm}\ra+10 \hat{C}_{n-1+m}^{\pm,0}\Big\} \quad \text{for }m \geq 0,
	\end{aligned}
\end{equation*}
where we have used $u_{0,3}^{\pm}=\bar{u}_{1,3}^{\pm}=0$ and $G_{m}, J_{m-1}^{\pm}$ are defined in \eqref{2.13} and \eqref{3.13} respectively. Therefore the proof Lemma \ref{lem4.3} is complete. $\hfill\square$

\begin{remark}\label{rem4.3}
1) Taking $k=1$ in \eqref{4.33} and using \eqref{3.3}, we get $u_{1,3}^{\pm}=-\bar{u}_{1,3}^{\pm,0}=0$. It is clear that $\hat{A}_{n+m}^{\pm}, \hat{C}_{n+m}^{\pm}$ on RHS of \eqref{4.34} depend on $\hat{f}_{m}^{\pm}$ and $\hat{f}_{m+2-n}^{\pm}$.

\smallskip
	
2) We point out that the terms on the left hand side of \eqref{4.37} depend only on $f_{0}, f_{m}, f_{k}$ and $\bar{f}_{k}^{\pm}$ with $1\leq k\leq m-1$. In fact $G_{m}$ depends only on $f_{k}$ with $1\leq k\leq m$ from \eqref{2.13}, $\bar{u}_{m,3}^{\pm}$ depends only on $\bar{f}_{m-1}^{\pm}$, $\bar{f}_{m+1-n}^{\pm}$ from \eqref{3.21}, and $\la \mathcal{A}_{3i}, J_{m-1}^{\pm,0}\ra$, $\la \mathcal{B}_{3}, J_{m-1}^{\pm,0}\ra$ depends only on $f_{0}, f_{m}, f_{k}$ and $\bar{f}_{k}^{\pm}$ with $1\leq k\leq m-1$ as explained in the last two paragraphs in the proof of Proposition \ref{prop3.1}.
\end{remark}

\begin{remark}
Taking $m=1$ in \eqref{4.37} and using \eqref{3.3}, \eqref{3.20-1} and \eqref{4.7-1}, we obtain
\begin{equation}\label{4.35}
\left\{\begin{aligned}
&\mp \pa_{y}\bar{u}_{1,i}^{\pm}(t,x_{\sp},0)=\frac{1}{\lambda}\la
\mathcal{A}_{3,i}, G_{1}^{\pm}\ra,\quad i=1,2,\\
&\mp \pa_{y}\bar{\t}_{1}^{\pm}(t,x_{\sp},0)=\frac{1}{\kappa}\la
\mathcal{B}_{3}, G_{1}^{\pm}\ra.
\end{aligned}
\right.
\end{equation}
It follows from \eqref{2.13} that
\begin{equation}\label{4.35-1}
\begin{aligned}
G_{1}=\mathbf{L}^{-1}\big\{&-(\mathbf{I}-\mathbf{P})v\cdot \nabla_{x}f_{0}+\Gamma(\mathbf{P}f_{1},\mathbf{P}f_{1})+\Gamma(\mathbf{P}f_{1},(\mathbf{I}-\mathbf{P})f_{1})+\Gamma((\mathbf{I}-\mathbf{P})f_{1},\mathbf{P}f_{1})\\
&+\Gamma((\mathbf{I}-\mathbf{P})f_{1},(\mathbf{I}-\mathbf{P})f_{1})+\Gamma(f_{0},(\mathbf{I}-\mathbf{P})f_{2})+\Gamma((\mathbf{I}-\mathbf{P})f_{2},f_{0})\big\}
\end{aligned}
\end{equation}
Noting $u_{0,3}^{\pm}=u_{1,3}^{\pm}\equiv 0$, and using \eqref{2.0}, \eqref{4.35-1} and similar arguments as in \eqref{3.20-1}, one can deduce
\begin{equation}\label{4.35-2}
\la \mathcal{A}_{3i}, G_{1}^{\pm}\ra=-\lambda \pa_{x_{3}}u_{0,i}^{\pm},\quad  \la \mathcal{B}_{3},G_{1}^{\pm}\ra =-\kappa \pa_{x_{3}}\theta_{0}^{\pm},\quad i=1,2.	
\end{equation}
Substituting \eqref{4.35-2} into \eqref{4.35}, we obtain
\begin{equation}\label{4.35-3}
	\left\{\begin{aligned}
		& \pa_{y}\bar{u}_{1,i}^{\pm}(t,x_{\sp},0)=\pm \pa_{x_{3}}u_{0,i}^{\pm},\quad i=1,2,\\
		& \pa_{y}\bar{\t}_{1}^{\pm}(t,x_{\sp},0)=\pm \pa_{x_{3}}\theta_{0}^{\pm}.
	\end{aligned}
	\right.
\end{equation}
\end{remark}

\begin{remark}\label{rem4.2}
	From \eqref{4.3}, we have
	\begin{equation}\nonumber
		\hat{S}_{i,1}^{\pm}=\hat{S}_{i,2}^{\pm}=\hat{f}_{i,1}^{\pm} \equiv 0, \quad \hat{f}_{n-1,1}^{\pm}=\hat{f}_{n,1}^{\pm}\equiv 0\quad \text {for } \quad 1 \leq i \leq n-2,
	\end{equation}
	Moreover, since $f_0=\mathbf{P}f_{0}$ and  $f_k=\mathbf{P}f_{k},\bar{f}_{k}^{\pm}=\mathbf{P}\bar{f}_{k}^{\pm}$ for $1\leq k\leq n-3$, we obtain from \eqref{4.33} that
	\begin{equation}\nonumber
		\begin{aligned}
			\hat{g}_{k}^{\pm}(t,x_{\sp},v_{\sp},v_3)&=[f_{k}^{\pm}(v)+(\bar{f}_{k}^{\pm}+\hat{f}_{k,1}^{\pm})(t,x_{\sp},0,v_{\sp},v_{3})]\\
			&\qquad -[f_{k}^{\pm}(R_{x}v)+(\bar{f}_{k,\pm}+\hat{f}_{k,\pm,1})(t,x_{\sp},0,v_{\sp},-v_{3})]\\
			&=2(u_{k,3}^{\pm}+\bar{u}_{k,3}^{\pm,0})v_{3}\equiv 0,\quad \text{for }\mp v_{3}<0,\quad 1\leq k\leq n-3.
		\end{aligned}
	\end{equation}
	Additionally, noting that $\bar{f}_{n-2}^{\pm}=\mathbf{P}\bar{f}_{n-2}^{\pm}$ and
	$$
	\begin{aligned}
	\mathbf{(I-P)}f_{n-2}&=\frac{1}{2}\mathbf{(I-P)}\{\frac{f_0f_0}{\sqrt{\mu}}\}\\
	&=\frac{1}{2}\sum\limits_{l,m=1}^3u_{0,\ell}u_{0,m}\mathcal{A}_{l m}+\sum\limits_{m=1}^3u_{0,m}\t_0\mathcal{B}_{m}+\frac{\t_0^2}{8}\mathbf{(I-P)}\{(|v|^2-5)^2\sqrt{\mu}\},
	\end{aligned}
	$$
	then from $u_{0,3}^{\pm}=0$, we have
	\begin{equation}\nonumber
		\begin{aligned}
			\hat{g}_{n-2}^{\pm}(t,x_{\sp},v_{\sp},v_3)&=[f_{n-2}^{\pm}(v)+(\bar{f}_{n-2}^{\pm}+\hat{f}_{n-2,1}^{\pm})(t,x_{\sp},0,v_{\sp},v_{3})]\\
			&\qquad -[f_{n-2}^{\pm}(R_{x}v)+(\bar{f}_{n-2}^{\pm}+\hat{f}_{n-2,1}^{\pm})(t,x_{\sp},0,v_{\sp},-v_{3})]\\
			&=2(u_{n-2,3}^{\pm}+\bar{u}_{n-2,3}^{\pm,0})v_{3}+\mathbf{(I-P)}f_{n-2}^{\pm}(v)-\mathbf{(I-P)}f_{n-2}^{\pm}(R_{x}v)\\
			&\equiv 0,\quad \text{for }\mp v_{3}<0.
		\end{aligned}
	\end{equation}
	then it follows from Lemma \ref{lem4.2}, $\eqref{4.3}_1$ and \eqref{4.15}-\eqref{4.16} that
	\begin{equation}\nonumber
		\hat{f}_{i,2}^{\pm} \equiv 0, \quad \text { for } \quad 1 \leq i \leq n-2.
	\end{equation}
	which implies $\hat{F}_{i}^{\pm} \equiv 0$ for $1 \leq i \leq n-2$. Furthermore, it follows from \eqref{4.3}-\eqref{4.3-1} that
	\begin{equation}\nonumber
		\hat{S}_{i,1}^{\pm}=\hat{f}_{i,1}^{\pm} \equiv 0, \quad \text { for }1 \leq i \leq 2 n-2,\qquad \hat{S}_{i,2}^{\pm}\equiv 0,\quad \text{for }1\leq i\leq 2n-4.
	\end{equation}
Hence, it also holds that $\hat{A}_{i}^{\pm}=\hat{B}_{i}^{\pm}=\hat{C}_{i}^{\pm}\equiv 0$ for $1\leq i\leq 2n-2$.
\end{remark}

\subsection{Criteria of uniqueness for  pressure terms}
It follows from \eqref{2.17} that $\nabla_{x}p$ is uniquely determined up to a function of $t$. However, we find that the pressure terms can be uniquely determined by the compatibility conditions for the equations of $u_{k}$. Indeed,  in view of $\eqref{2.17}_{1}$ and \eqref{4.33}-\eqref{4.34}, we have the following restrictions on pressure terms.
\begin{lemma}\label{lem2.2}
	To guarantee the compatibility conditions for equations of velocity in \eqref{2.17}, we  require the pressure terms $p_{k}$ in \eqref{2.17-1} to satisfy the followings:
    \begin{itemize}
    	\item [(1)] For $k=0$, $p_{0}(t)\equiv p_{0}(0)=(\rho_{0}+\theta_{0})(0)$. It is worth mentioning that $(\rho_{0}+\theta_{0})(0)$ is a constant since $\nabla_{x}(\rho_{0}+\theta_{0})\equiv 0$.
    	\item [(2)] For $1\leq k\leq n-3$ (if $n=3$, then skip this case), we impose
    	\begin{equation}\label{4.80}
    		\begin{aligned}
    		p_{k}(t)&=(\rho_{k}+\theta_{k})(0)-\frac{2}{3}\int_{0}^{t}\int_{\mathbb{T}^2}(\langle \mathcal{B}_{3}, G_{k-1}^{+}\rangle-\langle \mathcal{B}_{3}, G_{k-1}^{-}\rangle)(s,x_{\sp})dx_{\sp}ds\\
    			&\quad -\frac{5}{3}\int_{\mathbb{T}^2}\int_{0}^{\infty}[(\bar{\rho}_{k-1}^{+}+\bar{\rho}_{k-1}^{-})(t,x_{\sp},y)-\big(\bar{\rho}_{k-1}^{+}+\bar{\rho}_{k-1}^{-}\big)(0,x_{\sp},y)]dydx_{\sp}\\
    			&\quad -\frac{5}{3}\int_{0}^{t}\int_{\mathbb{T}^2}\int_{0}^{+\infty}[\operatorname{div}_{\sp}\bar{u}_{k-1,\sp}^{+}\theta_{0}^{+}+\operatorname{div}_{\sp}\bar{u}_{k-1,\sp}^{-}\theta_{0}^{-}]dydx_{\sp}ds.
    		\end{aligned}
    	\end{equation}
    \item [(3)] For $k=n-2$, we impose
    \begin{equation}\label{4.81}
    	\begin{aligned}
    		\int_{\Omega}p_{n-2}(t,x)dx&=\int_{\Omega}(\rho_{n-2}+\theta_{n-2})(0,x)dx-\frac{1}{3}\int_{\Omega}|u_{0}(t,x)|^2dx\\
    		&\quad -\frac{2}{3}\int_{0}^{t}\int_{\mathbb{T}^2}\big(\langle \mathcal{B}_{3}, G_{n-3}^{+}\rangle-\langle \mathcal{B}_{3}, G_{n-3}^{-}\rangle\big)(s,x_{\sp})dx_{\sp}ds\\
    		&\quad -\frac{5}{3}\int_{\mathbb{T}^2}\int_{0}^{\infty}[(\bar{\rho}_{n-3}^{+}+\bar{\rho}_{n-3}^{-})(t,x_{\sp},y)-(\bar{\rho}_{n-3}^{+}+\bar{\rho}_{n-3}^{-})(0,x_{\sp},y)]dydx_{\sp}\\
    		&\quad
    		-\frac{5}{3}\int_{0}^{t}\int_{\mathbb{T}^2}\int_{0}^{\infty}(\operatorname{div}_{\sp}\bar{u}_{n-3,\sp}^{+}\theta_{0}^{+}+\operatorname{div}_{\sp}\bar{u}_{n-3,\sp}^{-}\theta_{0}^{-})dydx_{\sp}ds\\
    		&:=\mathfrak{c}_{n-2}(t).
    	\end{aligned}
    \end{equation}
   \item [(4)] For $k\geq n-1$, denoting $k=n-2+m$ with $m\geq 1$, then we impose
   \begin{equation}\label{4.82}
   	\begin{aligned}
   		&\int_{\Omega}p_{n-2+m}(t,x)dx=\int_{\Omega}(\rho_{n-2+m}+\theta_{n-2+m})(0,x)dx-\frac{2}{3}\int_{\Omega}(u_{0}\cdot u_{m})(t,x)dx\\
   		&\quad -\frac{2}{3}\int_{0}^{t}\int_{\mathbb{T}^2}(\langle \mathcal{B}_{3}, G_{n-3+m}^{+}\rangle-\langle \mathcal{B}_{3}, G_{n-3+m}^{-}\rangle)(s,x_{\sp})dx_{\sp}ds\\
   		&\quad -\frac{5}{3}\int_{\mathbb{T}^2}\int_{0}^{\infty}[(\bar{\rho}_{n-3+m}^{+}+\bar{\rho}_{n-3+m}^{-})(t,x_{\sp},y)-(\bar{\rho}_{n-3+m}^{+}+\bar{\rho}_{n-3+m}^{-})(0,x_{\sp},y)]dydx_{\sp}\\
   		&\quad
   		-\frac{5}{3}\int_{0}^{t}\int_{\mathbb{T}^2}\int_{0}^{\infty}[\operatorname{div}_{\sp}\bar{u}_{n-3+m,\sp}^{+}\theta_{0}^{+}+\operatorname{div}_{\sp}\bar{u}_{n-3+m,\sp}^{-}\bar{\theta}_{0}^{-}]dydx_{\sp}ds\\
   		&\quad
   		-\frac{5}{3}\int_{0}^{t}\int_{\mathbb{T}^2}\int_{0}^{\infty}[\pa_{t}\bar{\rho}_{m-1}^{+}\theta_{0}^{+}+\pa_{t}\bar{\rho}_{m-1}^{-}\bar{\theta}_{0}^{-}]dydx_{\sp}ds\\
   		&\quad +\frac{5}{3}\int_{0}^{t}\int_{\mathbb{T}^2}[(\hat{A}_{2n-4+m}^{+,0}+5\hat{C}_{2n-4+m}^{+,0})-(\hat{A}_{2n-4+m}^{-,0}+5\hat{C}_{2n-4+m}^{-,0})](s,x_{\sp})dx_{\sp}ds\\
   		&\quad
   		+\frac{5}{3}\int_{0}^{t}\int_{\mathbb{T}^2}[(\hat{A}_{n-2+m}^{+,0}+5\hat{C}_{n-2+m}^{+,0})\theta_{0}^{+}-(\hat{A}_{n-2+m}^{-,0}+5\hat{C}_{n-2+m}^{-,0})\theta_{0}^{-}](s,x_{\sp})dx_{\sp}ds\\
   		&:=\mathfrak{c}_{n-2+m}(t),
   	\end{aligned}	
   \end{equation}
    \end{itemize}
where $G_{k}$ is the one defined in \eqref{2.13}.
\end{lemma}

\noindent\textbf{Proof. } From $\eqref{2.17}_{1}$ and \eqref{4.33}-\eqref{4.34}, one has
\begin{equation}\label{4.65}
	\begin{aligned}
		-\int_{\Omega}\pa_{t}\rho_{k+2-n}dx&=\int_{\Omega}\operatorname{div}_{x}u_{k}dx=\int_{\mathbb{T}^2}(u_{k,3}^{+}-u_{k,3}^{-})dx_{\sp}\\
		&=\int_{\mathbb{T}^2}[-\bar{u}_{k,3}^{+,0}-(\hat{A}_{k}^{+,0}+5\hat{C}_{k}^{+,0})+\bar{u}_{k,3}^{-,0}+(\hat{A}_{k}^{-,0}+5\hat{C}_{k}^{-,0})]dx_{\sp},
	\end{aligned}
\end{equation}
 which, together with \eqref{3.3-3} and \eqref{3.21}, yields
\begin{equation}\label{4.66}
	\begin{aligned}
		&-\int_{\Omega}\pa_{t}\rho_{k+2-n}dx=\int_{\mathbb{T}^2}(\hat{A}_{k}^{-,0}+5\hat{C}_{k}^{-,0})-(\hat{A}_{k}^{+,0}+5\hat{C}_{k}^{+,0})dx_{\sp}\\
		&\qquad+\int_{\mathbb{T}^2}\int_{0}^{\infty}[(\pa_{t}\bar{\rho}_{k+1-n}^{-}+\operatorname{div}_{\sp}\bar{u}_{k-1,\sp}^{-})+(\pa_{t}\bar{\rho}_{k+1-n}^{+}+\operatorname{div}_{\sp}\bar{u}_{k-1,\sp}^{+})]dydx_{\sp}\\
		&=\int_{\mathbb{T}^2}(\hat{A}_{k}^{-,0}+5\hat{C}_{k}^{-,0})-(\hat{A}_{k}^{+,0}+5\hat{C}_{k}^{+,0})dx_{\sp}+\int_{\mathbb{T}^2}\int_{0}^{\infty}(\pa_{t}\bar{\rho}_{k+1-n}^{-}+\pa_{t}\bar{\rho}_{k+1-n}^{+})dydx_{\sp}.
	\end{aligned}
\end{equation}
Now, it follows from \eqref{2.17-1} and \eqref{4.66} that
\begin{align}
		&\pa_{t}\int_{\Omega}p_{k}dx=\pa_{t}\int_{\Omega}\theta_{k}dx+\pa_{t}\int_{\Omega}\rho_{k}dx\nonumber\\
		&=\pa_{t}\int_{\Omega}\theta_{k}dx-\int_{\mathbb{T}^2}\int_{0}^{\infty}\big(\pa_{t}\bar{\rho}_{k-1}^{-}+\pa_{t}\bar{\rho}_{k-1}^{+}\big)dydx_{\sp}\nonumber\\
		&\quad +\int_{\mathbb{T}^2}[(\hat{A}_{k+n-2}^{+,0}+5\hat{C}_{k+n-2}^{+,0})-(\hat{A}_{k+n-2}^{-,0}+5\hat{C}_{k+n-2}^{-,0})]dx_{\sp},\quad \text{for }0\leq k\leq n-3,\label{4.67}\\
		&\pa_{t}\int_{\Omega}p_{n-2}dx=\pa_{t}\int_{\Omega}\theta_{n-2}dx+\pa_{t}\int_{\Omega}\rho_{n-2}dx-\frac{1}{3}\pa_{t}\int_{\Omega}|u_{0}|^2dx\nonumber\\
		&=\pa_{t}\int_{\Omega}\theta_{n-2}dx-\frac{1}{3}\pa_{t}\int_{\Omega}|u_{0}|^2dx-\int_{\mathbb{T}^2}\int_{0}^{\infty}(\pa_{t}\bar{\rho}_{n-3}^{-}+\pa_{t}\bar{\rho}_{n-3}^{+})dydx_{\sp}\nonumber\\
		&\quad +\int_{\mathbb{T}^2}[(\hat{A}_{2n-4}^{+,0}+5\hat{C}_{2n-4}^{+,0})-(\hat{A}_{2n-4}^{-,0}+5\hat{C}_{2n-4}^{-,0})]dx_{\sp},\label{4.68}\\
		&\pa_{t}\int_{\Omega}p_{n-2+m}dx=\pa_{t}\int_{\Omega}\theta_{n-2+m}dx+\pa_{t}\int_{\Omega}\rho_{n-2+m}dx-\frac{2}{3}\pa_{t}\int_{\Omega}u_{0}\cdot u_{m}dx\nonumber\\
		&=\pa_{t}\int_{\Omega}\theta_{n-2+m}dx-\frac{2}{3}\pa_{t}\int_{\Omega}u_{0}\cdot u_{m}dx-\int_{\mathbb{T}^2}\int_{0}^{\infty}(\pa_{t}\bar{\rho}_{n-3+m}^{-}+\pa_{t}\bar{\rho}_{n-3+m}^{+})dydx_{\sp}\nonumber\\
		&\quad-\int_{\mathbb{T}^2}[(\hat{A}_{2n-4+m}^{-,0}+5\hat{C}_{2n-4+m}^{-,0})-(\hat{A}_{2n-4+m}^{+,0}+5\hat{C}_{2n-4+m}^{+,0})]dx_{\sp},\quad \text{for }m\geq 1.\label{4.69}
\end{align}

We divide the rest calculations into three cases.

{\it Case 1.} For the case of $k=0$, it follows from $\eqref{2.8-1}_{3}$ that
\begin{equation*}
	\pa_{t}\int_{\Omega}\theta_{0}dx+\int_{\Omega}(u_{0}\cdot \nabla_{x})\theta_{0}dx=\frac{2}{5}\pa_{t}\int_{\Omega}p_{0}(t)dx,
\end{equation*}
which, together with $\operatorname{div}_{x}u_{0}=0$, \eqref{2.8-2} and integrating by parts, yields that
\begin{equation}\label{4.75}
	\pa_{t}\int_{\Omega}\theta_{0}dx=\frac{2}{5}\pa_{t}\int_{\Omega}p_{0}(t)dx.
\end{equation}
Substituting \eqref{4.75} into \eqref{4.67} and noting Remark \ref{rem4.2}, one has
\begin{equation*}
	\pa_{t}\int_{\Omega}p_{0}(t)dx=0,
\end{equation*}
which, together with \eqref{2.4}-\eqref{2.5-1}, implies that $p_{0}(t)\equiv p_{0}(0)=(\rho_{0}+\theta_{0})(0)$ is a constant, and depends only on the initial data $(\rho_{0}(0),\theta_{0}(0))$.

{\it Case 2.} For the case of $1\leq k\leq n-3$, it follows from $\eqref{2.17}_{3}$ and $\operatorname{div}_{x}u_{0}=\operatorname{div}_{x}u_{k}=0$ that
\begin{equation}\label{4.76}
	\begin{aligned}
		&\pa_{t}\int_{\Omega}\theta_{k}dx+\int_{\mathbb{T}^2}\big(u_{k,3}^{+}\theta_{0}^{+}-u_{k,3}^{-}\theta_{0}^{-}\big)dx_{\sp}
		\\&=\int_{\Omega}\mathfrak{q}_{k-1}dx=\frac{2}{5}\pa_{t}\int_{\Omega}p_{k}dx-\frac{2}{5}\int_{\mathbb{T}^2}\big(\langle \mathcal{B}_{3}, G_{k-1}^{+}\rangle-\langle \mathcal{B}_{3}, G_{k-1}^{-}\rangle\big) dx_{\sp},
	\end{aligned}
\end{equation}
where we have used \eqref{2.17-2} in the last equality. Using \eqref{4.33}-\eqref{4.34} and \eqref{3.21}, one has
\begin{equation}\label{4.77}
	\begin{aligned}
		&\int_{\mathbb{T}^2}(u_{k,3}^{+}\theta_{0}^{+}-u_{k,3}^{-}\theta_{0}^{-})dx_{\sp}\\
		&=\int_{\mathbb{T}^2}(-\bar{u}_{k,3}^{+}\theta_{0}^{+}+\bar{u}_{k,3}^{-}\theta_{0}^{-})dx_{\sp} +\int_{\mathbb{T}^2}[-(\hat{A}_{k}^{+,0}+5\hat{C}_{k}^{+,0})\theta_{0}^{+}+(\hat{A}_{k}^{-,0}+5\hat{C}_{k}^{-,0})\theta_{0}^{-}]dx_{\sp}\\
		&=\int_{\mathbb{T}^2}\int_{0}^{+\infty}[(\pa_{t}\bar{\rho}_{k+1-n}^{+}+\operatorname{div}_{\sp}\bar{u}_{k-1,\sp}^{+})\theta_{0}^{+}+(\pa_{t}
		\bar{\rho}_{k+1-n}^{-}+\operatorname{div}_{\sp}\bar{u}_{k-1,\sp}^{-})\theta_{0}^{-}]dydx_{\sp}\\
		&=\int_{\mathbb{T}^2}\int_{0}^{+\infty}[\operatorname{div}_{\sp}\bar{u}_{k-1,\sp}^{+}\theta_{0}^{+}+\operatorname{div}_{\sp}\bar{u}_{k-1,\sp}^{-}\theta_{0}^{-}]dydx_{\sp},
	\end{aligned}
\end{equation}
where we have used $\hat{A}_{k}^{\pm}=\hat{B}_{k}^{\pm}=\hat{C}_{\pm}^{k}\equiv 0$ for $1\leq k\leq 2n-2$ (see Remark \ref{rem4.2}) and the fact $\bar{\rho}_{j+1-n}^{\pm}\equiv 0$ for $1\leq j\leq n-1$. Substituting \eqref{4.77} into \eqref{4.76}, we obtain
\begin{equation}\label{4.78}
	\begin{aligned}
		\pa_{t}\int_{\Omega}\theta_{k}dx&=\frac{2}{5}\pa_{t}\int_{\Omega}p_{k}dx-\frac{2}{5}\int_{\mathbb{T}^2}\big(\langle \mathcal{B}_{3}, G_{k-1}^{+}\rangle-\langle \mathcal{B}_{3}, G_{k-1}^{-}\rangle\big) dx_{\sp}\\
		&\quad  -\int_{\mathbb{T}^2}\int_{0}^{+\infty}[\operatorname{div}_{\sp}\bar{u}_{k-1,\sp}^{+}\theta_{0}^{+}+\operatorname{div}_{\sp}\bar{u}_{k-1,\sp}^{-}\theta_{0}^{-}]dydx_{\sp},\\
	\end{aligned}
\end{equation}
which, together with \eqref{4.67}, implies that
\begin{equation}\label{4.79}
	\begin{aligned}
		\pa_{t}\int_{\Omega}p_{k}dx&=-\frac{2}{3}\int_{\mathbb{T}^2}\big(\langle \mathcal{B}_{3}, G_{k-1}^{+}\rangle-\langle \mathcal{B}_{3}, G_{k-1}^{-}\rangle\big) dx_{\sp}-\frac{5}{3}\int_{\mathbb{T}^2}\int_{0}^{\infty}(\pa_{t}\bar{\rho}_{k-1}^{+}+\pa_{t}\bar{\rho}_{k-1}^{-})dydx_{\sp}\\
		&\quad -\frac{5}{3}\int_{\mathbb{T}^2}\int_{0}^{+\infty}[\operatorname{div}_{\sp}\bar{u}_{k-1,\sp}^{+}\theta_{0}^{+}+\operatorname{div}_{\sp}\bar{u}_{k-1,\sp}^{-}\theta_{0}^{-}]dydx_{\sp}\\
		&\quad  +\frac{5}{3}\int_{\mathbb{T}^2}[(\hat{A}_{k+n-2}^{+,0}+5\hat{C}_{k+n-2}^{+,0})-(\hat{A}_{k+n-2}^{-,0}+5\hat{C}_{k+n-2}^{-,0})]dx_{\sp}\\
	  &=-\frac{2}{3}\int_{\mathbb{T}^2}\big(\langle \mathcal{B}_{3}, G_{k-1}^{+}\rangle-\langle \mathcal{B}_{3}, G_{k-1}^{-}\rangle\big) dx_{\sp}-\frac{5}{3}\int_{\mathbb{T}^2}\int_{0}^{\infty}(\pa_{t}\bar{\rho}_{k-1}^{+}+\pa_{t}\bar{\rho}_{k-1}^{-})dydx_{\sp}\\
	  &\quad -\frac{5}{3}\int_{\mathbb{T}^2}\int_{0}^{+\infty}[\operatorname{div}_{\sp}\bar{u}_{k-1,\sp}^{+}\theta_{0}^{+}+\operatorname{div}_{\sp}\bar{u}_{k-1,\sp}^{-}\theta_{0}^{-}]dydx_{\sp},
	\end{aligned}
\end{equation}
where we have used $\hat{A}_{k}^{\pm}=\hat{B}_{k}^{\pm}=\hat{C}_{\pm}^{k}\equiv 0$ for $1\leq k\leq 2n-2$ (see Remark \ref{rem4.2}). Integrating \eqref{4.79} over $[0,t]$ and noting $\eqref{2.17-1}_{1}$, we conclude \eqref{4.80}.

{\it Case 3. } For the case of $k=n-2$, using \eqref{4.68} and similar calculations as in case 2, one gets \eqref{4.81}. For the case of $k\geq n-1$, we denote $k=n-2+m$ with $m\geq 1$. Using \eqref{4.69} and similar calculations as in case 2, we deduce \eqref{4.82}. Thus, the proof of Lemma \ref{lem2.2} is complete. $\hfill\square$

\begin{remark}\label{rem2.3}
	It follows from \eqref{4.80} that $p_{k}(t)$ ($1\leq k\leq n-3$) depends only on the initial data $(\rho_{k}+\theta_{k})(0)$, $f_{0}$, $f_{j}$ and $\bar{f}_{j}^{\pm}$ with $1\leq j\leq k-1$. Then it can be regarded as a known function when we solve the equations of $(\rho_{k},u_{k},\theta_{k})$.
\end{remark}

\begin{remark}\label{rem2.3-1}
	It follows from \eqref{4.81}-\eqref{4.82},  \eqref{4.3}-\eqref{4.4} and \eqref{4.6} that $\mathfrak{c}_{n-2+k}(t)$ ($k\geq 0$) depends only on
	the initial data $(\rho_{n-2+k}+\theta_{n-2+k})(0)$, $f_{0}$, $f_{j}$, $\bar{f}_{j}^{\pm}$ with $1\leq j\leq n-3+k$ and $\hat{f}_{i}$ with $1\leq i\leq n-4+k$. That is $\int_{\Omega} p_{n-2+k}dx$ can be determined once one solves $f_{0}$, $f_{j}$, $\bar{f}_{j}^{\pm}$ with $1\leq j\leq n-3+k$ and $\hat{f}_{i}$ with $1\leq i\leq n-4+k$. This will not cause any contradictions since $p_{n-2+k}$ plays no role when one solves $f_{j},\bar{f}_{j}^{\pm},\hat{f}_{j}^{\pm}$ for $0\leq j<n-2+k$. Moreover, the restriction of $\mathfrak{c}_{n-2+k}(t)$ provides a criterion of  uniqueness for $p_{n-2+k}$.
	
	On the other hand, it follows from \eqref{2.17} that $\nabla_{x}p_{n-2+k}$ can be determined when we solve the equation of $(\rho_{k},u_{k},\theta_{k})$. Once $\nabla_{x}p_{n-2+k}$ is obtained, $p_{n-2+k}$ can be determined uniquely until we give the initial data of $(\rho_{n-2+k},u_{n-2+k},\theta_{n-2+k})$, see \eqref{4.82}.
	
	Finally, we point out that one must involve the pressure term  $p_{k+n-2}$  when we solve $(\rho_{k},u_{k},\theta_{k})$ with $N+3-n\leq k\leq N$ (see \eqref{2.17}), since we use the truncated Hilbert expansion \eqref{1.22-2} in the present paper, then the estimate of $\nabla_x p_{k+n-2}$ is enough and there is no need to determine $p_{k+n-2}$ for $(\rho_{k},u_{k},\theta_{k})$ with $N+3-n\leq k\leq N$.
\end{remark}

\begin{remark}\label{rem2.4}
	If $n=3$, it follows from $\eqref{4.81}$ and $G_{0}\equiv 0$ that
	$$
	\begin{aligned}
	\mathfrak{c}_{1}(t)=\int_{\Omega}p_{1}(t,x)dx
	&=\int_{\Omega}(\rho_{1}+\theta_{1})(0,x)dx-\frac{1}{3}\int_{\Omega}|u_{0}(0,x)|^2dx,
	\end{aligned}
	$$
	where we have used the fact that $\int_{\Omega}|u_{0}(t,x)|^2dx=\int_{\Omega}|u_{0}(0,x)|^2dx$. Hence, $\mathfrak{c}_{1}(t)$ is a constant.
\end{remark}

\

\section{Existence of solutions for linear systems}

\subsection{Existence of solutions for a kind of linear Euler equations}\label{LE}
In view of Lemma \ref{lem2.2} and Remarks \ref{rem2.3}-\ref{rem2.3-1}, to solve \eqref{2.17}, we consider the following linear problem for $(\tilde{\rho},\tilde{u},\tilde{\theta})(t,x)$
\begin{equation}\label{4.40}
	\left\{\begin{aligned}
		&\mathrm{div}_{x}\tilde{u}=\mathfrak{r},\quad \tilde{\rho}+\tilde{\theta}=\mathfrak{p},\\
		&\partial_{t}\tilde{u}+(u_{0}\cdot \nabla_{x})\tilde{u}+(\tilde{u}\cdot \nabla_{x})u_{0}+\nabla_{x} \tilde{p}=\mathfrak{h},\\
		&\partial_{t}\tilde{\theta}+(u_{0}\cdot \nabla_{x})\tilde{\theta}+(\tilde{u}\cdot \nabla_{x})\theta_0=\mathfrak{q},\\
		&\int_{\Omega}\tilde{p}(t,x)dx=\mathfrak{c}(t),
	\end{aligned}
	\right.	
\end{equation}
with $(t,x)\in (0,\tau)\times \Omega$. Here, $\mathfrak{r},\mathfrak{p},\mathfrak{h},\mathfrak{q}, \mathfrak{c}(t)$ are known functions, $(\rho_{0},u_{0},\theta_{0})$ is the local smooth solution of incompressible Euler system \eqref{2.8-1} and $\tau$ is the lifespan given in Lemma \ref{lem1.1}. Throughout of this section, all functions are considered to be periodic in $x_{1}$ and $x_{2}$. Noting \eqref{4.33}, we should supplement \eqref{4.40} with following boundary conditions
\begin{equation}\label{4.41}
	\tilde{u}_{3}(t,x_{\sp},0)=-d_{0}(t,x_{\sp}),\quad \tilde{u}_{3}(t,x_{\sp},1)=d_{1}(t,x_{\sp}), \quad \forall (t,x_{\sp})\in (0,\tau)\times \mathbb{T}^2,
\end{equation}
and the initial data
\begin{equation}\label{4.42}
	(\tilde{\rho},\tilde{u},\tilde{\theta})(0,x)=(\tilde{\rho}_{0},\tilde{u}_{0},\tilde{\theta}_{0})(x).
\end{equation}

\begin{lemma}\label{lem4.4}
	Let $m\geq 3$ and $\tilde{\rho}_{0},\tilde{\theta}_{0}\in H^{m}(\Omega), \tilde{u}_{0}\in (H^{m}(\Omega))^3$ satisfy the compatibility conditions. (Here the compatibility conditions means the $(\tilde{\rho}_{0}, \tilde{\theta}_{0}, \tilde{u}_{0})$ satisfies $\eqref{4.40}_{1}$ and \eqref{4.41}, and the time-derivatives of initial data are defined through \eqref{4.40} inductively). Assume the following holds:
	$$
	\int_{\mathbb{T}^2}(d_{1}(t,x_{\sp})+d_{0}(t,x_{\sp})dx_{\sp}=\int_{\Omega}\operatorname{div}_{x}\tilde{u}dx=\int_{\Omega}\mathfrak{r}dx.
	$$
	Let
	$$
	\mathfrak{h}\in L^{\infty}(0,\tau; (H^{m}(\Omega))^3), \quad \mathfrak{p},\,\mathfrak{q},\,\pa_{t}\mathfrak{r}\in L^{\infty}(0,\tau; H^{m}(\Omega)),\quad \mathfrak{r}\in L^{\infty}(0,\tau; H^{m+1}(\Omega)),\,\,\mathfrak{c}(t)\in L^{\infty}(0,\tau),
	$$
	and
	$$
	d_{0},\,\,d_1\in L^{\infty}(0,\tau; H^{m+1}(\mathbb{T}^2)), \quad \pa_{t}d_{0},\,\,\pa_{t}d_1\in L^{\infty}(0,\tau; H^{m}(\mathbb{T}^2)).
	$$
	Then there is a unique solution $(\tilde{\rho},\tilde{u},\tilde{\theta},\tilde{p})$ of systems \eqref{4.40}-\eqref{4.42} with $\int_{\Omega}\tilde{p}(t,x)dx=0$ such that
	$$
	\tilde{\rho},\,\,\tilde{\theta}\in L^{\infty}(0,\tau; H^{m}(\Omega)),\quad \tilde{u}\in L^{\infty}(0,\tau;(H^{m}(\Omega))^3),\quad \tilde{p}\in L^{\infty}(0,\tau;H^{m+1}(\Omega)).
	$$
	Moreover, it holds
	\begin{equation}\label{4.64}
		\begin{aligned}
			&\sup_{t\in [0,\tau]}\|(\tilde{\rho},\tilde{u},\tilde{\theta})(t)\|_{H^{m}(\Omega)}+\sup_{t\in [0,\tau]}\| \nabla \tilde{p}(t)\|_{H^{m}(\Omega)}\\
			&\leq C(\tau, E_{m+2})\sup_{t\in [0,\tau]}\Big\{\|(\tilde{u}_{0},\tilde{\theta}_{0})\|_{H^{m}(\Omega)}+\|(\mathfrak{q},\mathfrak{h}
			,\mathfrak{r},\mathfrak{p})(t)\|_{H^{m}(\Omega)}+\|\partial_{t}\mathfrak{r}(t)\|_{H^{m-1}(\Omega)}\\
			&\quad +\|(d_{0},d_{1})(t)\|_{H^{m+1}(\mathbb{T}^2)}+\|(\partial_{t}d_{0},\partial_{t}d_{1})(t)\|_{H^{m}(\mathbb{T}^2)}\Big\},
		\end{aligned}
	\end{equation}
    and
    \begin{equation}\label{4.64-1}
    \begin{aligned}
    &\sup_{t\in [0,\tau]}\|\tilde{p}\|_{L^2(\Omega)}\\
    &\leq C(\tau, E_{m+2})\sup_{t\in [0,\tau]}\Big\{\|(\tilde{u}_{0},\tilde{\theta}_{0})\|_{H^{m}(\Omega)}+\|(\mathfrak{q},\mathfrak{h}
    ,\mathfrak{r},\mathfrak{p})(t)\|_{H^{m}(\Omega)}+\|\partial_{t}\mathfrak{r}(t)\|_{H^{m-1}(\Omega)}\\
    &\quad +\|(d_{0},d_{1})(t)\|_{H^{m+1}(\mathbb{T}^2)}+\|(\partial_{t}d_{0},\partial_{t}d_{1})(t)\|_{H^{m}(\mathbb{T}^2)}\Big\}+\|\mathfrak{c}(t)\|_{L^{\infty}},
    \end{aligned}
    \end{equation}
	where $E_{m}:=\sup_{t\in [0,\tau]}\|(u_{0},\theta_{0})(t)\|_{H^{m}(\Omega)}$. We emphasize that the estimate in \eqref{4.64} is independent of $\mathfrak{c}(t)$, which is important for us.
\end{lemma}

\noindent\textbf{Proof.} We divide the proof into four steps.

{\it Step 1.} We consider
\begin{equation}\label{4.44}
	\left\{\begin{aligned}
		&\operatorname{div}_{x}\hat{u}=\mathfrak{r}, \quad x\in \mathbb{T}^{2}\times [0,1],\\
		&\hat{u}_{3}(t,x_{\sp},0)=-d_{0}(t,x_{\sp}),\quad \hat{u}_{3}(t,x_{\sp},1)=d_{1}(t,x_{\sp}).
	\end{aligned}
	\right.
\end{equation}
We search the solution of the form $\hat{u}=\nabla_{x} q$ with $\int_{\Omega}q(t,x)dx=0$, then $q$ satisfies
\begin{equation}\label{4.45}
	\left\{\begin{aligned}
		&\Delta_{x} q=\mathfrak{r},\\
		&\partial_{x_3}q(t,x_{\sp},0)=-d_{0}(t,x_{\sp}),\quad \partial_{x_3}q(t,x_{\sp},1)=d_{1}(t,x_{\sp}).
	\end{aligned}
	\right.
\end{equation}
Using the classical theory of elliptic equation, we can solve \eqref{4.45}. In fact, there is a constant $C$ depending only on $m$ and $\Omega$ such that for $t\in [0,\tau]$,
\begin{equation}\label{4.46}
	\begin{aligned}
	\|\hat{u}(t)\|_{H^m(\Omega)}&=\|\nabla_{x} q(t)\|_{H^m(\Omega)}\leq C\left\{\|\mathfrak{r}(t)\|_{H^{m-1}(\Omega)}+\|(d_{0},d_{1})(t)\|_{H^{m}(\mathbb{T}^2)}\right\}.
	\end{aligned}
\end{equation}

{\it Step 2.} Let $v=\tilde{u}-\hat{u}$ and $p=\tilde{p}-\mathfrak{c}(t)$, it follows from \eqref{4.40} and \eqref{4.44} that
\begin{equation}\label{4.47}
	\left\{
	\begin{aligned}
		&\operatorname{div}_{x} v=0,\quad\int_{\Omega}p(t,x)dx=0,\\
		&\partial_{t}v+(u_{0}\cdot \nabla_{x})v+(v\cdot \nabla_{x})u_{0}+\nabla_{x} p=\mathfrak{h}-\partial_{t}\hat{u}-(u_{0}\cdot \nabla_{x})\hat{u}-(\hat{u}\cdot \nabla_{x})u_{0}:=\hat{\mathfrak{h}},\\
		&v_{3}(t,x_{\sp},0)=v_{3}(t,x_{\sp},1)=0.
	\end{aligned}
	\right.
\end{equation}
Applying $\operatorname{div}_{x}$ on both sides of $\eqref{4.47}_2$, one has
\begin{equation}\label{4.48}
	\left\{\begin{aligned}
		&\Delta_{x} p=\operatorname{div}_{x}\mathfrak{h}-2\sum\limits_{i,j=1}^3\partial_{i}u_{0,j}\partial_{j}v_{i}-\partial_{t}\mathfrak{r}-2\sum\limits_{i,j=1}^3\partial_{i}u_{0,j}\partial_{j}\hat{u}_{i}-\sum\limits_{i=1}^3u_{0,i}\partial_{i}\mathfrak{r},\\
		&\partial_{x_3}p(t,x_{\sp},0)=\mathfrak{h}_{3}(t,x_{\sp},0)+(u_{0,\sp}^{-}\cdot \nabla_{\sp})d_{0}(t,x_{\sp})+d_{0}(t,x_{\sp})\partial_{x_3}u_{0,3}^{-},\\
		&\partial_{x_3}p(t,x_{\sp},1)=\mathfrak{h}_{3}(t,x_{\sp},1)-(u_{0,\sp}^{+}\cdot \nabla_{\sp})d_{1}(t,x_{\sp})-d_{1}(t,x_{\sp})\partial_{x_3}u_{0,3}^{+}
		\\&\int_{\Omega}p(t,x)dx=0,
	\end{aligned}
	\right.
\end{equation}
Using the classical theory of elliptic equation and \eqref{4.46}, we have
\begin{equation}\label{4.49}
	\begin{aligned}
		\|\nabla_{x} p(t)\|_{H^{m}(\Omega)}
		&\leq C\Big\{\|\mathfrak{h}(t)\|_{H^{m}(\Omega)}+\|\partial_{t}\mathfrak{r}(t)\|_{H^{m-1}(\Omega)}+\|u_{0}(t)\|_{H^{m}(\Omega)}\|(v,\mathfrak{r})(t)\|_{H^{m}(\Omega)}\\
		&\qquad\quad  +\|u_0(t)\|_{H^{m+2}(\Omega)}\|(d_{0},d_{1})(t)\|_{H^{m+1}(\mathbb{T}^2)}
		\Big\},
	\end{aligned}
\end{equation}

Using the classical Poincar\'{e} inequality, one has
$$
\|\tilde{p}-\mathfrak{c}(t)\|_{L^2(\Omega)}\leq C\|\nabla_{x} \tilde{p}\|_{L^2(\Omega)}=C\|\nabla_{x} p\|_{L^2(\Omega)},
$$
which, together with \eqref{4.49}, yields that
\begin{equation}\label{4.49-1}
\begin{aligned}
\|\tilde{p}(t)\|_{L^2(\Omega)}&\leq C\|\nabla_{x} p(t)\|_{H^{m}(\Omega)}+\|\mathfrak{c}(t)\|_{L^{\infty}}\\
&\leq C\Big\{\|\mathfrak{h}(t)\|_{H^{m}(\Omega)}+\|\partial_{t}\mathfrak{r}(t)\|_{H^{m-1}(\Omega)}+\|u_{0}(t)\|_{H^{m}(\Omega)}\|(v,\mathfrak{r})(t)\|_{H^{m}(\Omega)}\\
&\qquad\quad  +\|u_0(t)\|_{H^{m+2}(\Omega)}\|(d_{0},d_{1})(t)\|_{H^{m+1}(\mathbb{T}^2)}
\Big\}+\|\mathfrak{c}(t)\|_{L^{\infty}}.
\end{aligned}
\end{equation}

{\it Step 3.} Let $\alpha=(\alpha_{1},\alpha_2,\alpha_{3})$ be a multi-index with $|a|=\alpha_1+\alpha_2+\alpha_3\leq m$. Applying  $D^{\alpha}:=\partial_{x_1}^{\alpha_{1}}\partial_{x_{2}}^{\alpha_{2}}\partial_{x_3}^{\alpha_3}$ to $\eqref{4.47}_2$ and multiplying the resultant equation by $D^{\alpha} v$, then we can obtain
\begin{equation}\label{4.50}
	\frac{1}{2}\frac{d}{dt}\|v\|_{H^{m}(\Omega)}^2=-\la (u_{0}\cdot \nabla_{x})v,v\ra_{H^{m}}-\la (v\cdot \nabla_{x})u_{0},v\ra_{H^{m}}-\la \nabla_{x} p,v\ra_{H^{m}}-\la \hat{\mathfrak{h}},v\ra_{H^{m}}.
\end{equation}
For the last term on the right hand side (RHS) of \eqref{4.50}, it follows from \eqref{4.46} that
\begin{equation}\label{4.51}
	\begin{aligned}
		\vert\la \hat{\mathfrak{h}},v\ra_{H^{m}}\vert&\leq \|\hat{\mathfrak{h}}\|_{H^{m}(\Omega)}\|v\|_{H^{m}(\Omega)}\leq (\|(\mathfrak{h},\partial_{t}\hat{u})\|_{H^{m}(\Omega)}+\|u_{0}\|_{H^{m+1}(\Omega)}\|\hat{u}\|_{H^{m+1}(\Omega)})\|v\|_{H^{m}(\Omega)}\\
		&\leq C(\|\mathfrak{h}\|_{H^{m}(\Omega)}+\|\partial_{t}\mathfrak{r}\|_{H^{m-1}(\Omega)}+\|(\partial_{t}d_{0},\partial_{t}d_{1})\|_{H^{m}(\mathbb{T}^2)})\|v\|_{H^{m}(\Omega)}\\
		&\quad +C\|u_{0}\|_{H^{m+1}(\Omega)}\left(\|\mathfrak{r}\|_{H^{m}(\Omega)}+\|(d_{0},d_{1})\|_{H^{m}(\mathbb{T}^2)}\right)\|v\|_{H^{m}(\Omega)}.
	\end{aligned}
\end{equation}
For the third term on the RHS of \eqref{4.50}, we have from \eqref{4.49} that
\begin{equation}\label{4.52}
	\begin{aligned}
	\vert\la \nabla_{x} p,v\ra_{H^{m}}\vert&\leq \|\nabla_{x} p\|_{H^{m}(\Omega)}\|v\|_{H^{m}(\Omega)}\\
		&\leq C\|u_{0}\|_{H^{m}(\Omega)}\|v\|_{H^{m}(\Omega)}^2+C\big[\|\mathfrak{h}\|_{H^{m}(\Omega)}+\|\partial_{t}\mathfrak{r}\|_{H^{m-1}(\Omega)}\\
		&\quad +\|u_{0}\|_{H^{m+2}(\Omega)}\left(\|\mathfrak{r}\|_{H^{m}(\Omega)}+\|(d_{0},d_{1})\|_{H^{m+1}(\mathbb{T}^2)}\right)\big]\|v\|_{H^{m}(\Omega)}.
	\end{aligned}
\end{equation}
For the second term on the RHS of \eqref{4.50}, it is clear that
\begin{equation}\label{4.53}
\vert\la (v\cdot \nabla_{x})u_{0},v\ra_{H^{m}}\vert\leq \|u_{0}\|_{H^{m+1}(\Omega)}\|v\|_{H^{m}(\Omega)}^2.
\end{equation}
It remains to estimate the first term on the RHS of \eqref{4.50}. We notice that
$$
D^{\alpha}((u_{0}\cdot \nabla_{x})v)=(u_{0}\cdot \nabla_{x})D^{\alpha}v+\sum\limits_{|\beta|>0\atop |\beta|+|\gamma|=|\alpha|}c_{\beta,\gamma}(D^{\beta}u_{0}\cdot \nabla_{x})D^{\gamma}v,
$$
which, together with $\operatorname{div}_{x}u_{0}=0$ and integrating by parts, yields that
\begin{align}\label{4.54}
 \vert\la (u_{0}\cdot \nabla_{x})v,v\ra_{H^{m}}\vert&\leq \sum\limits_{|\alpha|\leq m}\vert\la D^{\alpha}((u_{0}\cdot \nabla_{x})v),D^{\alpha}v\ra_{L^2}\vert\nonumber\\
 &=\sum\limits_{|\alpha|\leq m}\vert\la(u_{0}\cdot \nabla_{x})D^{\alpha}v,D^{\alpha}v\ra_{L^2}\vert+\sum\limits_{|\beta|>0,|\alpha|\leq m\atop |\beta|+|\gamma|=|\alpha|}\vert c_{\beta,\gamma}\la (D^{\beta}u_{0}\cdot \nabla_{x})D^{\gamma}v,D^{\alpha}v\ra_{L^2}\vert\nonumber\\
 &\leq C\|u_{0}\|_{H^{m}(\Omega)}\|v\|_{H^{m}(\Omega)}^2.
\end{align}
Substituting \eqref{4.51}-\eqref{4.54} into \eqref{4.50}, one has
\begin{equation}\label{4.55}
	\begin{aligned}
		\frac{d}{dt}\|v\|_{H^{m}(\Omega)}&\leq C\|u_{0}\|_{H^{m+1}(\Omega)}\|v\|_{H^{m}(\Omega)} +C\|\mathfrak{h}\|_{H^{m}(\Omega)}+\|\partial_{t}\mathfrak{r}\|_{H^{m-1}(\Omega)}\\
		&+\|(\partial_{t}d_{0},\partial_{t}d_{1})\|_{H^{m}(\mathbb{T}^2)} +C\|u_{0}\|_{H^{m+2}(\Omega)}(\|\mathfrak{r}\|_{H^{m}(\Omega)}+\|(d_{0},d_{1})\|_{H^{m+1}(\mathbb{T}^2)}).
	\end{aligned}
\end{equation}
Hence, using the Gr\"{o}nwall inequality, we deduce from \eqref{4.55} that
\begin{equation}\label{4.56}
	\begin{aligned}
		\sup_{t\in [0,\tau]}\|v(t)\|_{H^{m}(\Omega)}&\leq C(\tau, E_{m+2})\sup_{t\in [0,\tau]}\Big\{\|v_{0}\|_{H^{m}(\Omega)}+\|\hat{\mathfrak{h}}(t)\|_{H^{m}(\Omega)}+\|\partial_{t}\mathfrak{r}(t)\|_{H^{m-1}(\Omega)}\\
		&+\|\mathfrak{r}(t)\|_{H^{m}(\Omega)}+\|(\partial_{t}d_{0},\partial_{t}d_{1})(t)\|_{H^{m}(\mathbb{T}^2)}+\|(d_{0},d_{1})(t)\|_{H^{m+1}(\mathbb{T}^2)}\Big\}\\
		&\leq C(\tau, E_{m+2})\sup_{t\in [0,\tau]}\Big\{\|\tilde{u}_{0}\|_{H^{m}(\Omega)}+\|\hat{\mathfrak{h}}(t)\|_{H^{m}(\Omega)}+\|\partial_{t}\mathfrak{r}(t)\|_{H^{m-1}(\Omega)}\\
		&+\|\mathfrak{r}(t)\|_{H^{m}(\Omega)}+\|(\partial_{t}d_{0},\partial_{t}d_{1})(t)\|_{H^{m}(\mathbb{T}^2)}+\|(d_{0},d_{1})(t)\|_{H^{m+1}(\mathbb{T}^2)}\Big\}
	\end{aligned}
\end{equation}
where $v_{0}:=v(0,x)=\tilde{u}_{0}-\hat{u}(0,x)$ is the initial data of $v$, and we have used \eqref{4.46} and
$$
\begin{aligned}
\|v_{0}\|_{H^{m}(\Omega)}
&\leq \|\tilde{u}_{0}\|_{H^{m}(\Omega)}+C\sup_{t\in [0,\tau]}\{\|\mathfrak{r}(t)\|_{H^{m-1}(\Omega)}+\|(d_{0},d_{1})(t)\|_{H^{m}(\mathbb{T}^2)}\}
\end{aligned}
$$
in the last inequality.

Due to the priori estimates on $p$ and $v$ obtained in \eqref{4.49} and \eqref{4.56} respectively, the existence and uniqueness of smooth solution $(v,p)$ to linear systems \eqref{4.47} over $[0,\tau]\times \Omega$ can be obtained by standard arguments, which are omitted for simplicity of presentation.

Once we have obtained the existence and uniqueness of $v$, the existence and uniqueness of $\tilde{u}$ and $\tilde{p}$ can be directly obtained by the constructions $\tilde{u}=v+\hat{u}$ and $\tilde{p}=p+\mathfrak{c}(t)$, respectively.  Moreover, combining \eqref{4.46} with \eqref{4.56}, one deduces that
\begin{equation}\label{4.57}
	\begin{aligned}
		\sup_{t\in [0,\tau]}\|\tilde{u}(t)\|_{H^{m}(\Omega)}&\leq C(\tau, E_{m+2})\sup_{t\in [0,\tau]}\Big\{\|\tilde{u}_{0}\|_{H^{m}(\Omega)}+\|(\mathfrak{h},\mathfrak{r})(t)\|_{H^{m}(\Omega)}+\|\partial_{t}\mathfrak{r}(t)\|_{H^{m-1}(\Omega)}\\
		&+\|(\partial_{t}d_{0},\pa_{t}d_{1})(t)\|_{H^{m}(\mathbb{T}^2)} +\|(d_{0},d_{1})(t)\|_{H^{m+1}(\mathbb{T}^2)}\Big\}.
	\end{aligned}
\end{equation}

{\it Step 4.} Since $\tilde{u}$ is already a known function, the existence of $\tilde{\theta}$ is trivial by classical transport theory. For the energy estimate of $\tilde{\theta}$, it follows from \eqref{4.40} that
\begin{equation}\label{4.58}
	\frac{1}{2}\frac{d}{dt}\|\tilde{\theta}\|_{H^{m}(\Omega)}^2=-\la (u_{0}\cdot \nabla_{x})\tilde{\theta},\tilde{\theta}\ra_{H^{m}}-\la (\tilde{u}\cdot \nabla_{x})\theta_{0},\tilde{\theta}\ra_{H^{m}}+\la \mathfrak{q},\theta\ra_{H^{m}}.
\end{equation}
It is clear that
\begin{equation}\label{4.59}
\vert\la \mathfrak{q},\tilde{\theta}\ra_{H^{m}}\vert+ \vert\la (\tilde{u}\cdot \nabla_{x})\theta_{0},\tilde{\theta}\ra_{H^{m}}\vert\leq (\|\mathfrak{q}\|_{H^{m}(\Omega)}+ \|\tilde{u}\|_{H^{m}(\Omega)}\|\theta_{0}\|_{H^{m+1}(\Omega)})\|\tilde{\theta}\|_{H^{m}(\Omega)}.
\end{equation}
Using a similar argument as in \eqref{4.54}, one gets
\begin{equation}\label{4.60}
	\vert\la (u_{0}\cdot \nabla_{x})\tilde{\theta},\tilde{\theta}\ra_{H^{m}}\vert\leq C\|u_{0}\|_{H^{m}(\Omega)}\|\tilde{\theta}\|_{H^{m}(\Omega)}^2.
\end{equation}
Substituting \eqref{4.59}-\eqref{4.60} into \eqref{4.58} and using Gr\"{o}nwall inequality, one deduces that
\begin{equation}\label{4.61}
	\begin{aligned}
		&\sup_{t\in [0,\tau]}\|\tilde{\theta}(t)\|_{H^{m}(\Omega)}\leq C(\tau,E_{m+1})\sup_{t\in [0,\tau]}(\|\tilde{\theta}_{0}\|_{H^{m}(\Omega)}+\|\tilde{u}(t)\|_{H^{m}(\Omega)}+\|\mathfrak{q}(t)\|_{H^{m}(\Omega)})\\
		&\leq C(\tau,E_{m+2})\sup_{t\in [0,\tau]}\Big\{\|\tilde{\theta}_{0}\|_{H^{m}(\Omega)}+\|(\mathfrak{q},\mathfrak{h},\mathfrak{r})(t)\|_{H^{m}(\Omega)}+\|\partial_{t}\mathfrak{r}(t)\|_{H^{m-1}(\Omega)}\\
		&\quad  +\|(\partial_{t}d_{0},\partial_{t}d_{1})(t)\|_{H^{m}(\mathbb{T}^2)}+\|(d_{0},d_{1})(t)\|_{H^{m+1}(\mathbb{T}^2)}\Big\},
	\end{aligned}
\end{equation}
where we have used \eqref{4.57} in the last inequality.

Finally, recalling $\tilde{\rho}+\tilde{\theta}=\mathfrak{p}$, we obtain from \eqref{4.61} that
\begin{align}\label{4.62}
  \sup_{t\in [0,\tau]}\|\tilde{\rho}(t)\|_{H^{m}(\Omega)}
  &\leq  C(\tau,E_{m+2})\sup_{t\in [0,\tau]}\Big\{\|\tilde{\theta}_{0}\|_{H^{m}(\Omega)}+\|(\mathfrak{q},\mathfrak{h},\mathfrak{r},\mathfrak{p})(t)\|_{H^{m}(\Omega)}+\|\partial_{t}\mathfrak{r}(t)\|_{H^{m-1}(\Omega)}\nonumber\\
  &\quad +\|(\partial_{t}d_{0},\partial_{t}d_{1})(t)\|_{H^{m}(\mathbb{T}^2)}+\|(d_{0},d_{1})(t)\|_{H^{m+1}(\mathbb{T}^2)}\Big\},
\end{align}
which, together with \eqref{4.49}-\eqref{4.49-1}, \eqref{4.57} and \eqref{4.61}, yields \eqref{4.64}--\eqref{4.64-1}. Therefore the proof of Lemma \ref{lem4.4} is complete. $\hfill\square$

\subsection{Existence of solutions for a kind of  linear parabolic system}
As mentioned in Remark \ref{rem4.3}, we have $u_{1,3}^{\pm}=0$. Then $\eqref{3.12}_{1}-\eqref{3.12}_{3}$ can be written as following the linear parabolic system of $(u,\theta):=(u_1,u_2,\t)$
\begin{equation}\label{5.1}
	\left\{\begin{aligned}
		&\partial_{t}u_{i}+(u_{0,\sp}^{\pm}\cdot \nabla_{\sp})u_{i}+(y\partial_{x_3}u_{0,3}^{\pm})\partial_{y}u_{i}+(u\cdot \nabla_{\sp})u_{0,i}^{\pm}=\lambda\partial_{yy}u_{i}+\mathfrak{f}_{i},\;i=1,2,\\
		&\partial_{t}\theta+(u_{0,\sp}^{\pm}\cdot \nabla_{\sp})\t+(y\partial_{x_3}u_{0,3}^{\pm})\partial_{y}\t+(u\cdot \nabla_{\sp})\t_{0}^{\pm}=\frac{2}{5}\k\partial_{yy}\t+\mathfrak{g},
	\end{aligned}\right.
\end{equation}
where $(t, x_{\sp}, y) \in[0, \tau] \times \mathbb{T}^{2} \times \R_{+}$, and $(u_{0}^{\pm}, \theta_{0}^{\pm})$ are the corresponding boundary values of incompressible Euler solutions, which are independent of $y \in \mathbb{R}_{+}$. Recalling \eqref{4.37}, we supplement system \eqref{5.1} with the non-homogeneous Neumann boundary conditions, i.e.,
\begin{equation}\label{5.2}
	\left\{\begin{aligned}
		&\partial_{y} u_{i}\left(t, x_{\sp}, y\right)\vert_{y=0}=b_{i}\left(t, x_{\sp}\right),\quad i=1,2,\\
		&\partial_{y} \theta\left(t, x_{\sp}, y\right)\vert_{y=0}=a\left(t, x_{\sp}\right), \quad\lim _{y \rightarrow+\infty}(u, \theta)\left(t, x_{\sp}, y\right)=0,
	\end{aligned}\right.
\end{equation}
and the initial data
\begin{equation}\label{5.3}
	u\left(t, x_{\sp}, y\right)\vert_{t=0}:=\mathfrak{u}_{0}\left(x_{\sp}, y\right),\quad \theta\left(t, x_{\sp}, y\right)\vert_{t=0}:=\mathfrak{t}_{0}\left(x_{\sp}, y\right).
\end{equation}
The initial data $\left(\mathfrak{u}_{0}, \mathfrak{t}_{0}\right)$ should satisfy the corresponding compatibility conditions. Let $l \geq 0$, we define the notations
\begin{equation}\label{5.3-1}
	\|f\|_{L_{l}^{2}}^{2}=\int_{\R_{+}}\int_{\mathbb{T}^2}(1+y)^{l}|f(x_{\sp}, y)|^{2} d x_{\sp} d y,
\end{equation}
and
\begin{equation}\label{5.3-2}
\bar{x}:=(x_{\sp}, y), \quad \nabla_{\bar{x}}:=(\nabla_{\sp}, \partial_{y}) \equiv(\partial_{x_{1}}, \partial_{x_{2}}, \partial_{y}).
\end{equation}

The existence theory and uniform estimates for system \eqref{5.1} had been given in \cite{Guo-Huang-Wang}.
\begin{lemma}[{\cite[Lemma 4.1]{Guo-Huang-Wang}}]\label{lem5.1}
	Let $l\geq 0, k\geq 3$. Let the compatibility condition for the initial data \eqref{5.3} be satisfied (Here the compatibility condition means that the initial data \eqref{5.3} satisfies the boundary condition \eqref{5.2}, and the time-derivatives of initial data $(\mathfrak{u}_0,\mathfrak{t}_0)$ are defined through system \eqref{5.1} inductively). Assume that
	$$
	\sup _{t \in[0, \tau]}\Big\{\sum_{|\beta|+2 \gamma \leq k+2}\|\partial_{t}^{\gamma}\nabla_{\sp}^{\beta} (a, b)(t)\|_{L^{2}\left(\mathbb{T}^{2}\right)}^{2}+\sum_{j=0}^{k} \sum_{|\beta|+2 \gamma=j}\|\partial_{t}^{\gamma}\nabla_{\bar{x}}^{\beta} (\mathfrak{f}, \mathfrak{g})(t)\|_{L_{l_{j}}^{2}}^{2}\Big\}<\infty,
	$$
	with $l_{j}:=l+2(k-j), 0 \leq j \leq k .$ Then there exists a unique smooth solution $(u, \theta)$ of $\eqref{5.1}-\eqref{5.3}$ over $t \in[0, \tau]$, which satisfies
	\begin{equation}\nonumber
		\begin{aligned}
			&\sum_{j=0}^{k} \sum_{|\beta|+2 \gamma=j} \sup _{t \in[0, \tau]}\Big\{\|\partial_{t}^{\gamma} \nabla_{\bar{x}}^{\beta}(u, \theta)(t)\|_{L_{l_{j}}^{2}}^{2}+\int_{0}^{t}\|\partial_{t}^{\gamma} \nabla_{\bar{x}}^{\beta} \partial_{y}(u, \theta)(s)\|_{L_{l_{j}}^{2}}^{2} d s\Big\} \\
			&\leq C(\tau, E_{k+3})\sup _{t \in[0, \tau]}\Big\{\sum_{j=0}^{k} \sum_{|\beta|+2 \gamma=j}\|\partial_{t}^{\gamma} \nabla_{\bar{x}}^{\beta}(u, \theta)(0)\|_{L_{l_{j}}^{2}}^{2}\\
			&\quad +\sum_{|\beta|+2 \gamma \leq k+2}\|\nabla_{\sp}^{\beta} \partial_{t}^{\gamma}(a, b)(t)\|_{L^{2}(\mathbb{T}^{2})}^{2} +\sum_{j=0}^{k} \sum_{|\beta|+2 \gamma=j}\|\partial_{t}^{\gamma}\nabla_{\bar{x}}^{\beta} (\mathfrak{f}, \mathfrak{g})(t)\|_{L_{l_{j}}^{2}}^{2}\Big\},
		\end{aligned}
	\end{equation}
	where $E_{k}$ is defined in Lemma \ref{lem4.4}.
\end{lemma}

\section{The solutions of expansions}
For later use, we introduce following  weight functions
\begin{align}\label{6.1}
	\tw_{\k_i}(v)=w_{\k_i}(v) \mu(v)^{-\fa},
	\quad \fw_{\bar{\k}_i}(v)=w_{\bar{\k}_i}(v) \mu(v)^{-\fa} \,\quad
	\mbox{and}\,\quad
	\fw_{\hat{\k}_i}(v)=w_{\hat{\k}_i}(v) \mu(v)^{-\fa},
\end{align}
for constants $\k_i, \bar{\k}_i, \hat{\k}_i\geq 0,\, 1\leq i\leq N$ and $0\leq \fa <\frac{1}{2}$. We denote
\begin{align}\nonumber
 \hat{x}=(x_{\sp} ,\eta)\in \mathbb{T}^2\times \R_{+}, \quad \nabla_{\hat{x}}:=(\nabla_{\sp} ,\pa_{\eta}).
\end{align}
Recall \eqref{5.3-2} and the weighted $L_{l}^2-$norm in \eqref{5.3-1}.

\begin{proposition} \label{prop5.1}
	Let $0\leq \fa<\frac12$ in \eqref{6.1}. Let  $s_0, s_i, \bar{s}_i, \hat{s}_i\in \mathbb{N}_+$, $\k_i, \bar{\k}_{i}, \hat{\k}_i \in \R_+$ for $1\leq i\leq N$; and define  $l_j^i:=\bar{l}_i+2(\bar{s}_i-j)$ for $1\leq i \leq N, \  0\leq j\leq \bar{s}_i$. For these parameters, we assume the restrictions  \eqref{5.84}-\eqref{5.86} hold. Let the initial data $(\rho_{i}, u_i, \theta_i)(0)$ of IBVP \eqref{2.17}, \eqref{4.33}-\eqref{4.34}, and  initial data $(\bar{u}_{i,\sp}^{\pm},\bar{\theta}_{i}^{\pm})(0)$ of IBVP \eqref{3.12}, \eqref{4.37} satisfy
	\begin{align}\label{5.0}
		\sum_{i=0}^{N}\Big\{\sum_{\gamma+|\beta|\leq s_{i}}\|\pa_t^\gamma\nabla_x^{\beta}(\rho_{i}, u_{i},  \theta_{i})(0)\|_{L^2_{x}}+\sum_{j=0}^{\bar{s}_{i}} \sum_{j=2\gamma+|\beta|}   \|\pa_t^{\gamma} \nabla_{\bar{x}}^{\beta} (\bar{u}_{i,\sp}^{\pm},\bar{\theta}_{i}^{\pm})(0)\|^2_{L^2_{l_j^{i}}}
		\Big\}<\infty.
	\end{align}
	And we also assume that the compatibility conditions for initial data $(\rho_{i}, u_i, \theta_i)(0)$ and $(\bar{u}_{i,\sp}^{\pm},\bar{\theta}_{i}^{\pm})(0)$ are satisfied.
	Then there exist solutions $F_i=\sqrt{\mu}f_i, \, \bar{F}_{i}^{\pm}=\sqrt{\mu}\bar{f}_{i}^{\pm}, \, \hat{F}_{i}^{\pm}=\sqrt{\mu}\hat{f}_{i}^{\pm}$ to \eqref{1.7-1}, \eqref{1.14} and \eqref{1.19-1} over the time interval  $t\in[0,\tau]$ respectively. Moreover, we have
	\begin{align}\label{5.0-1}
		\begin{split}
			&\sup_{t\in[0,\tau]} \sum_{i=1}^N\Big\{\sum_{\gamma+|\beta|\leq s_i} \|\tw_{\k_i}\pa_t^\gamma \nabla_x^{\beta} f_{i}(t)\|_{L^2_xL^\infty_v}
			+\sum_{j=0}^{\bar{s}_i}\sum_{j=2\gamma+|\beta|} \|\fw_{\bar{\k}_i}\pa_t^\gamma \nabla_{\bar{x}}^{\beta} \bar{f}_{i}^{\pm}(t)\|_{L^2_{l^i_j}L^\infty_v}\\
			&\qquad\qquad\qquad+\sum_{\gamma+|\beta|\leq \hat{s}_i} \| e^{\zeta_i \eta}\fw_{\hat{\k}_i}\pa_t^\gamma \nabla_{\sp} ^{\beta} \hat{f}_{i}^{\pm}(t)\|_{L^\infty_{\hat{x},v}}\Big\}\\
			&\leq C\Big( \tau,\sum_{i=0}^{N}\Big[\sum_{\gamma+|\beta|\leq s_{i}}\|\pa_t^\gamma\nabla_x^{\beta}(\rho_{i}, u_{i},  \theta_{i})(0)\|_{L^2_{x}}+\sum_{j=0}^{\bar{s}_{i}} \sum_{j=2\gamma+|\beta|}   \|\pa_t^{\gamma} \nabla_{\bar{x}}^{\beta} (\bar{u}_{i,\sp}^{\pm},\bar{\theta}_{i}^{\pm})(0)\|^2_{L^2_{l_j^{i}}}\Big]
			\Big),
		\end{split}
	\end{align}
	where the positive constants $\zeta_i>0 \, (i=1,\cdots,N)$ satisfying $\zeta_{i+1}\leq \frac12\zeta_{i}$ and $\zeta_1=1$.
\end{proposition}

\begin{remark}
	Since the Knudsen boundary layer $\hat{f}_{k}^{\pm}$ is indeed a stationary problem with $(t,x_{\sp} )$ as parameters, hence there is no necessary to give initial data for $\hat{f}_{k}^{\pm}$.
\end{remark}

\noindent\textbf{Proof.} For simplicity of presentations, we consider the case $n=3$ throughout the proof. For the case $n\geq 4$, the proof is similar and slightly easier. We point out that the framework of the proof is similar to \cite[Proposition 5.1]{Guo-Huang-Wang}. We divide the proof into three steps since the proof is long.

\noindent{\it{Step 1. Construction of solutions $f_1$, $\bar{f}_{1}^{\pm}$ and $\hat{f}_{1}^{\pm}$.}}

{\it{Step 1.1 Construction of $f_1$.}} It follows from \eqref{2.6-1} that the microscopic part of $f_1$ is given by
\begin{equation}\label{5.40-1}
\begin{aligned}
(\mathbf{I-P})f_1
&=\frac{1}{2}\sum\limits_{l,m=1}^{3}u_{0,l}u_{0,m}\mathcal{A}_{lm}+\sum\limits_{m=1}^3u_{0,m}\t_{0}\mathcal{B}_{m}+\frac{\t_{0}^2}{8}(\mathbf{I-P})\{(|v|^2-5)^2\sqrt{\mu}\},
\end{aligned}
\end{equation}
which, together with Lemma \ref{lem1.1}, yields that
\begin{equation}\label{5.41-2}
\sup_{t\in [0,\tau]}\sum\limits_{\g+|\beta|\leq s_1}\|\tw_{\k_1}\pa_{t}^{\g}\nabla_{x}(\mathbf{I-P})f_{1}(t)\|_{L_{x}^2L_{v}^2}\leq C\big(\tau, \sum\limits_{\g+|\beta|\leq s_{1}+2}\|\pa_{t}^{\g}\nabla_{x}^{\beta}(\rho_{0},u_{0},\t_{0})(0)\|_{L_{x}^2}\big).
\end{equation}

For the macroscopic part $(\rho_{1},u_{1},\theta_{1})$, it follows from Lemmas \ref{lem2.1}, \ref{lem2.2} and Remark \ref{rem4.3} that
\begin{equation}\label{5.40}
	\left\{\begin{aligned}
		&\operatorname{div}_{x}u_1=-\pa_{t}\rho_{0},\quad \rho_{1}+\theta_{1}-\frac{1}{3}|u_{0}|^2=p_{1},\\
		&\partial_{t}u_{1}+(u_{0}\cdot \nabla_{x})u_{1}+(u_{1}\cdot \nabla_{x})u_{0}+\nabla_{x} p_{2}=\mathfrak{h}_{0}=\pa_{t}\rho_{0}u_{0},\\
		&\partial_{t}\theta_{1}+(u_{0}\cdot \nabla_{x})\theta_{1}+(u_{1}\cdot \nabla_{x})\theta_{0}=\mathfrak{q}_{0}=\frac{2}{5}\pa_{t}p_{1}+\frac{2}{15}\pa_{t}|u_{0}|^2+\pa_{t}\rho_{0}\theta_{0},\\
		&u_{1,3}(t,x_{\sp},0)=u_{1,3}(t,x_{\sp},1)=0,\,\,\,\int_{\Omega}p_{2}dx=\mathfrak{c}_{2}(t),
	\end{aligned}
	\right.
\end{equation}
where $p_{1}$ is the pressure determined by the incompressible Euler system \eqref{2.8-1}-\eqref{2.8-4}.

It follows from Lemma \ref{lem1.1}, Remark \ref{rem2.4} and the Poincar\'{e} inequality that
\begin{equation}\label{5.40-2}
\begin{aligned}
\sup_{t\in [0,\tau]}\sum\limits_{\gamma+|\beta|\leq s_{1}}\|\pa_{t}\nabla_{x}^{\beta}p_{1}\|_{L_{x}^2}&=\sup_{t\in [0,\tau]}\big(\sum\limits_{\gamma\leq s_{1}}\|\pa_{t}^{\gamma}p_{1}(t)\|_{L_{x}^2}+\sum\limits_{\gamma+|\beta|\leq s_{1}\atop |\beta|\geq 1}\|\pa_{t}^{\gamma}\nabla_{x}^{\beta-1}(\nabla_{x} p_{1})(t)\|_{L_{x}^2}\big)\\
&\leq\sum\limits_{\gamma\leq s_{1}}\|\pa_{t}^{\gamma}\mathfrak{c}_{1}\|_{L_{t}^{\infty}}+\sup_{t\in [0,\tau]}\sum\limits_{\gamma+|\beta|\leq s_{1}\atop |\beta|\geq 1}\|\pa_{t}^{\gamma}\nabla_{x}^{\beta-1}(\nabla_{x} p_{1})(t)\|_{L_{x}^2}\\
&\leq C\Big(\tau, \sum\limits_{i=0}^{1}\sum\limits_{\gamma+|\beta|\leq s_{i}}\|\partial_{t}^{\gamma}\nabla_{x}^{\beta}(\rho_{i},u_i,\theta_{i})(0)\|_{L_{x}^2}\Big),\quad\text{for } s_{0}\geq s_{1}.
\end{aligned}
\end{equation}

Using \eqref{5.40-2} and Lemma \ref{lem4.4}, we obtain the existence of solution $(\r_1,u_1,\t_1)$ and
\begin{equation}\label{5.41}
	\begin{aligned}
		&\sup_{t\in [0,\tau]}\sum\limits_{\gamma+|\beta|\leq s_{1}}\|\pa_{t}^{\gamma}\nabla_{x}^{\beta}(\rho_{1},u_{1},\theta_{1})(t))\|_{L_{x}^2}+\sup_{t\in [0,\tau]}\sum\limits_{\gamma+|\beta|\leq s_{1}}\|\partial_{t}^{\gamma}\nabla_{x}^{\beta}(\nabla_{x}p_{2})\|_{L_{x}^2}\\
		&
		\leq C\Big(\tau, \sum\limits_{i=0}^{1}\sum\limits_{\gamma+|\beta|\leq s_{i}}\|\partial_{t}^{\gamma}\nabla_{x}^{\beta}(\rho_{i},u_{i},\theta_{i})(0)\|_{L_{x}^2}\Big).
	\end{aligned}
\end{equation}
We point out that we cannot deduce the estimate on $\|\partial_{t}^{\gamma}p_{2}\|_{L_{x}^2}$ directly from \eqref{5.41} and the classical Poincar\'{e} inequality, since $\int_{\Omega}p_{2}(t)dx=\mathfrak{c}_{2}(t)$ may depend on $\bar{\rho}_{1}^{\pm}$ from \eqref{4.82}. Combining \eqref{5.41} with \eqref{5.41-2}, one deduces that for $s_0\geq s_1+2$, it holds that
\begin{equation}\label{5.42}
	\begin{aligned}
		\sup_{t\in [0,\tau]}\sum\limits_{\gamma+|\beta|\leq s_1}\|\tw_{\k_{1}}\partial_{t}^{\gamma}\nabla_{x}f_{1}(t)\|_{L_{x}^2L_{v}^{\infty}}\leq C\Big(\tau, \sum\limits_{i=0}^{1}\sum\limits_{\gamma+|\beta|\leq s_{i}}\|\partial_{t}^{\gamma}\nabla_{x}^{\beta}(\rho_{i},u_{i},\theta_{i})(0)\|_{L_{x}^2}\Big).
	\end{aligned}
\end{equation}

{\it{Step 1.2. Construction of $\bar{f}_{1}^{\pm}$.}} Noting from \eqref{3.1} that $\bar{f}_{1}^{\pm}\in \mathcal{N}$, we only need to construct the macroscopic part $(\bar{\r}_{1}^{\pm},\bar{u}_{1}^{\pm},\bar{\t}_{1}^{\pm})$. Recalling \eqref{4.35-3}, the boundary conditions for $(\bar{u}_{1,\sp}^{\pm},\bar{\t}_{1}^{\pm})$ are
\begin{equation}\label{5.43}
	\left\{
	\begin{aligned}
	&\partial_{y}\bar{u}_{1,i}^{\pm}(t,x_{\sp},0)=\pm \pa_{x_{3}}u_{0,i}^{\pm},\quad i=1,2,\\ &\partial_{y}\bar{\theta}_{1}^{\pm}(t,x_{\sp},0)=\pm \pa_{x_{3}}\theta_{0}^{\pm}.
	\end{aligned}
\right.
\end{equation}
It follows from Proposition \ref{prop3.1} and \eqref{3.20-1} that $(\bar{\rho}_{1}^{\pm},  \bar{u}_{1}^{\pm},\bar{\t}_{1}^{\pm})$ satisfies
\begin{equation}\label{5.44}
	\left\{
	\begin{aligned}
		&\partial_{t}\bar{u}_{1,\sp}^{\pm}+y\partial_{x_3}u_{0,3}^{\pm}\partial_{y}\bar{u}_{1,\sp}^{\pm}+(u_{0,\sp}^{\pm}\cdot \nabla_{\sp})\bar{u}_{1,\sp}^{\pm}+(\bar{u}_{1,\sp}^{\pm}\cdot \nabla_{\sp})u_{0,\sp}^{\pm}=\lambda\partial_{y}^2\bar{u}_{1,\sp}^{\pm},\\
		&\partial_{t}\bar{\t}_{1}^{\pm}+y\partial_{x_3}u_{0,3}^{\pm}\partial_{y}\bar{\theta}_{1}^{\pm}+(u_{0,\sp}^{\pm}\cdot \nabla_{\sp})\bar{\t}_{1}^{\pm}+(\bar{u}_{1,\sp},\nabla_{\sp})\t_{0}^{\pm}=\frac{2}{5}\k \partial_{y}^2\bar{\t}_{1}^{\pm},\\
		&\bar{u}_{1,3}^{\pm}=0,\quad \bar{\rho}_{1}^{\pm}+\bar{\theta}_{1}^{\pm}=0.
	\end{aligned}
	\right.
\end{equation}

Hence applying Lemma \ref{lem5.1}, we can establish the existence of smooth solution $(\rho_{1}^{\pm}, \bar{u}_{1}^{\pm},\bar{\t}_{1}^{\pm}) $ of \eqref{5.43}-\eqref{5.44} with
\begin{equation}\label{5.45}
	\begin{aligned}
		&\sum\limits_{j=0}^{\bar{s}_{1}}\sum\limits_{|\beta|+2\gamma=j}\Big\{\|\partial_{t}^{\gamma}\nabla_{\bar{x}}^{\beta}(\bar{\rho}_{1}^{\pm}, \bar{u}_{1,\sp}^{\pm},\bar{\t}_{1}^{\pm})(t)\|_{L_{l_{j}^1}^2}^2+\int_{0}^{\tau}\|\partial_{t}^{\gamma}\nabla_{\bar{x}}^{\beta}(\bar{\rho}_{1}^{\pm}, \bar{u}_{1,\sp}^{\pm},\bar{\t}_{1}^{\pm})\|_{L_{l_{j}^1}^2}^2dt\Big\}\\
		&\leq C(\tau, E_{\bar{s}_1+3})\Big\{\sum\limits_{j=0}^{\bar{s}_{1}}\sum\limits_{|\beta|+2\gamma=j}\|\partial_{t}^{\gamma}\nabla_{\bar{x}}^{\beta}(\bar{u}_{1,\sp}^{\pm},\bar{\t}_{1}^{\pm})(0)\|_{L_{l_{1}^2}}^2\\
		&\qquad +\sup_{t\in [0,\tau]}\sum\limits_{|\beta|+2\gamma\leq \bar{s}_{1}+2}\|\partial_{t}^{\gamma}\partial_{\sp}^{\beta}(\pa_{x_{3}}u_{0,\sp}^{\pm},\pa_{x_{3}}\t_{0}^{\pm})(t)\|_{L^2(\mathbb{T}^2)}^2\Big\}\\
		&\leq C\Big(\tau, E_{\bar{s}_{1}+4},\sum\limits_{j=0}^{\bar{s}_{1}}\sum\limits_{|\beta|+2\gamma}\|\partial_{t}^{\gamma}\nabla_{\bar{x}}^{\beta}(\bar{u}_{1,\sp}^{\pm},\bar{\t}_{1}^{\pm})(0)\|_{L_{l_{j}^1}^2}^2\Big),
	\end{aligned}
\end{equation}
where we have denoted $l_{j}^1=\bar{l}_{1}+2(\bar{s}_{1}-j)$ with $\bar{l}_{1}\gg 1$ and $s_{0}\geq 4+\bar{s}_{1}$. Recalling \eqref{3.3-1}, we get from \eqref{5.45} that
\begin{equation}\label{5.46}
	\begin{aligned}
		&\sup_{t\in [0,\tau]}\sum\limits_{j=0}^{\bar{s}_1}\sum\limits_{|\beta|+2\gamma=j}\|\fw_{\bar{\k}_1}\partial_{t}^{\gamma}\nabla_{\bar{x}}\bar{f}_{1}^{\pm}(t)\|_{L_{l_{j}^1}^2L_{v}^{\infty}}\\
		&\leq C\Big(\tau,\sum\limits_{i=0}^1\sum\limits_{\gamma+|\beta|\leq s_{i}}\|\partial_{t}^{\gamma}\nabla_{x}^{\beta}(\rho_{0},u_{0},\theta_{0})(0)\|_{L_{x}^2}, \sum\limits_{j=0}^{\bar{s}_1}\sum\limits_{|\beta|+2\gamma=j}\|\partial_{t}^{\gamma}\nabla_{\bar{x}}^{\beta}(\bar{u}_{1,\sp}^{\pm},\bar{\t}_{1}^{\pm})(0)\|_{L_{l_{j}^1}^2}\Big).
	\end{aligned}
\end{equation}

{\it{Step 1.3. Construction of $\hat{f}_{1}^{\pm}$.}} As mentioned in Remark \ref{rem4.2}, we know that
$
\hat{f}_{1}^{\pm}\equiv 0.
$
This is reasonable since the Knudsen layer is used to mend the boundary condition at higher orders.

\noindent{\it{Step 2. Construction of solutions $f_{k},\bar{f}_{k}^{\pm}$ and $\hat{f}_{k}^{\pm}$.}} We shall use induction argument. Suppose we have already proved the existence of $f_{i},\bar{f}_{i}^{\pm}$ and $\hat{f}_{i}^{\pm}$ for $1\leq i\leq k$ such that
\begin{equation}\nonumber
	\begin{aligned}
		D_{k}+\bar{D}_{k}^{\pm}+\hat{D}_{k}^{\pm}\leq C\Big(\tau, &\sum\limits_{i=0}^{k}\Big[\sum\limits_{\gamma+|\beta|\leq s_{i}}\|\partial_{t}^{\gamma}\nabla_{x}^{\beta}(\rho_{i},u_{i},\theta_{i})(0)\|_{L_{x}^2}+\sum\limits_{j=0}^{\bar{s}_{i}}\sum\limits_{|\beta|+2\gamma=j}\|\partial_{t}^{\gamma}\nabla_{\bar{x}}^{\beta}(\bar{u}_{i,\sp}^{\pm},\bar{\t}_{i}^{\pm})(0)\|_{L_{l_{j}^{i}}^2}\Big]\Big)
	\end{aligned}
\end{equation}
where
\begin{align*}
	&D_{k}:=\sup_{t\in [0,\tau]}\Big\{\sum\limits_{i=1}^{k}\sum\limits_{\gamma+|\beta|\leq s_{i}}(\|\tw_{\k_{i}}\partial_{t}^{\gamma}\nabla_{x}^{\beta}f_{i}\|_{L_{x}^2L_{v}^{\infty}}+\|\pa_{t}^{\gamma}\nabla_{x}^{\beta}p_{i}\|_{L_{x}^2}) +\sum\limits_{\gamma+|\beta|\leq s_{k}}\|\pa_{t}^{\gamma}\nabla_{x}^{\beta}(\nabla_{x}p_{k+1})\|_{L_{x}^2}\Big\},\\
	&\bar{D}_{k}^{\pm}:=\sup_{t\in [0,\tau]}\Big\{\sum\limits_{i=1}^{k}\sum\limits_{j=0}^{\bar{s}_{i}}\sum\limits_{2\gamma+|\beta|=j}\|\fw_{\k_{i}}\partial_{t}^{\gamma}\nabla_{\bar{x}}^{\beta}\bar{f}_{i}^{\pm}\|_{L_{l_j^i}^2L_{v}^{\infty}}\Big\},\\
	&\hat{D}_{k}^{\pm}:=\sup_{t\in [0,\tau]}\Big\{\sum\limits_{i=1}^{k}\sum\limits_{\gamma+|\beta|\leq \hat{s}_{i}}\|e^{\z_{i}\eta}\fw_{\hat{\kappa}_{i}}\partial_{t}^{\gamma}\nabla_{\sp}^{\beta}\hat{f}_{i}^{\pm}\|_{L_{\hat{x},v}^{\infty}}\Big\}.
\end{align*}
for some $s_{i}>\bar{s}_{i}>\hat{s}_{i}\geq s_{i+1}>\bar{s}_{i+1}>\hat{s}_{i+1}\gg 1$, $\k_{i}\gg \bar{\k}_{i}\gg \hat{\k}_{i}\gg \k_{i+1}\gg \bar{\k}_{i+1}\gg \hat{\k}_{i+1}\gg 1$ with $1\leq i\leq k-1$, and $l_{j}^{i}=\bar{l}_{i}+2(\bar{s}_{i}-j), \bar{l}_{i}\gg 1$ with $1\leq i\leq k$ and $0\leq j\leq \bar{s}_{i}$.

{\it{Step 2.1. Construction of $f_{k+1}$.}} For the microscopic part of $f_{k+1}$, noting $\mathbf{L}^{-1}$ preserves the decay property of $v$, we obtain from \eqref{2.0} and Sobolev inequality that
\begin{equation}\label{5.48}
	\sum\limits_{\gamma+|\beta|\leq s_{k+1}}\|\tw_{\k_{k+1}}\partial_{t}^{\gamma}\nabla_{x}^{\beta}(\mathbf{I-P})f_{k+1}(t)\|_{L_{x}^2L_{v}^{\infty}}\leq C(D_{k}), \quad \text{for }s_{k}\geq s_{k+1}+1. 
\end{equation}

For the macroscopic part of $f_{k+1}$. Firstly, By using Lemma \ref{lem2.2} and similar arguments as in \eqref{5.40-2}, we obtain that for $s_{k}\geq s_{k+1}$ and $\bar{s}_{k},\hat{s}_{k-1}\geq s_{k+1}+1$,
\begin{equation}\label{5.50-4}
	\begin{aligned}
	&\sup_{t\in [0,\tau]}\sum\limits_{\gamma+|\beta|\leq s_{k+1}}\|\pa_{t}^{\g}\nabla_{x}^{\beta}p_{k+1}\|_{L_{x}^2}\\
	&\leq C\big(\tau, D_{k}, \bar{D}_{k}^{\pm},\hat{D}_{k-1}^{\pm},\sum\limits_{|\beta|+\gamma\leq s_{k+1}}\|\pa_{t}^{\gamma}\nabla_{x}^{\beta}(\rho_{k+1},u_{k+1},\t_{k+1})(0)\|_{L_{x}^2}\big).
	\end{aligned}
\end{equation}
Then, using \eqref{5.50-4}, \eqref{4.33}-\eqref{4.34}, Lemmas \ref{lem2.1}, \ref{lem4.4} and similar arguments in \cite[Proposition 5.1]{Guo-Huang-Wang}, we can establish the existence of  $(\rho_{k+1},u_{k+1},\theta_{k+1})$ satisfying
\begin{equation}\label{5.57}
	\begin{aligned}
		&\sup_{t\in [0,\tau]}\sum\limits_{|\beta|+\gamma\leq s_{k+1}}\|\pa_{t}^{\g}\pa_{x}^{\beta}(\rho_{k+1},u_{k+1},\t_{k+1})(t)\|_{L_{x}^2}+\sup_{t\in [0,\tau]}\sum\limits_{|\beta|+\gamma\leq s_{k+1}}\|\pa_{t}^{\g}\nabla_{x}^{\beta}(\nabla_{x}p_{k+2})(t)\|_{L_{x}^2}\\
		&\leq C\Big(\tau, D_{k}+\bar{D}_{k}^{\pm}+\hat{D}_{k-1}^{\pm}, \sum\limits_{|\beta|+\gamma\leq s_{k+1}}\|\pa_{t}^{\gamma}\nabla_{x}^{\beta}(\rho_{k+1},u_{k+1},\t_{k+1})(0)\|_{L_{x}^2}\Big),
	\end{aligned}
\end{equation}
which, together with \eqref{5.48}, yields that
\begin{equation}\label{5.58}
	\begin{aligned}
		&\sup_{t\in [0,\tau]}\sum\limits_{|\beta|+\gamma\leq s_{k+1}}(\|\tw_{\k_{k+1}}\pa_{t}^{\gamma}\nabla_{x}^{\beta}f_{k+1}(t)\|_{L_{x}^2L_{v}^{\infty}}+\sum\limits_{|\beta|+\gamma\leq s_{k+1}}\|\pa_{t}^{\g}\nabla_{x}^{\beta}(\nabla_{x} p_{k+2})(t)\|_{L_{x}^2})\\
		&\leq C\Big(\tau, D_{k}+\bar{D}_{k}^{\pm}+\hat{D}_{k-1}^{\pm},\sum\limits_{|\beta|+\gamma\leq s_{k+1}}\|\pa_{t}^{\gamma}\nabla_{x}^{\beta}(\rho_{k+1},u_{k+1},\t_{k+1})(0)\|_{L_{x}^2}\Big),
	\end{aligned}
\end{equation}
provided $l_{j}^{i}\geq l_{j}^{i+1}$, $\bar{l}_{k}\gg 1$, $s_{k}\geq s_{k+1}+2\geq 5$, $\bar{s}_{k}\geq 2s_{k+1}+4$ and $\hat{s}_{k-2}>s_{k+1}+2$.

Since we have obtained $f_{i}$ with $1\leq i\leq k+1$, by using \eqref{2.0}, it holds that
\begin{equation}\label{5.59}
	\begin{aligned}
		&\sup_{t\in [0,\tau]}\sum\limits_{|\beta|+\gamma\leq s_{k+1}-3}\|\tw_{\k_{k+2}}\pa_{t}^{\gamma}\nabla_{x}^{\beta}(\mathbf{I-P})f_{k+2}(t)\|_{L_{x}^2L_{v}^{\infty}}\\
		&\leq C\Big(\tau, D_{k}+\bar{D}_{k}^{\pm}+\hat{D}_{k-1}^{\pm},\sum\limits_{|\beta|+\gamma\leq s_{k+1}}\|\pa_{t}^{\gamma}\nabla_{x}^{\beta}(\rho_{k+1},u_{k+1},\t_{k+1})(0)\|_{L_{x}^2}\Big),
	\end{aligned}
\end{equation}
which will be used when we consider the boundary conditions \eqref{4.37} and the source terms in the equations \eqref{3.12} of  $(\bar{u}_{k+1,\sp}^{\pm},\bar{\t}_{k+1}^{\pm})$ in the following.

{\it{Step 2.2. Construction of $\bar{f}_{k+1}^{\pm}$.}} For the microscopic part of $f_{k+1}$, it follows from \eqref{3.15}, \eqref{3.18} and similar arguments in \cite[Proposition 5.1]{Guo-Huang-Wang} that


\begin{equation}\label{5.62}
	\begin{aligned}
		&\sum\limits_{2\g+|\beta|\leq j}\|\fw_{\bar{\k}_{k+1}}\pa_{t}^{\gamma}\nabla_{\bar{x}}(\mathbf{I-P})\bar{f}_{k+1}^{\pm}\|_{L_{l}^2L_{v}^{\infty}}\\
		&\leq C\Big(\sum\limits_{i=0}^{k-1}\sum\limits_{2\gamma+|\beta|\leq j+\fb+2}\|\tw_{\k_{i}}\pa_{t}^{\gamma}\nabla_{x}^{\beta}f_{i}(t)\|_{L_{x}^2L_{v}^{\infty}}+\sum\limits_{i=1}^{k}\sum\limits_{2\g+|\beta|\leq j+2}\|\fw_{\bar{\k}_{i}}\pa_{t}^{\g}\nabla_{\bar{x}}^{\beta}\bar{f}_{i}^{\pm}(t)\|_{L_{l+2\fb}^2L_{v}^{\infty}}\Big).
	\end{aligned}
\end{equation}


For the macroscopic part of $\bar{f}_{k+1}^{\pm}$, recalling \eqref{3.13}, from the first glance, we need to know the information on $\mathbf{P}\bar{f}_{k+1}$, which is obviously an unknown quantity. Fortunately, due to the explicit formulas of the source term $\mathfrak{f}_{k,\sp}^{\pm},\mathfrak{g}_{k}^{\pm}$ \eqref{3.12} and boundary conditions  \eqref{4.37}, we only need to estimate $\la \mathcal{A}_{3i},J_{k}^{\pm}\ra (i=1,2) $ and $\la \mathcal{B}_{3},J_{k}^{\pm}\ra$, which depend only on $\bar{u}_{k+1,3}^{\pm}$ and $f_{i}(0\leq i\leq k+1), \bar{f}_{j}^{\pm}(1\leq j\leq k)$ as explained in the last two paragraphs in the proof of Proposition \ref{prop3.1}. Hence, we can use \eqref{2.13}, \eqref{4.37}, \eqref{5.59}, Lemma \ref{lem5.1} and similar arguments in \cite[Proposition 5.1]{Guo-Huang-Wang} to establish the existence of $(\bar{\rho}_{k+1}^{\pm},\bar{u}_{k+1}^{\pm},\bar{\theta}_{k+1}^{\pm})$, and moreover it holds that

\begin{align}\label{5.74}
	&\sup_{t\in [0,\tau]}\sum_{j=0}^{\bar{s}_{k+1}} \sum_{j=2\gamma+|\beta|} \|\fw_{\hat{\k}_{k+1}}\pa_t^{\gamma} \nabla_{\bar{x}}^{\beta} \bar{f}_{k+1}^{\pm}(t)\|^2_{L^2_{l_j^{k+1}}L^\infty_v} \nonumber\\
	&\leq C\Big(\tau, D_k+\bar{D}_k^{\pm}+\hat{D}_k^{\pm},  \sum_{j=0}^{\bar{s}_{k+1}} \sum_{j=2\gamma+|\beta|}   \|\pa_t^{\gamma} \nabla_{\bar{x}}^{\beta} (\bar{u}_{k+1,\sp}^{\pm},\bar{\theta}_{k+1}^{\pm})(0)\|^2_{L^2_{l_j^{k+1}}}\nonumber\\
	&\qquad+\sum_{\gamma+|\beta|\leq s_{k+1}}\|\pa_t^\gamma\nabla_x^{\beta}(\rho_{k+1}, u_{k+1},  \theta_{k+1})(0)\|_{L^2_{x}}\Big)<\infty,
\end{align}
provided that $\bar{s}_{k}\geq \bar{s}_{k+1}+6$, $s_{k}>\bar{s}_{k+1}+8+\fb$ and $l_{j}^{i}\geq 2l_{j}^{i+1}+18+2\fb$ for $1\leq i\leq k$ with $l_{j}^{k+1}=\bar{l}_{k+1}+2(\bar{s}_{k+1}-j)$ and $\bar{l}_{k+1}\gg 1, l_{j}^{k+1}\geq 2\fb+6$ .

{\it{Step 2.3. Construction of $\hat{f}_{k+1}^{\pm}$.}} By using \eqref{4.3}--\eqref{4.3-1}, Lemmas \ref{lem4.1}--\ref{lem4.2}, Remark \ref{rem4.4} and similar arguments in \cite[Proposition 5.1]{Guo-Huang-Wang}, we can establish the existence of $\hat{f}_{k+1}^{\pm}$ and obtain

\begin{align*}
	&\sup_{t\in [0,\tau]} \sum_{|\beta|+\gamma\leq \hat{s}_{k+1}} \Big\{\|\fw_{\hat{\k}_{k+1}}e^{\zeta_{k+1}\eta}\partial_t^{\gamma}\nabla^{\beta}_{\sp} \hat{f}_{k+1}^{\pm}(t)\|_{L^\infty_{x_{\sp} ,\eta,v}\cap L^2_{x_{\sp} }L^\infty_{\eta,v}}+\|\fw_{\hat{\k}_{k+1}}\partial_t^{\gamma}\nabla^{\beta}_{\sp} \hat{f}_{k+1}^{\pm,0}(t)\|_{L^\infty_{x_{\sp} ,v}\cap L^2_{x_{\sp} }L^\infty_{v}}\Big\}\nonumber\\
	&\leq  C\Big(\tau, D_{k}+\bar{D}_{k}^{\pm}+\hat{D}_{k}^{\pm},  \sum_{j=0}^{\bar{s}_{k+1}} \sum_{j=2\gamma+|\beta|}   \|\pa_t^{\gamma} \nabla_{\bar{x}}^{\beta} (\bar{u}_{k+1,\sp}^{\pm},\bar{\theta}_{k+1}^{\pm})(0)\|^2_{L^2_{l_j^{k+1}}}\nonumber\\
	&\qquad+\sum_{\gamma+|\beta|\leq s_{k+1}}\|\pa_t^\gamma\nabla_x^{\beta}(\rho_{k+1}, u_{k+1},  \theta_{k+1})(0)\|_{L^2_{x}}\Big),
\end{align*}
provided that $\hat{s}_{k}\geq 2+\hat{s}_{k+1}$, $\bar{s}_{k}\geq 2(\hat{s}_{k+1}+\fb+2)$ and  $\zeta_{k}\geq 2\zeta_{k+1}$.

\noindent{\it{Step 3.}} By induction, we establish the existence of solutions $f_i, \bar{f}_i, \hat{f}_i, \, i=1,\cdots, N$ with
\begin{align*}
	D_N+\bar{D}_{N}^{\pm}+\hat{D}_{N}^{\pm}&\leq C\Big(\tau, \sum_{i=0}^{N}\sum_{j=0}^{\bar{s}_{i}} \sum_{j=2\gamma+|\beta|}   \|\pa_t^{\gamma} \nabla_{\bar{x}}^{\beta} (\bar{u}_{i,\sp}^{\pm},\bar{\theta}_{i}^{\pm})(0)\|^2_{L^2_{l_j^{i}}}\nonumber\\
	&\qquad\qquad+\sum_{i=0}^{N}\sum_{\gamma+|\beta|\leq s_{i}}\|\pa_t^\gamma \nabla_x^{\beta}(\rho_{i}, u_{i},  \theta_{i})(0)\|_{L^2_{x}}\Big),
\end{align*}
where we have chosen $s_i,\,\,\bar{s}_i,\,\,\hat{s}_i,\,\,\zeta_{i}$ satisfying 
\begin{align}\label{5.84}
		& s_0>s_i>\bar{s}_i>\hat{s}_i> s_{i+1}>\bar{s}_{i+1}>\hat{s}_{i+1}> \cdots \geq 3,\,0<2\zeta_{i+1}\leq \zeta_{i} ,\,\,\,\,\,\mbox{for}\ i=1,\cdots,  N-1;\nonumber\\
		&s_1\leq s_0-2,\,s_{i+1}\leq s_i-2,\,\bar{s}_{i+1}\leq \bar{s}_{i}-6,\,\hat{s}_{i+1}\leq \hat{s}_{i}-2, \quad\qquad\quad\quad\mbox{for}\,i=1,\cdots,  N-1;\nonumber\\
		&s_{i+1}+2\leq \min\{\frac{1}{2}\bar{s}_{i},\hat{s}_i\},\,\bar{s}_{i+1}\leq s_{i}-8-\fb,\,2\hat{s}_{i+1}\leq  \bar{s}_{i}-4-2\fb,\,\quad\,\,\mbox{for}\,i=1,\cdots, N-1;
\end{align}
and taken $l_j^i=\bar{l}_j+2(\bar{s}_i-j)$ with $0\leq j\leq \bar{s}_i$ so that
\begin{equation}\label{5.85}
	l^{N}_{j}\geq 2\fb+6 \quad \mbox{and} \quad  l_j^i\geq 2 l_j^{i+1} +18+2\fb,\quad \mbox{for}\, 1\leq i\leq  N-1.
\end{equation}
Here we don't give the precise relations between the velocity weight functions $\k_i, \bar{\k}_i$ and $\hat{\k}_i$, but we require
\begin{equation}\label{5.86}
	\k_i\gg \bar{\k}_i\gg  \hat{\k}_i\gg \k_{i+1}\gg \bar{\k}_{i+1}\gg  \hat{\k}_{i+1}\gg1,
\end{equation}
since the functions  $F_k, \bar{F}_k^{\pm}$ and $\hat{F}_k^{\pm}$ indeed decay exponentially with respect to particle velocity $v$. Moreover, we mention that $f_k, \bar{f}_k^{\pm}$ are smooth, but $\hat{f}_{k}^{\pm}$ is only continuous away from the grazing set $[0,\tau]\times \gamma_0^{i}$ with $i=0,1$. Therefore the proof of Proposition \ref{prop5.1} is complete. $\hfill\Box$

\section{$L^2-L^{\infty}$ method: proof of Theorem \ref{theorem1.1}}
In this section, we apply $L^2-L^{\infty}$ framework to estimate the remainder term $F_{R}^{\v}$ in \eqref{1.22-2} over $[0,\tau]\times \mathbb{T}^2\times (0,1)$. Recalling the boundary conditions in Section \ref{sec2.4}, it is easy to check that $F_{R}^{\v}$ satisfies the specular reflection boundary condition, i.e.,
\begin{equation}\label{7.1}
	F_{R}^{\v}(t,x_{\sp},x_{3},v_{\sp},v_{\sp},v_{3})\vert_{\gamma_{-}^{i}}=F_{R}^{\v}(t,x_{\sp},i,v_{\sp},-v_{3})\quad \text{with } i=0,1.
\end{equation}
\subsection{$L^2$-energy estimate for remainder} Noting \eqref{1.23}, \eqref{1.31} and \eqref{1.32}, we can rewrite the equation of $f_{R}^{\v}$ as
\begin{equation}\label{7.2}
	\begin{aligned}
		&\partial_{t}f_{R}^{\v}+\frac{1}{\v^{n-2}}v\cdot \nabla_{x} f_{R}^{\v}-\frac{1}{\v^{2n-2}}\frac{1}{\mathcal{M}_{\v}}[Q(\mathcal{M}_{\v},\sqrt{\mathcal{M}_{\v}}f_{R}^{\v})+Q(\sqrt{\mathcal{M}_{\v}}f_{R}^{\v},\mathcal{M}_{\v})]\\
		&=-\frac{1}{\sqrt{\mathcal{M}_{\v}}}[\partial_{t}\sqrt{\mathcal{M}_{\v}}+\v^{2-n}v\cdot \nabla_{x}\sqrt{\mathcal{M}_{\v}}]f_{R}^{\v}-\frac{1}{\v^{2}}\frac{1}{\sqrt{\mathcal{M}_{\v}}}[Q(r_0,\sqrt{\mathcal{M}_{\v}}f_{R}^{\v})+Q(\sqrt{\mathcal{M}_{\v}}f_{R}^{\v},r_{0})]\\
		&\quad +\sum\limits_{i=1}^N\v^{i-n}\frac{1}{\sqrt{\mathcal{M}_{\v}}}\Big[Q(\sqrt{\mathcal{M}_{\v}}f_{R}^{\v},F_{i}+\sum\limits_{\pm}[\phi(\v y)\bar{F}_{i}^{\pm}+\phi(\v^{n}\eta)\hat{F}_{i}^{\pm}])\\
		&\qquad \qquad +Q(F_{i}+\sum\limits_{\pm}[\phi(\v y)\bar{F}_{i}^{\pm}+\phi(\v^{n}\eta)\hat{F}_{i}^{\pm}],\sqrt{\mathcal{M}_{\v}}f_{R}^{\v})\Big]\\
		&\quad +\v^{k_0-n}\frac{1}{\sqrt{\mathcal{M}_{\v}}}[Q(\sqrt{\mathcal{M}_{\v}}f_{R}^{\v},\sqrt{\mathcal{M}_{\v}}f_{R}^{\v})]+\frac{1}{\v^{2n-4+k_0}}\frac{1}{\sqrt{\mathcal{M}_{\v}}}(R^{\v}+\sum\limits_{\pm}(\bar{R}^{\pm,\v}+\hat{R}^{\pm,\v})),
	\end{aligned}
\end{equation}
where $R^{\v}, \bar{R}^{\pm,\v},\hat{R}^{\pm,\v}$ are the ones defined in \eqref{1.24}-\eqref{1.28}. From \eqref{7.1}, we know $f_{R}^{\v}$ satisfies specular reflection boundary conditions:
\begin{equation}\label{7.3}
	f_{R}^{\v}(t,x_{\sp},0,v_{\sp},v_{3})\vert_{v_{3}>0}=f_{R}^{\v}(t,x_{\sp},0,v_{\sp},-v_{3}),\quad f_{R}^{\v}(t,x_{\sp},1,v_{\sp},v_{3})\vert_{v_{3}<0}=f_{R}^{\v}(t,x_{\sp},1,v_{\sp},-v_{3})
\end{equation}

\begin{lemma}\label{lem7.1}
	Assume $k_0\geq 1$, $\beta>\frac{9}{2}+\frac{k_0-1}{n-1}, N\geq 2n-3+k_0$ and $\fb\geq n+k_0-2.5$. Let $\tau>0$ be the lifespan of incompressible Euler equation \eqref{2.8-1}, then there exists a suitably small constant $\v_0>0$ such that for all $\v\in (0,\v_0)$ and $t\in [0,\tau]$, it holds that
	\begin{equation}\label{7.9}
		\frac{d}{dt}\|f_{R}^{\v}(t)\|_{L^2}^2+\frac{\tilde{c}_0}{2\v^{2n-2}}\|(\mathbf{I}-\mathbf{P}_{\v})f_{R}^{\v}(t)\|_{\nu}^2\leq C(\tau)\{1+\v^{2k_0-2}\|h_{R}^{\v}(t)\|_{L^{\infty}}^2\}\cdot (1+\|f_{R}^{\v}\|_{L^2}^2).
	\end{equation}
\end{lemma}
\textbf{Proof.} We define the linearized collision operator $\mathbf{L}_{\v}$ by
$$
\mathbf{L}_{\v}g:=-\frac{1}{\sqrt{\mathcal{M}_{\v}}}[Q(\mathcal{M}_{\v},\sqrt{\mathcal{M}_{\v}}g)+Q(\sqrt{\mathcal{M}_{\v}}g,\mathcal{M}_{\v})],
$$
then, similar as \eqref{1.7-2}, there exists two positive numbers $\hat{c}_0, \tilde{c}_{0}>0$ such that for any function $g$
$$
\langle\mathbf{L}_{\v}g,g\rangle\geq \hat{c}_0\|(\mathbf{I}-\mathbf{P}_{\v})g\|_{\nu_{\v}}^2\geq \tilde{c}_0\|(\mathbf{I}-\mathbf{P}_{\v})g\|_{\nu}^2,
$$
where $\mathbf{P}_{\v}$ is the $L_{v}^2$ projection with respect to the orthogonal bases of null space of $\mathbf{L}_{\v}$ and
\begin{equation}\label{7.9-1}
\nu_{\v}:=\int_{\R^3}\int_{\mathbb{S}^2}B(v-u,\theta)\mathcal{M}_{\v}(u)dud\omega\cong (1+|v|)\cong\nu(v).
\end{equation}
Multiplying \eqref{7.2} by $f_{R}^{\v}$ and integrating the resultant equation over $\Omega\times \R^3$, one obtains that
\begin{equation}\label{7.10}
	\begin{aligned}
		&\frac{1}{2}\frac{d}{dt}\|f_{R}^{\v}\|_{L^2}^2+\frac{\tilde{c}_0}{\v^{2n-2}}\|(\mathbf{I}-\mathbf{P}_{\v})f_{R}^{\v}\|_{\nu}^2\\
		&\quad +\frac{1}{2}\Big(\int_{\mathbb{T}^2}\int_{\R^3}v_{3}|f_{R}^{\v}(t,x_{\sp},1,v)|^2dvdx_{\sp}-\int_{\mathbb{T}^2}\int_{\R^3}v_{3}|f_{R}^{\v}(t,x_{\sp},0,v)|^2dvdx_{\sp}\Big)\\
		&\leq -\int_{\Omega}\int_{\R^3}\frac{1}{\sqrt{\mathcal{M}_{\v}}}[\partial_{t}\sqrt{\mathcal{M}_{\v}}+\v^{2-n}v\cdot \nabla_{x}\sqrt{\mathcal{M}_{\v}}]|f_{R}^{\v}|^2dvdx\\
		&\quad -\int_{\Omega}\int_{\R^3}\frac{1}{\v^{2}}\frac{1}{\sqrt{\mathcal{M}_{\v}}}[Q(r_0,\sqrt{\mathcal{M}_{\v}}f_{R}^{\v})+Q(\sqrt{\mathcal{M}_{\v}}f_{R}^{\v},r_{0})]f_{R}^{\v}dvdx\\
		&\quad+\int_{\Omega}\int_{\R^3}\v^{k_0-n}\frac{1}{\sqrt{\mathcal{M}_{\v}}}[Q(\sqrt{\mathcal{M}_{\v}}f_{R}^{\v},\sqrt{\mathcal{M}_{\v}}f_{R}^{\v})]f_{R}^{\v}dvdx\\
		&\quad+\int_{\Omega}\int_{\R^3}\sum\limits_{i=1}^{N}\v^{i-n}\frac{1}{\sqrt{\mathcal{M}_{\v}}}\Big[Q(\sqrt{\mathcal{M}_{\v}}f_{R}^{\v},F_{i}+\sum\limits_{\pm}[\phi(\v y)\bar{F}_{i}^{\pm}+\phi(\v^n\eta)\hat{F}_{i}^{\pm}])\\
		&\qquad \qquad +Q(F_{i}+\sum\limits_{\pm}[\phi(\v y)\bar{F}_{i}^{\pm}+\phi(\v^n\eta)\hat{F}_{i}^{\pm}],\sqrt{\mathcal{M}_{\v}}f_{R}^{\v})\Big]f_{R}^{\v}dvdx\\
		&\quad +\int_{\Omega}\int_{\R^3}\frac{1}{\v^{2n-4+k_0}}\frac{1}{\sqrt{\mathcal{M}_{\v}}}(R^{\v}+\sum\limits_{\pm}(\bar{R}^{\pm,\v}+\hat{R}^{\pm,\v}))f_{R}^{\v}dvdx:=\sum\limits_{j=1}^{5}I_{j}.
	\end{aligned}
\end{equation}
Using \eqref{7.3}, it is clear that
\begin{equation}\label{7.11}
	\int_{\mathbb{T}^2}\int_{\R^3}v_{3}|f_{R}^{\v}(t,x_{\sp},i,v)|^2dvdx_{\sp}=0,\quad \text{for }i=0,1.
\end{equation}
For any $\delta>0$, using similar arguments as in \cite{Guo Jang Jiang-1}, $I_1$ can be decomposed as
\begin{equation}\label{7.12}
	\begin{aligned}
		I_{1}=&-\int_{\Omega}\int_{|v|\geq \frac{\delta}{\v^{n-1}}}\frac{1}{\sqrt{\mathcal{M}_{\v}}}[\partial_{t}\sqrt{\mathcal{M}_{\v}}+\v^{2-n}v\cdot \nabla_{x}\sqrt{\mathcal{M}_{\v}}]|f_{R}^{\v}|^2dvdx
		\\&-\int_{\Omega}\int_{|v|\leq \frac{\delta}{\v^{n-1}}}\frac{1}{\sqrt{\mathcal{M}_{\v}}}[\partial_{t}\sqrt{\mathcal{M}_{\v}}+\v^{2-n}v\cdot \nabla_{x}\sqrt{\mathcal{M}_{\v}}]|f_{R}^{\v}|^2dvdx\\
		:=&I_{1,1}+I_{1,2}.
	\end{aligned}
\end{equation}
It follows from \eqref{1.32}, \eqref{1.35}-\eqref{1.36} that $f_{R}^{\v}\leq C(1+|v|^2)^{-\frac{\beta}{2}}h_{R}^{\v}$, then it holds that
\begin{equation}\label{7.13}
	\begin{aligned}
		|I_{1,1}|&\leq C\|\nabla_{x}({\rho_0},u_{0},\theta_{0})\|_{L^2}\cdot\|(1+|v|^2)^{\frac{\beta}{2}}f_{R}^{\v}\|_{L^{\infty}}\cdot \|f_{R}^{\v}\|_{L^2}\cdot \Big(\int_{|v|\geq \frac{\delta}{\v^{n-1}}}(1+|v|^2)^{3-\beta}\,dv\Big)^{\frac{1}{2}}\\
		&\leq C_{\delta}\v^{(n-1)(2\beta-9)}\|h_{R}\|_{\infty}\|f_{R}^{\v}\|_{L^2},
	\end{aligned}
\end{equation}
and
\begin{equation}\label{7.14}
	\begin{aligned}
		|I_{1,2}|&\leq C\|\nabla_{x}({\rho_0},u_{0},\theta_{0})\|_{L^{\infty}}\cdot\|\mathbf{1}_{|v|\leq \frac{\delta}{\v^{n-1}}}(1+|v|^2)^{\frac{3}{4}}f_{R}^{\v}\|_{L^2}^2\\
		&\leq C\|\mathbf{1}_{|v|\leq \frac{\delta}{\v^{n-1}}}(1+|v|^2)^{\frac{3}{4}}\mathbf{P}_{\v}f_{R}^{\v}\|_{L^2}^2+C\|\mathbf{1}_{|v|\leq \frac{\delta}{\v^{n-1}}}(1+|v|^2)^{\frac{3}{4}}(\mathbf{I}-\mathbf{P}_{\v})f_{R}^{\v}\|_{L^2}^2\\
		&\leq \frac{C\delta^2}{\v^{2n-2}}\|(\mathbf{I-P})f_{R}^{\v}\|_{\nu}^2+C\|f_{R}^{\v}\|_{L^2}^2.
	\end{aligned}
\end{equation}
Hence, substituting \eqref{7.13}-\eqref{7.14} into \eqref{7.12}, one obtains
\begin{equation}\label{7.15}
	|I_1|\leq |I_{1,1}|+|I_{1,2}|\leq C_{\delta}(\v^{(n-1)(\beta-3)}\|h_{R}\|_{L^{\infty}}\|f_{R}^{\v}\|_{L^2}+\|f_{R}^{\v}\|_{L^2}^2)+\frac{C\delta^2}{\v^{2n-2}}\|(\mathbf{I}-\mathbf{P}_{\v})f_{R}^{\v}\|_{\nu}^2.
\end{equation}
For $I_2$, noting the definition of $r_{0}$ in \eqref{1.31}, we have
$$
|r_{0}|\leq Cp(v)\mu_{M}^{\epsilon}\quad \text{for some } \epsilon\in (\frac{1}{2},1),
$$
where $p(v)$ is a polynomial with respect to $v$, thus a similar calculation in \cite{Guo-Huang-Wang} shows that
\begin{equation}\label{7.16}
	\begin{aligned}
		|I_2|&\leq \frac{C}{\v^2}\|(\mathbf{I}-\mathbf{P}_{\v})f_{R}^{\v}\|_{\nu}\|f_{R}^{\v}\|_{\nu}\cdot\Big(\int_{\R^3}(1+|v|)\Big(\frac{r_{0}}{\sqrt{\mathcal{M}_{\v}}}\Big)^2dv\Big)^{\frac{1}{2}}\\
		&\leq \frac{\delta+\v^{2n-4}}{\v^{2n-2}}\|(\mathbf{I}-\mathbf{P}_{\v})f_{R}^{\v}\|_{\nu}^2+C_{\delta}\v^{2n-6}\|f_{R}^{\v}\|_{L^2}^2.
	\end{aligned}
\end{equation}

For $I_3$, using the similar calculations in \cite{Guo Jang Jiang-1}, we have
\begin{equation}\label{7.17}
	\begin{aligned}
		I_{3}&=\v^{k_0-n}\int_{\Omega}\int_{\R^3}\frac{1}{\sqrt{\mathcal{M}_{\v}}}\Big[Q(\sqrt{\mathcal{M}_{\v}}f_{R}^{\v},\sqrt{\mathcal{M}_{\v}}f_{R}^{\v})\Big](\mathbf{I}-\mathbf{P}_{\v})f_{R}^{\v}dvdx\\
		&\leq C\v^{k_0-n}\|(\mathbf{I}-\mathbf{P}_{\v})f_{R}^{\v}\|_{\nu}\|h_{R}^{\v}\|_{L^{\infty}}\|f_{R}^{\v}\|_{L^2}\\
		&\leq \frac{\delta}{\v^{2n-2}}\|(\mathbf{I}-\mathbf{P}_{\v})f_{R}^{\v}\|_{\nu}^2+C_{\delta}\v^{2k_0-2}\|h_{R}^{\v}\|_{L^{\infty}}^2\|f_{R}^{\v}\|_{L^2}^2
	\end{aligned}
\end{equation}

Recall from \eqref{5.84} that
\begin{equation*}
	s_N>\bar{s}_N\geq 2\fb+4+\hat{s}_N,\quad \hat{s}_N\geq 1,
\end{equation*}
which, together with \eqref{5.0-1} \eqref{5.85}-\eqref{5.86} and Sobolev imbedding theorem, yields  that, for $1\leq i\leq N$ and $t\in[0,\tau]$,
\begin{align}\label{7.17-1}
	\begin{split}
		&\sum_{k=0}^{2\fb+2}\{\|\tw_{\k_i}(v)   \nabla^k_{t,x}f_i(t)\|_{L^2_{x,v}}+\|\tw_{\k_i}\nabla^k_{t,x} f_i(t)\|_{L^\infty_{x,v}}\}\leq C(\tau),\\
		&\sum_{k=0}^{\fb+2}\{\|\fw_{\bar{\k}_i} (1+y)^{\fb+9} \nabla^k_{t,\bar{x}}\bar{F}_{i}^{\pm}(t)\|_{L^2_{\bar{x},v}}+\|\fw_{\bar{\k}_i} (1+y)^{\fb+9} \nabla^k_{t,\bar{x}}\bar{F}_{i}^{\pm}(t)\|_{L^\infty_{\bar{x},v}} \}\leq C(\tau), \\
		&\sum_{k=0,1}\{\|\fw_{\hat{\k}_i} e^{\frac{1}{2^N}\eta} \nabla_{t,x_{\sp} }^{k}\hat{F}_{i}^{\pm}(t)\|_{L^2_{\hat{x},v}}+ \|\fw_{\hat{\k}_i} e^{\frac{1}{2^N}\eta} \nabla_{t,x_{\sp} }^{k}\hat{F}_{i}^{\pm}(t)\|_{L^\infty_{\hat{x},v}}\}\leq C(\tau),
	\end{split}
\end{align}
where $C(\tau)$ depending only on $\tau$ and the initial data,
\begin{align}
	C(\tau):=C\Big( &\tau,\sum_{i=0}^{N}\Big[\sum_{\gamma+|\beta|\leq s_{i}}\|\pa_t^\gamma\nabla_x^{\beta}(\rho_{i}, u_{i},  \theta_{i})(0)\|_{L^2_{x}}+\sum_{j=0}^{\bar{s}_{i}} \sum_{2\gamma+|\beta|=j}   \|\pa_t^{\gamma} \nabla_{\bar{x}}^{\beta} (\bar{u}_{i,\sp}^{\pm},\bar{\theta}_{i}^{\pm})(0)\|^2_{L^2_{l_j^{i}}}
	\Big]\Big).\nonumber
\end{align}

Recalling \eqref{1.30} and \eqref{1.35}, for $1\leq i\leq N$, it holds that
\begin{align}\label{7.17-2}
	\begin{split}
		|w_{\beta}(v) \frac{\sqrt{\mu}}{\sqrt{\mathcal{M}_{\v}}} {f}_{i}(t,x_{\sp} ,x_3,v)|
		&\leq C|w_{\beta}(v) \mu^{-\fa} {f}_i(t,x_{\sp} ,x_{3},v)|\cdot (\mu_{M})^{(\frac12+\fa)\alpha-\frac12},\\
		|w_{\beta}(v) \frac{\sqrt{\mu}}{\sqrt{\mathcal{M}_{\v}}} \bar{f}_{i}^{\pm}(t,x_{\sp} ,y,v)|
		&\leq C |w_{\beta}(v) \mu^{-\fa} \bar{f}_{i}^{\pm}(t,x_{\sp} ,y,v)|\cdot (\mu_{M})^{(\frac12+\fa)\alpha-\frac12},\\
		|w_{\beta}(v) \frac{\sqrt{\mu}}{\sqrt{\mathcal{M}_{\v}}} \hat{f}_{i}^{\pm}(t,x_{\sp} ,\eta,v)|
		&\leq C |w_{\beta}(v) \mu^{-\fa} \hat{f}_i(t,x_{\sp} ,\eta,v)|\cdot (\mu_{M})^{(\frac12+\fa)\alpha-\frac12}.
	\end{split}
\end{align}
Taking $0<\frac{1}{2\alpha}(1-\alpha)<\fa<\frac12$, we have
$ (\frac12+\fa)\alpha-\frac12>0$. It follows from \eqref{6.1}, \eqref{7.17-1}-\eqref{7.17-2} that
\begin{equation}\label{7.18}
	\begin{aligned}
		I_{4}&\leq C\|(\mathbf{I}-\mathbf{P}_{\v})f^{\v}_R \|_{\nu} \|f^\v_R\|_{\nu}\\
		&\quad \times \sum_{i=1}^N \v^{i-n}  \{\|w_{\beta} \frac{\sqrt{\mu}}{\sqrt{\mathcal{M}_{\v}}} {f}_{i}\|_{L^\infty_{x,v}}+\|w_{\beta} \frac{\sqrt{\mu}}{\sqrt{\mathcal{M}_{\v}}} \bar{f}_{i}^{\pm}\|_{L^\infty_{\bar{x},v}}+\|w_{\beta} \frac{\sqrt{\mu}}{\sqrt{\mathcal{M}_{\v}}} \hat{f}_{i}^{\pm}\|_{L^\infty_{\hat{x},v}}\}\\
		&\leq C\frac{1}{\v^{n-1}}\|(\mathbf{I}-\mathbf{P}_{\v})f_{R}^{\v}\|_{\nu}\|f_{R}^{\v}\|_{\nu}\leq \frac{C}{\v^{n-1}}\|(\mathbf{I}-\mathbf{P}_{\v})f_{R}^{\v}\|_{\nu}^2+\frac{C}{\v^{n-1}}\|(\mathbf{I}-\mathbf{P}_{\v})f_{R}^{\v}\|_{\nu}\cdot \|\mathbf{P}_{\v}\|_{\nu}\\
		&\leq (\delta+C\v^{n-1})\frac{1}{\v^{2n-2}}\|(\mathbf{I}-\mathbf{P}_{\v})f_{R}^{\v}\|_{\nu}^2+C_{\delta}\|f_{R}^{\v}\|_{L^2}^2.
	\end{aligned}
\end{equation}

Using \eqref{1.22-1}, \eqref{1.24}-\eqref{1.28}, \eqref{7.17-1} and noting the definition of $\phi$ in \eqref{1.22-1}, we get
\begin{equation}\label{7.18-1}
\begin{aligned}
\Big(\int_{\Omega}\int_{\R^3}\Big|\frac{1}{\sqrt{\mathcal{M}_{\v}}}R^{\v}\Big|^2dvdx\Big)^{\frac{1}{2}}&\leq C(\tau)\v^{N-1},\\
\Big(\int_{\Omega}\int_{\R^3}\Big|\frac{1}{\sqrt{\mathcal{M}_{\v}}}\bar{R}^{\pm,\v}\Big|^2dvdx\Big)^{\frac{1}{2}}&\leq
C\v^{n-\frac{1}{2}}\Big(\int_{\frac{1}{4\v}}^{\frac{1}{2\v}}(1+y)^{-2\fb-18}dy\Big)^{\frac{1}{2}} +C(\tau)(\v^{N-\frac{1}{2}}+\v^{\fb+n-\frac{3}{2}})\\
&\leq  C(\tau)(\v^{N-\frac{1}{2}}+\v^{\fb+n-\frac{3}{2}}),\\
\Big(\int_{\Omega}\int_{\R^3}\Big|\frac{1}{\sqrt{\mathcal{M}_{\v}}}\hat{R}^{\pm,\v}\Big|^2dvdx\Big)^{\frac{1}{2}}&\leq
	C\v^{n-1+\frac{n}{2}}\Big(\int_{\frac{1}{4\v}}^{\frac{1}{2\v}}e^{-\frac{1}{2N}\eta}d\eta\Big)^{\frac{1}{2}}\\
	&\quad +C(\tau)(\v^{N-1+\frac{n}{2}}+\v^{(n-1)(\fb+1)+n-2+\frac{n}{2}})\\
&\leq  C(\tau)(\v^{N-1+\frac{n}{2}}+\v^{(n-1)(\fb+1)+n-2+\frac{n}{2}}),
\end{aligned}
\end{equation}
which implies that
\begin{equation}\label{7.19}
	\begin{aligned}
		I_{5}\leq C(\v^{N-2n-k_0+3}+\v^{\fb-n-k_0+2.5})\|f_{R}^{\v}\|_{L^2}.
	\end{aligned}
\end{equation}
Substituting \eqref{7.11}, \eqref{7.15}-\eqref{7.17}, \eqref{7.18}-\eqref{7.19} into \eqref{7.10}, and taking $\delta$ small enough, we conclude the proof of \eqref{7.9}. $\hfill\square$

\subsection{Weighted $L^{\infty}$-estimate for remainder}
Given $(t,x,v)$, let $[X(s),V(s)]$ be the backward bi-characteristics of the scaled Boltzmann equation, that is
\begin{equation}\nonumber
	\left\{\begin{aligned}
		&\frac{dX(s)}{ds}=\frac{1}{\v^{n-2}}V(s),\quad \frac{dV(s)}{ds}=0,\\
		&[X(t),V(t)]=[x,v].
	\end{aligned}
	\right.
\end{equation}
The solution is then given by
$$
[X(s),V(s)]=[X(s;t,x,v),V(s;t,x,v)]=[x-\frac{1}{\v^{n-2}}(t-s)v,v].
$$

For each $(x,v)$ with $x\in \mathbb{T}^2\times [0,1]$ and $v_{3}\neq 0$, we define its {\it{backward exit time}} $t_{\mathbf{b}}\geq 0$ to be the last moment at which  $[X(s;0,x,v)]$ remains in $\mathbb{T}^2\times [0,1]$:
$$
t_{\mathbf{b}}(x,v)=\sup\{\lambda\geq 0:x-\frac{\lambda}{\v^{n-2}}v\in \mathbb{T}^2\times [0,1]\},
$$
where we regard $(x_{1}-\frac{\lambda v_{1}}{\v^{n-2}}, x_{2}-\frac{\lambda v_{2}}{\v^{n-2}})\in \R^2$ belongs to $\mathbb{T}^2$ with redefining them in $[0,1]^2$. We therefore have $x-\frac{1}{\v^{n-2}}t_{\mathbf{b}}v\in \partial\Omega$. Also, we define
$$
x_{\mathbf{b}}(x,v)=x(t_{\mathbf{b}})=x-\frac{1}{\v^{n-2}}t_{\mathbf{b}}v\in \mathbb{T}^2\times \{0,1\}.
$$

Now let $x\in\mathbb{T}^2\times [0,1]$, $(x,v)\not\in \gamma_{0}\cup \gamma_{-}^{0}\cup \gamma_{-}^{1}$ and $(t_{0},x_{0},v_{0})=(t,x,v)$, the back-time cycle is defined as
\begin{equation}\label{7.20-1}
\left\{\begin{aligned}
	&X_{cl}(s;t,x,v)=\sum\limits_{k}\mathbf{1}_{[t_{k+1},t_{k})}(s)\{x_{k}-\frac{1}{\v^{n-2}}(t_{k}-s)v_{k}\},\\
	&V_{cl}(s;t,x,v)=\sum\limits_{k}\mathbf{1}_{[t_{k+1},t_{k})}(s)v_{k},
\end{aligned}
\right.
\end{equation}
where $(t_{k+1},x_{k+1},v_{k+1})$ is defined as
$$
(t_{k+1},x_{k+1},v_{k+1}):=(t_{k}-t_{\mathbf{b}}(x_{k},v_{k}),x_{\mathbf{b}}(x_{k},v_{k}),R_{x_{k+1}}v_{k}).
$$
Clearly, for $k\geq 1$, it holds
\begin{equation}\label{7.21}
	t_{k}-t_{k+1}=\frac{\v^{n-2}}{|v_{0,3}|},\quad x_{k,3}=\frac{1-(-1)^k}{2}x_{1,3}+\frac{1+(-1)^k}{2}x_{2,3},\quad  v_{k,\sp}=v_{0,\sp},\quad v_{k,3}=(-1)^{k}v_{0,3}
\end{equation}

As in \cite{Guo-Huang-Wang}, we denote
$$
L_{M}g=-\frac{1}{\sqrt{\mu_{M}}}\left\{Q(\mathcal{M}_{\v},\sqrt{\mu_{M}}g)+Q(\sqrt{\mu_{M}}g,\mathcal{M}_{\v})\right\}=\nu_{\v}g-Kg,
$$
where the frequency $\nu_{\v}$ is defined in \eqref{7.9-1} and $Kg=K_2g-K_1g$ with
$$
\begin{aligned}
&K_1g=\int_{\R^3}\int_{\mathbb{S}^2}B(\theta)|v-u|\sqrt{\mu_{M}(u)}\frac{\mathcal{M}_{\v}(v)}{\sqrt{\mu_{M}(v)}}g(u)dud\omega,\\
&\begin{aligned}
	K_2g=&\int_{\R^3}\int_{\mathbb{S}^2}B(\theta)|v-u|\sqrt{\mu_{M}(u')}\frac{\mathcal{M}_{\v}(v')}{\sqrt{\mu_{M}(v)}}g(v')dud\omega\\
	&+\int_{\R^3}\int_{\mathbb{S}^2}B(\theta)|v-u|\sqrt{\mu_{M}(v')}\frac{\mathcal{M}_{\v}(u')}{\sqrt{\mu_{M}(v)}}g(u')dud\omega.
\end{aligned}
\end{aligned}
$$
From \cite[Lemma 3]{Guo2010}, it holds that $Kg(v)=\int_{\R^3}k(v,v')g(v')dv'$ where the kernel $k(v,v')$ satisfies
	\begin{equation}\label{7.22}
		|k(v,v')|\leq C\frac{\exp\{-c|v-v'|^2\}}{|v-v'|}.
	\end{equation}
Denoting $K_{\omega}g:=\omega_{\beta}K(\frac{g}{\omega_{\beta}})$, it follows from \eqref{1.23} and \eqref{1.36} that
\begin{equation}\label{7.23}
	\begin{aligned}
		&\partial_{t}h_{R}^{\v}+\frac{1}{\v^{n-2}}\cdot \nabla_{x}h_{R}^{\v}+\frac{1}{\v^{2n-2}}\nu_{\v}h_{R}^{\v}-\frac{1}{\v^{2n-2}}K_{\omega}h_{R}^{\v}\\
		&= -\frac{1}{\v^2}\frac{\omega_{\beta}(v)}{\sqrt{\mu_{M}}}\Big[Q\Big(r_{0},\frac{\sqrt{\mu_{M}}h_{R}^{\v}}{\omega_{\beta}(v)}\Big)+Q\Big(\frac{\sqrt{\mu_{M}}h_{R}^{\v}}{\omega_{\beta}(v)},r_{0}\Big)\Big]\\
		&\quad +\frac{1}{\v^{n-k_{0}}}\frac{\omega_{\beta}(v)}{\sqrt{\mu_{M}}}Q\Big(\frac{\sqrt{\mu_{M}}h_{R}^{\v}}{\omega_{\beta}(v)},\frac{\sqrt{\mu_{M}}h_{R}^{\v}}{\omega_{\beta}(v)}\Big)\\
		&\quad+\sum\limits_{i=1}^{N}\frac{1}{\v^{n-i}}\frac{\omega_{\beta}(v)}{\sqrt{\mu_{M}}}\Big[Q\Big(\frac{\sqrt{\mu_{M}}h_{R}^{\v}}{\omega_{\beta}(v)},F_{i}+\sum\limits_{\pm}[\phi(\v y)\bar{F}_{i}^{\pm}\phi(\v^{n}\eta)\hat{F}_{i}^{\pm}]\Big)\\
		&\qquad\qquad+Q\Big(F_{i}+\sum\limits_{\pm}[\phi(\v y)\bar{F}_{i}^{\pm}+\phi(\v^{n}\eta)\hat{F}_{i}^{\pm}],\frac{\sqrt{\mu_{M}}h_{R}^{\v}}{\omega_{\beta}(v)}\Big)\Big]\\
		&\quad +\frac{1}{\v^{2n-4+k_0}}\frac{\omega_{\beta}(v)}{\sqrt{\mu_{M}}}(R^{\v}+\sum\limits_{\pm}\bar{R}^{\pm,\v}+\sum\limits_{\pm}\hat{R}^{\pm,\v}):=\sum\limits_{j=1}^{4}S_{R,j}
	\end{aligned}
\end{equation}
\begin{lemma}\label{lem7.3}
Let $k_0>2n-1$. There exists $\v_{0}>0$ small such that for $t\in [0,\tau]$ and $\v\in (0,\v_0)$, it holds that
	\begin{equation}\label{7.23-1}
		\begin{aligned}
		\sup_{0\leq s\leq t}\|\v^{3n-3}h_{R}^{\v}(s)\|_{L^{\infty}}\leq & C(\tau)\big(\|\v^{3n-3}h_{R}^{\v}(0)\|_{L^{\infty}}+\v^{N+3n-k_0-2}+\v^{\fb+4n-k_0-3}\big)\\
	    &+C\sup_{0\leq s\leq t}\|f_{R}^{\v}(s)\|_{L^2}.
	    \end{aligned}
	\end{equation}
\end{lemma}

\noindent\textbf{Proof.} We divide the proof  into five steps.

{\it{Step 1.}} For any $t>0, (x,v)\not\in \g_{0}\cup \g_{-}^{0}\cup \g_{-}^{1}$, integrating \eqref{7.23} along the backward trajectory, and using \eqref{7.20-1}, one has that
\begin{equation}\label{7.24}
	\begin{aligned}
		h_{R}^{\v}(t,x,v)&=\exp\big\{-\frac{1}{\v^{2n-2}}\int_{0}^{t}\nu_{\v}(\xi,X_{cl}(\xi),V_{cl}(\xi))d\xi\big\}h_{R}^{\v}\big(0,x_{k}-\frac{1}{\v^{n-2}}v_{k}t_{k},v_{k}\big)\\
		&+\sum\limits_{l=0}^{k}\int_{\max\{t_{l+1},0\}}^{t_{l}}\exp\big\{-\frac{1}{\v^{2n-2}}\int_{s}^{t}\nu_{\v}(\xi,X_{cl}(\xi),V_{cl}(\xi))d\xi\big\}\\
		&\qquad \times \sum\limits_{j=1}^{4}S_{R,j}\big(s,x_{l}-\frac{v_{l}}{\v^{n-2}}(t_{l}-s),v_{l}\big)ds\\
		&+\frac{1}{\v^{2n-2}}\sum\limits_{l=0}^{k}\int_{\max\{t_{l+1},0\}}^{t_{l}}\exp\big\{-\frac{1}{\v^{2n-2}}\int_{s}^{t}\nu_{\v}\big(\xi,X_{cl}(\xi),V_{cl}(\xi)\big)d\xi\big\}\\
		&\qquad \times K_{\omega}h_{R}^{\v}\big(s,x_{l}-\frac{v_{l}}{\v^{n-2}}(t_{l}-s),v_{l}\big)ds\\
		&:=\sum\limits_{k=1}^{3}H_{k},
	\end{aligned}
\end{equation}
where $k$ is the collision number such that $0\in [t_{k+1},t_{k})$.

{\it{Step 2.}} Since $|V_{cl}(s)|=|v_0|=|v|$, $\nu_{\v}(\xi, X_{cl}(\xi), V_{cl}(\xi))\cong (1+|v|)\cong \nu(v)$,
then it holds that
\begin{equation}\label{7.26-1}
	\begin{aligned}
		&\int_{0}^{t}\exp\big\{-\frac{1}{\v^{2n-2}}\int_{s}^{t}\nu_{\v}(\xi, X_{cl}(\xi), V_{cl}(\xi))d\xi\big\}\nu_{\v}(s, X_{cl}(s),V_{cl}(s))ds\\
		&\leq  C\int_{0}^{t}\exp\big\{-\frac{\nu(v)(t-s)}{C\v^{2n-2}}\big\}\nu(v)ds\leq C\v^{2n-2},
	\end{aligned}
\end{equation}
and
\begin{equation}\label{7.27}
	|H_1|\leq C\exp\big\{-\frac{\nu(v)t}{C\v^{2n-2}}\big\}\|h_{R}^{\v}(0)\|_{L^{\infty}}.
\end{equation}

{\it{Step 3.}} By \cite[Lemma 5]{Guo2010}, we have
\begin{equation}\nonumber
	\vert\frac{\omega_{\beta}}{\sqrt{\mu_{M}}}Q(\sqrt{\mu_{M}}g_1,\sqrt{\mu_{M}}g_2)\vert\leq C\nu(v)\|\omega_{\beta}g_1\|_{L^{\infty}}\|\omega_{\beta}g_2\|_{L^{\infty}},
\end{equation}
which yields that
\begin{equation}\label{7.29}
	\begin{aligned}
		&\sum\limits_{j=1}^3|S_{R,j}|\leq C(\tau)\nu^{1-n}\nu(v)\|h_{R}^{\v}\|_{L^{\infty}}+C\v^{k_0-n}\nu(v)\|h_{R}^{\v}\|_{L^{\infty}}^2. 
	\end{aligned}
\end{equation}
Hence it follows from \eqref{7.26-1} and \eqref{7.29} that
\begin{align}\label{7.30}
	&\Big\vert\sum\limits_{l=0}^{k}\int_{\max\{t_{l+1},0\}}^{t_{l}}\exp\big\{-\frac{1}{\v^{2n-2}}\int_{s}^{t}\nu_{\v}(\xi,X_{cl}(\xi),V_{cl}(\xi))d\xi\big\}\sum\limits_{j=1}^3S_{R,j}ds\Big|\nonumber\\
	&\leq C\int_{0}^{t}\exp\big\{-\frac{\nu(v)(t-s)}{C\v^{2n-2}}\big\}\sum\limits_{j=1}^3|S_{R,j}|ds\nonumber\\
		&\leq C(\tau)\v^{n-1}\sup_{0\leq s\leq t}\|h_{R}^{\v}(s)\|_{L^{\infty}}+C\v^{k_0+n-2}\sup_{0\leq s\leq t}\|h_{R}^{\v}(s)\|_{L^{\infty}}^2.
\end{align}
For $S_{R,4}$, using \eqref{7.17-1} and a similar argument as in \eqref{7.18-1}, one has
$$
|S_{R,4}|\leq C(\tau)(\v^{N+3-k_0-2n}+\v^{\fb+2-n-k_0}),
$$
which yields that
\begin{equation}\label{7.31}
	\begin{aligned}
		&\Big\vert\sum\limits_{l=0}^{k}\int_{\max\{t_{l+1},0\}}^{t_{l}}\exp\big\{-\frac{1}{\v^{2n-2}}\int_{s}^{t}\nu_{\v}(\xi,X_{cl}(\xi),V_{cl}(\xi))d\xi\big\}S_{R,4}ds\Big|\\
		&\leq C(\tau) \int_{0}^{t}\exp\big\{-\frac{\nu(v)(t-s)}{C\v^{2n-2}}\big\}|S_{R,4}|ds\leq C(\tau)(\v^{N+1-k_0}+\v^{\fb+n-k_0}).
	\end{aligned}
\end{equation}
Combining \eqref{7.30}-\eqref{7.31} with \eqref{7.24}, we obtain
\begin{equation}\label{7.32}
	|H_2|\leq C(\tau)\big[\v^{n-1}\sup_{0\leq s\leq t}\|h_{R}^{\v}(s)\|_{L^{\infty}}+\v^{k_0+n-2}\sup_{0\leq s\leq t}\|h_{R}^{\v}(s)\|_{L^{\infty}}^2+\v^{N+1-k_0}+\v^{\fb+n-k_0}\big].
\end{equation}

{\it{Step 4.}}  Let $k_{\omega}(v,v')$ be the corresponding kernel associated with $K_{\omega}$. Using \eqref{7.22}, we have
\begin{equation}\label{7.33}
	|k_{\omega}(v,v')|\leq C\frac{\omega_{\beta}(v)\exp\{-c|v-v'|^2\}}{\omega_{\beta}(v')|v-v'|}\leq C\frac{\exp\{-\frac{3}{4}c|v-v'|^2\}}{|v-v'|}.
\end{equation}
Then $H_3$ can be bounded as
\begin{equation}\label{7.34}
	\begin{aligned}
		|H_3|&\leq \frac{1}{\v^{2n-2}}\sum\limits_{l=0}^{k}\int_{\max\{t_{l+1},0\}}^{t_{l}}\exp\big\{-\frac{1}{\v^{2n-2}}\int_{s}^{t}\nu_{\v}(\xi,X_{cl}(\xi),V_{cl}(\xi))d\xi\big\}\\
		&\quad \times \int_{\R^3}\vert k_{\omega}(V_{cl}(s),v')h_{R}^{\v}(s,x_{l}-\frac{v_{l}}{\v^{n-2}}(t_{l}-s),v')\vert dv'ds
	\end{aligned}
\end{equation}
We denote $V_{cl}'(s_1)=V_{cl}(s_{1};s,X_{cl}(s),v')$ and  $X_{cl}'(s_1)=X_{cl}(s_1;s,X_{cl}(s),v')$. Noting \eqref{7.27}, \eqref{7.32}, and using \eqref{7.24} and Vidav's's iteration in \eqref{7.34}, we obtain
\begin{align}
	&|H_3|\leq \frac{1}{\v^{4n-4}}\sum\limits_{l=0}^{k}\int_{\max\{t_{l+1},0\}}^{t_{l}}\exp\big\{-\frac{1}{\v^{2n-2}}\int_{s}^{t}\nu_{\v}(v)(\xi,X_{cl}(\xi),V_{cl}(\xi))d\xi\big\}\nonumber\\
	&\quad \times \Big(\int_{\R^3}\left\vert k_{\omega}(V_{cl}(s),v')\right\vert dv'ds\Big) \sum\limits_{j=0}^{k_{l}}\int_{\max\{t_{j+1}',0\}}^{t_{j}'}\exp\big\{-\frac{1}{\v^{2n-2}}\int_{s_{1}}^{s}\nu_{\v}(v)(\xi,X_{cl}'(\xi),V_{cl}'(\xi))d\xi'\big\}\nonumber\\
		&\quad \times \int_{\R^3}\left\vert k_{\omega}(V_{cl}'(s_1),v'')h_{R}^{\v}(s_1,X_{cl}'(s_1),v'')\right\vert dv''ds_1\nonumber\\
		&\quad +C(\tau)\big(\|h_{R}^{\v}(0)\|_{L^{\infty}}+\v^{N+1-k_0}+\v^{\fb+n-k_0}+\v^{n-1}\sup_{0\leq s\leq t}\|h_{R}^{\v}(s)\|_{L^{\infty}}+\v^{k_0+n-2}\sup_{0\leq s\leq t}\|h_{R}^{s}\|_{L^{\infty}}^2\big)\nonumber\\
		&:=H_{3,1}+C(\tau)\big(\|h_{R}^{\v}(0)\|_{L^{\infty}}+\v^{N+1-k_0}+\v^{\fb+n-k_0}\nonumber\\
		&\qquad\qquad\quad+\v^{n-1}\sup_{0\leq s\leq t}\|h_{R}^{\v}(s)\|_{L^{\infty}} +\v^{k_0+n-2}\sup_{0\leq s\leq t}\|h_{R}^{s}\|_{L^{\infty}}^2\big),\label{7.35}
\end{align}
where $k_{l}$ is the collision number such that $0\in [t_{k_{l}+1}',t_{k_{l}}')$, and we have used the fact
\begin{equation}\label{7.36}
	\int_{\R^3}\left\vert k_{\omega}(v,v')\right\vert dv'\leq C(1+|v|)^{-1}.
\end{equation}

We still need to control $H_{3,1}$ carefully. As in \cite{Guo-Huang-Wang}, we divide it into three cases.

{\it Case 1.} For $|v|\geq \mathfrak{N}$, using \eqref{7.36}, a direct calculation shows that
\begin{equation}\label{7.37}
	\begin{aligned}
		|H_{3,1}|&\leq \frac{C}{\v^{4n-4}}\sup_{0\leq s\leq t}\|h_{R}^{\v}(s)\|_{L^{\infty}} \sum\limits_{l=0}^{k}\int_{\max\{t_{l+1},0\}}^{t_{l}} \exp\big\{-\frac{\nu(v)(t-s)}{C\v^{2n-2}}\big\}ds\\
		&\qquad \times \int_{\R^3}|k_{\omega}(V_{cl}(s),v')|dv'\int_{0}^{s}\exp\big\{-\frac{\nu(v')(s-s_1)}{C\v^{2n-2}}\big\}ds_1\\
		&\leq \frac{C}{\mathfrak{N}\v^{2n-2}}\sup_{0\leq s\leq t}\|h_{R}^{\v}(s)\|_{L^{\infty}}\int_{0}^{t}\exp\big\{-\frac{\nu(v)(t-s)}{C\v^{2n-2}}\big\}ds\\
		&\leq \frac{C}{\mathfrak{N}}\sup_{0\leq s\leq t}\|h_{R}^{\v}(s)\|_{L^{\infty}}.
	\end{aligned}
\end{equation}

{\it Case 2.} For either $|v|\leq \mathfrak{N},|v'|\geq 2\mathfrak{N}$ or $|v'|\leq 2\mathfrak{N}, |v''|\geq 3\mathfrak{N}$. Noting $|V_{cl}(s)|=|v|$ and $|V_{cl}'(s_1)|=|v'|$, it holds that either $|V_{cl}(s)-v'|\geq \mathfrak{N}$ or $|V_{cl}'(s_1)-v''|\geq \mathfrak{N}$, then either one following holds for some small positive constant $0<c_{1}\leq \frac{c}{32}$ (where $c$ is the one in \eqref{7.22}):
$$
\begin{aligned}
	&|k_{\omega}(V_{cl}(s),v')|\leq e^{-c_1\mathfrak{N}^2}|k_{\omega}(V_{cl}(s),v')|\exp(c_1|V_{cl}(s_1)-v'|^2),\\
	&|k_{\omega}(V_{cl}'(s_{1}),v'')|\leq e^{-c_1\mathfrak{N}^2}|k_{\omega}(V_{cl}'(s_{1}),v'')|\exp(c_1|V_{cl}'(s_1)-v''|^2),
\end{aligned}
$$
which, together with \eqref{7.33}, yields
\begin{equation}\nonumber
	\int_{\R^3}|k_{\omega}(v,v')|e^{c_1|v-v'|^2}dv'\leq \frac{C}{1+|v|}\quad \text{and}\quad  \int_{\R^3}|k_{\omega}(v',v'')|e^{c_1|v'-v''|^2}dv''\leq \frac{C}{1+|v'|}.
\end{equation}
Hence, the corresponding part of $H_{3,1}$ is bounded by
\begin{equation}\label{7.39}
	\begin{aligned}
	& \frac{C}{\v^{4n-4}}e^{-c_1\mathfrak{N}^2}\sup_{0\leq s\leq t}\|h_{R}^{\v}(s)\|_{L^{\infty}}\sum\limits_{l=0}^{k}\int_{\max\{t_{l+1},0\}}^{t_{l}}\exp\big\{-\frac{\nu(v)(t-s)}{C\v^{2n-2}}\big\}ds\\
		&\quad \times \int_{\R^3}|k_{\omega}(V_{cl}(s),v')|e^{c_{1}(|V_{cl}(s)-v'|^2)}dv'\int_{0}^{s}\exp\big\{-\frac{\nu(v')(s-s_1)}{C\v^{2n-2}}\big\}ds_1\\
		&\leq Ce^{-c_1\mathfrak{N}^2}\sup_{0\leq s\leq t}\|h_{R}^{\v}(s)\|_{L^{\infty}}.
	\end{aligned}
\end{equation}

{\it Case 3.} For $|v|\leq \mathfrak{N}, |v'|\leq  2\mathfrak{N}, |v''|\leq 3\mathfrak{N}$, which is the last remaining case. It follows from \eqref{7.21} that
\begin{equation}\label{7.39-1}
|k_{l}|\leq \frac{s|v'|}{\v^{n-2}}\leq \frac{2\mathfrak{N}s}{\v^{n-2}}\leq \frac{2\mathfrak{N}t}{\v^{n-2}},\quad \text{for }s\in [\max\{t_{l+1},0\},t_{l}).
\end{equation}
We first consider $t_{j}'-s_1\leq \delta \v^{3n-4}$ for some small $\delta>0$ determined later.
Noting $\nu(v)\geq \nu_{0}>0$ where $\nu_0$ is a positive constant independent of $v$, and using \eqref{7.39-1}, then the corresponding part of $H_{3,1}$ is bounded by
\begin{equation}\label{7.40}
	\begin{aligned}
&\frac{C}{\v^{4n-4}}\sup_{0\leq s\leq t}\|h_{R}^{s}(s)\|_{L^{\infty}}\times \sum\limits_{l=0}^{k}\int_{\max\{t_{l+1},0\}}^{t_{l}}\exp\big\{-\frac{\nu_0(t-s)}{C\v^{2n-2}}\big\}\int_{|v'|\leq 2\mathfrak{N}}|k_{\omega}(V_{cl}(s),v')|dv'ds\\
		&\quad \times \sum\limits_{j=0}^{k_{l}}\int_{t_{j}'-\delta\v^{3n-4}}^{t_{j}'}\exp\big\{-\frac{\nu_0(s-s_1)}{C\v^{2n-2}}\big\}\int_{|v''|\leq 3\mathfrak{N}} |k_{\omega}(V_{cl}'(s_{1}),v'')|dv''ds_1\\
		&\leq \frac{C(\mathfrak{N})\tau\delta}{\v^{2n-2}}\sup_{0\leq s\leq t}\|h_{R}^{s}(s)\|_{L^{\infty}}\int_{0}^{t}\exp\big\{-\frac{\nu_0(t-s)}{C\v^{2n-2}}\big\}ds\\
		&\leq C(\mathfrak{N})\tau\delta\sup_{0\leq s\leq t}\|h_{R}^{\v}(s)\|_{L^{\infty}}.
	\end{aligned}
\end{equation}
Next, we consider the case $t_{j}'-s_1\geq \delta\v^{3n-4}$. We denote $D:=\{(v',v'')\in \R^6: |v'|\leq 2\mathfrak{N}, |v''|\leq 3\mathfrak{N}\}$  and note that
$$
\int_{R^3}|k_{\omega}(v,v')|^2dv'\leq C.
$$
Hence, using Cauchy inequality, the corresponding part of $H_{3,1}$ can be bounded by
\begin{equation}\label{7.41}
	\begin{aligned}
	& \frac{1}{\v^{4n-4}}\sum\limits_{l=0}^{k}\int_{\max\{t_{l+1},0\}}^{t_{l}}\exp\big\{-\frac{\nu_{0}(t-s)}{C\v^{2n-2}}\big\}ds\\
		&\quad
		\times \Big(\iint_{D}\sum_{j=0}^{k_{l}}\int_{\max\{t_{j+1}',0\}}^{t_{j}'-\delta\v^{2n-2}}\exp\big\{-\frac{\nu_{0}(s-s_1)}{C\v^{2n-2}}\big\}|h_{R}^{\v}(s_1,X_{cl}'(s_1),v''|^2dv''dv'ds_1\Big)^{\frac{1}{2}}\\
		&\quad  \times\Big(\iint_{D}\sum_{j=0}^{k_{l}}\int_{\max\{t_{j+1}',0\}}^{t_{j}'-\delta\v^{3n-4}}\exp\big\{-\frac{\nu_{0}(s-s_1)}{C\v^{2n-2}}\big\}|k_{\omega}(V_{cl}(s),v')k_{\omega}(V_{cl}'(s_1),v'')|^2dv''dv'ds_1\Big)^{\frac{1}{2}}\\
		&\leq \frac{C({\mathfrak{N}})}{\v^{3n-3}}\sum\limits_{l=0}^{k}\int_{\max\{t_{l+1},0\}}^{t_{l}}\exp\big\{-\frac{\nu_{0}(t-s)}{C\v^{2n-2}}\big\}ds\\
		& \quad \times \Big(\iint_{D}\sum_{j=0}^{k_{l}}\int_{\max\{t_{j+1}',0\}}^{t_{j}'-\delta\v^{3n-4}}\exp\big\{-\frac{\nu_{0}(s-s_1)}{C\v^{2n-2}}\big\}|f_{R}^{\v}(s_1,X_{cl}'(s_1),v'')|^2dv''dv'ds_1\Big)^{\frac{1}{2}}.
	\end{aligned}
\end{equation}
We denote $v_{0}':=v'$ and $y_{j}:=X_{cl}'(s_1)$. Then for $s_{1}\in (\max\{t_{j+1}',0\},t_{j}')$, we have $y_{j}=x_{j}'-\frac{1}{\v^{n-2}}v_{j}'(t_{j}'-s_1)$ with $t_{j}'=t_{j}'(s_1;s, X_{cl}(s),v_{0}')$. We note that
\begin{equation}\nonumber
	s-t_{j}'=\left\{
	\begin{aligned}
		&\Big(\frac{\v^{n-2}X_{cl,3}(s)}{|v_{0,3}'|}+(j-1)\frac{\v^{n-2}}{|v_{0,3}'|}\Big)\mathbf{1}_{j\geq 1 }\qquad \quad\,\,\,\text{for }v_{0,3}'>0,\\
		&\Big(\frac{\v^{n-2}(1-X_{cl,3}(s))}{|v_{0,3}'|}+(j-1)\frac{\v^{n-2}}{|v_{0,3}'|}\Big)\mathbf{1}_{j\geq 1 }\quad \text{for }v_{0,3}'<0,
	\end{aligned}
	\right.
\end{equation}
then a direct calculation shows that
\begin{equation}\label{7.43}
	\left\{\begin{aligned}
		&y_{j,1}=x_{0,1}'-\frac{1}{\v^{n-2}}v_{0,1}'(s-s_1),\\
		&y_{j,1}=x_{0,2}'-\frac{1}{\v^{n-2}}v_{0,2}'(s-s_1),\\
		&y_{j,3}=x_{j,3}'-\frac{1}{\v^{n-2}}(-1)^{j}[v_{0,3}'(s-s_1)-\v^{n-2}\mathbf{1}_{j\geq 1}(X_{cl,3}(s)+(j-1))],
	\end{aligned}
	\right.\quad \text{for } v_{0,3}'>0,
\end{equation}
and
\begin{equation}\label{7.44}
	\left\{\begin{aligned}
		&y_{j,1}=x_{0,1}'-\frac{1}{\v^{n-2}}v_{0,1}'(s-s_1),\\
		&y_{j,1}=x_{0,2}'-\frac{1}{\v^{n-2}}v_{0,2}'(s-s_1),\\
		&y_{j,3}=x_{j,3}'-\frac{1}{\v^{n-2}}(-1)^{j}[v_{0,3}'(s-s_1)+\v^{n-2}\mathbf{1}_{j\geq 1}(j-X_{cl,3}(s))],
	\end{aligned}
	\right. \qquad \quad\,\,\,\,\text{for } v_{0,3}'<0,
\end{equation}
where $x_{0}'=X_{cl}(t_{0}';s,X_{cl}(s),v_{0}')=X_{cl}(s)$. We note that $x_{j,3}'\in \{0,1\}$ for $1\leq j\leq k_{l}$. For $s_1\in [\max\{t_{j+1}',0\},t_{j}'-\delta\v^{3n-4}]$, we have
$$
s-s_1\geq s-t_{j}'+\delta\v^{3n-4}\geq \delta\v^{3n-4},\quad \text{for }j=0,1,\cdots,k_{l},
$$
which, together with \eqref{7.43}-\eqref{7.44}, yields that
\begin{equation}\label{7.45}
	\Big\vert\mathrm{det}\Big(\frac{\partial y_{j}}{\partial v'}\Big)\Big\vert=\frac{1}{\v^{3n-6}}(s-s_1)^3\geq \delta^3\v^{6n-6}.
\end{equation}
Using \eqref{7.45} and making a change of variable $v'\mapsto y_{j}=X_{cl}'(s_1)$ implies that
\begin{equation}\nonumber
	\int_{|v'|\leq 2\mathfrak{N}}|f_{R}^{\v}(s_1, X_{cl}(s_1),v'')|^2dv'\leq \frac{C}{\delta^3\v^{6n-6}}\int_{\Omega}|f_{R}^{\v}(s_1,z,v'')|^2dz
\end{equation}
for $s_1\in [\max\{t_{j+1}',0\},t_{j}'-\delta\v^{3n-4}]$. Hence \eqref{7.41} can be bounded by
\begin{equation}\label{7.47}
 \frac{C({\mathfrak{N},\delta})}{\v^{5n-5}}\int_{0}^{t}\exp\big\{-\frac{\nu_{0}(t-s)}{C\v^{2n-2}}\big\}ds\sup_{0\leq s\leq t}\|f_{R}^{\v}(s)\|_{L^2}\leq \frac{C(\mathfrak{N},\delta)}{\v^{3n-3}}\sup_{0\leq s\leq t}\|f_{R}^{\v}(s)\|_{L^2}.
\end{equation}

{\it{Step 5.}} Combining \eqref{7.27}, \eqref{7.32}, \eqref{7.35}, \eqref{7.37}, \eqref{7.39}, \eqref{7.40} and \eqref{7.47}, we get
\begin{equation}\label{7.48}
	\begin{aligned}
		\sup_{0\leq s\leq t}\|\v^{3n-3}h_{R}^{\v}(s)\|_{L^{\infty}}&\leq C(\tau)(\|\v^{3n-3}h_{R}(0)\|_{L^{\infty}}+\v^{N+3n-k_0-2}+\v^{\fb+4n-k_0-3})\\
		&+C(\tau)\big(\frac{1}{\mathfrak{N}}+C(\mathfrak{N})\delta+\v^{n-1}\big)\sup_{0\leq s\leq t}\|\v^{3n-3}h_{R}^{s}\|_{L^{\infty}}\\
		&+C(\tau)\v^{k_0-2n+1}\sup_{0\leq s\leq t}\|\v^{3n-3}h_{R}^{s}\|_{L^{\infty}}^2+C({\mathfrak{N},\delta})\sup_{0\leq s\leq t}\|f_{R}^{\v}(s)\|_{L^2}.
	\end{aligned}
\end{equation}
First taking $\mathfrak{N}\gg 1$ large enough and then $\delta>0$ small enough, and finally choosing $\v_{0}$ small enough in \eqref{7.48}, we deduce \eqref{7.23-1}. This completes the proof of Lemma \ref{lem7.3}.$\hfill\square$

\subsection{Proof of Theorem \ref{theorem1.1}}
In this subsection, we complete the proof of Theorem \ref{theorem1.1}. With the help of Lemmas \ref{lem7.1} and \ref{lem7.3}, the proof is the same as \cite{Guo Jang Jiang-1}. Indeed, combining \eqref{7.23-1} with \eqref{7.9}, we obtain
\begin{equation}\nonumber
	\begin{aligned}
		&\frac{d}{dt}\|f_{R}^{\v}(t)\|_{L^2}^2+\frac{c_{0}}{2\v^{2n-2}}\|(\mathbf{I}-\mathbf{P}_{\v})f_{R}^{\v}(t)\|_{\nu}^2\\
		&\leq C(\tau)(1+\v^{2k_0-6n+4}\|\v^{3n-3}h_{R}^{\v}(0)\|_{L^{\infty}}^2)(1+\|f_{R}^{\v}\|_{L^2}^2)\\
		&\leq C(\tau)(1+\v^{2k_0-6n+4}\sup_{0\leq s\leq \tau}\|f_{R}^{\v}(s)\|_{L^2}^2+\v^{2k_0-6n+4}\|\v^{3n-3}h_{R}^{\v}(0)\|_{L^{\infty}}^2)(1+\|f_{R}(t)\|_{L^2}^2),
	\end{aligned}
\end{equation}
which, together with the Gr\"{o}nwall inequality, implies that
\begin{equation}\label{7.50}
	\begin{aligned}
&\sup_{0\leq s\leq \tau}\|f_{R}^{\v}(s)\|_{L^2}^2\leq (1+\|f_{R}^{\v}(0)\|_{L^2}^2)\\
&\qquad \times \exp\big\{C(\tau)(1+\v^{2k_0-6n+4}\|\v^{3n-3}h_{R}(0)\|_{L^{\infty}}^2+\v^{2k_0-6n+4}\sup_{0\leq s\leq \tau}\|f_{R}(s)\|_{L^2}^2)\big\}.
\end{aligned}
\end{equation}
Noting $k_{0}>3n-2$ and taking $\v$ small enough in \eqref{7.50}, one deduces that
\begin{equation}\nonumber
\sup_{0\leq s\leq \tau}\|f_{R}^{\v}(s)\|_{L^2}\leq C(\tau)(1+\|f_{R}^{\v}(0)\|_{L^2}+\|\v^{3n-3}h_{R}^{\v}(0)\|_{L^{\infty}}),
\end{equation}
which, together with \eqref{1.37-1} and \eqref{7.23-1}, yields \eqref{1.22-3}. Therefore the proof of Theorem \ref{theorem1.1} is complete. $\hfill\square$

\vspace{1.5mm}

\appendix
\renewcommand{\appendixname}{Appendix~\Alph{section}}

\section{Proof of Lemma \ref{lem2.1}}\label{AppendixA}
It is clear that \eqref{2.0} follows from \eqref{2.1}. We only need to prove \eqref{2.17}, and we divide the proof into two steps.

{\it Step 1.} For $1\leq k\leq n-3$, multiplying $\eqref{1.7-1}_5$ by $1, v$ and integrating the resultant equation over $\R^3$ with respect to $v$, one deduces that
\begin{equation}\label{2.5}
	\operatorname{div}_{x}u_{k}=0,\qquad \nabla_{x}(\r_k+\t_k)=0.
\end{equation}
Thus, we may assume
\begin{equation}\nonumber
	p_{k}:=\r_{k}+\t_{k}\equiv p_{k}(t)\qquad \text{for }1\leq k\leq n-3,
\end{equation}
where $p_{k}(t)$ is a function independent of $x$, and will be determined in Lemma\ref{lem2.2} below to guarantee the compatibility conditions for the equations of velocity.

{\it Step 2.} We consider the case of $k \geq n-2$. Multiplying $\eqref{1.7-1}_{5}$ by the five collision invariant and integrating the resultant equation over $\R^3$ with respect to $v$, we have
\begin{equation}\nonumber
	\left\{\begin{aligned}
		&\partial_{t} \int_{\mathbb{R}^{3}} F_{k+2-n} d v+\operatorname{div}_{x} \int_{\mathbb{R}^{3}} v F_{k} d v=0, \\
		&\partial_{t} \int_{\mathbb{R}^{3}} v_{i} F_{k+2-n} d v+\operatorname{div}_{x} \int_{\mathbb{R}^{3}} v_{i} vF_{k} d v=0\quad (i=1,2,3),\\
		&\partial_{t} \int_{\mathbb{R}^{3}}|v|^{2} F_{k+2-n} d v+\operatorname{div}_{x} \int_{\mathbb{R}^{3}}|v|^{2} vF_{k} d v=0,
	\end{aligned}\right.
\end{equation}
which, together with \eqref{2.2}, yields that
\begin{equation}\label{2.7}
	\left\{\begin{aligned}
		&\partial_{t} \rho_{k+2-n}+\operatorname{div}_{x} u_{k}=0, \\
		&\partial_{t} u_{k+2-n, i}+\partial_{x_i}\left(\rho_{k}+\theta_{k}\right)+\sum_{j=1}^{3} \partial_{x_j}\left\langle\mathcal{A}_{i j},(\mathbf{I}-\mathbf{P}) f_{k}\right\rangle=0,\quad i=1,2,3,\\
		&3 \partial_{t}\left(\rho_{k+2-n}+\theta_{k+2-n}\right)+5 \operatorname{div}_{x} u_{k}+2 \sum\limits_{j=1}^3\pa_{x_j}\left\langle\mathcal{B}_{j},(\mathbf{I}-\mathbf{P}) f_{k}\right\rangle=0.
	\end{aligned}\right.
\end{equation}

{\it Step 2.1.} If $k=n-2$, it follows from $\eqref{1.13}_2$ and \cite[p. 649]{Guo2006} that
\begin{equation}\label{2.6-1}
\begin{aligned}
	(\mathbf{I-P})f_{n-2}&=\frac{1}{2}(\mathbf{I-P})(\frac{\mathbf{P}f_{0}\mathbf{P}f_0}{\sqrt{\mu}})\\
	&=\frac{1}{2}\sum\limits_{l,s=1}^3u_{0,l}u_{0,s}\mathcal{A}_{ls}+\sum\limits_{l=1}^{3}\t_{0}u_{0,l}\mathcal{B}_{l}+\frac{\t_{0}^2}{8}(\mathbf{I-P})\big\{(|v|^2-5)^2\sqrt{\mu}\big\},
\end{aligned}
\end{equation}
which implies that
\begin{equation}\label{2.7-0}
	\begin{aligned}
		\la \mathcal{A}_{ii}, (\mathbf{I-P})f_{n-2}\ra&=\frac{2}{3}|u_{0,i}|^2-\frac{1}{3}\sum\limits_{j\neq i}|u_{0,j}|^2,\quad i=1,2,3,\\
		\la \mathcal{A}_{ij}, (\mathbf{I-P})f_{n-2}\ra&=u_{0,i}u_{0,j},\quad j\neq i,\\
		\la \mathcal{B}_{i}, (\mathbf{I-P})f_{n-2}\ra&=\frac{5}{2}\t_{0}u_{0,i},\quad i=1,2,3,
	\end{aligned}
\end{equation}
where we have used the fact that $\la \mathcal{A}_{ij},\mathcal{B}_{k}\ra=0$ and
\begin{equation}\label{2.7-6}
\begin{aligned}
	\la \mathcal{A}_{ij},\mathcal{A}_{ij}\ra=1,\quad\la \mathcal{A}_{ii},\mathcal{A}_{ii}\ra=\frac{4}{3},\quad \la \mathcal{B}_{i},\mathcal{B}_{i}\ra=\frac{5}{2},\quad \la \mathcal{A}_{ii},\mathcal{A}_{jj}\ra =-\frac{2}{3},\,\,i\neq j.
\end{aligned}
\end{equation}
Substituting \eqref{2.7-0} into \eqref{2.7}, and using \eqref{2.4}-\eqref{2.5-1}, we get  \eqref{2.8-1}.

{\it Step 2.2.} If $k\geq n-1$, we denote $k=n-2+m$ with $m\geq 1$. Using $\eqref{2.0}_2$, we have
\begin{equation}\label{2.7-1}
	\begin{aligned}
		\la \mathcal{A}_{ij}, (\mathbf{I-P})f_{m+n-2}\ra&=\la \mathcal{A}_{ij},\mathbf{L}^{-1}\{\Gamma(f_{0},\mathbf{P}f_{m})+\Gamma(\mathbf{P}f_{m},f_{0})\}\ra+\la \mathcal{A}_{ij},G_{m-1}\ra,\\
		\la \mathcal{B}_{i}, (\mathbf{I-P})f_{m+n-2}\ra&=\la \mathcal{B}_{i},\mathbf{L}^{-1}\{\Gamma(f_{0},\mathbf{P}f_{m})+\Gamma(\mathbf{P}f_{m},f_{0})\}\ra+\la \mathcal{B}_{i},G_{m-1}\ra
	\end{aligned}
\end{equation}
where $G_{m-1}$ is the one defined in \eqref{2.13}, which depends only on $f_{i}$ with $1\leq i\leq m-1$.

From \cite[p. 649]{Guo2006}, one obtains
\begin{align}\label{2.7-3}
		&\mathbf{L}^{-1}\big\{\Gamma(f_{0},\mathbf{P}f_{m})+\Gamma(\mathbf{P}f_{m},f_{0})\big\}=(\mathbf{I-P})(\frac{f_{0}\mathbf{P}f_{m}}{\sqrt{\mu}})\nonumber\\
		&=\sum\limits_{l,s=1}^3u_{0,l}u_{m,s}\mathcal{A}_{ls}+\sum\limits_{l=1}^3\theta_{0}u_{m,l}\mathcal{B}_l+\sum\limits_{l=1}^3\theta_{m}u_{0,l}\mathcal{B}_{l}+\frac{1}{4}\theta_{0}\theta_{m}(\mathbf{I-P})\big\{(|v|^2-5)^2\sqrt{\mu}\big\}.
\end{align}
Thus, using \eqref{2.7-6}, we get
\begin{equation}\label{2.7-4}
	\begin{aligned}
	\la \mathcal{A}_{ii},\mathbf{L}^{-1}\{[\Gamma(f_{0},\mathbf{P}f_{m})+\Gamma(\mathbf{P}f_{m},f_{0})]\}\ra&=\frac{4}{3}u_{0,i}u_{m,i}-\frac{2}{3}\sum\limits_{j\neq i}u_{0,j}u_{m,j},\quad i=1,2,3,\\
		\la \mathcal{A}_{ij},\mathbf{L}^{-1}\{[\Gamma(f_{0},\mathbf{P}f_{m})+\Gamma(\mathbf{P}f_{m},f_{0})]\}\ra&=u_{0,i}u_{m,j}+u_{m,i}u_{0,j},\quad i\neq j,\\
		\la\mathcal{B}_{i},\mathbf{L}^{-1}\{[\Gamma(f_{0},\mathbf{P}f_{m})+\Gamma(\mathbf{P}f_{m},f_{0})]\}\ra&=\frac{5}{2}\t_{0}u_{m,i}+\frac{5}{2}\t_{m}u_{0,i},\quad i=1,2,3.\\
	\end{aligned}
\end{equation}
Substituting \eqref{2.7-1}, \eqref{2.7-4} into \eqref{2.7} and using $\operatorname{div}_{x}u_0=0$, one deduces \eqref{2.17}. Therefore the proof of Lemma \ref{lem2.1} is complete. $\hfill\square$

\section{Proof of Proposition \ref{prop3.1}}\label{AppendixD}
For simplicity, we only prove the formulation of the lower viscous boundary expansion $\bar{F}_{i}^{-}$, since the upper viscous boundary expansions $\bar{F}_{i}^{+}$ can be formulated similarly. For this reason, we drop the superscript ``-'' of the lower boundary layers for simplicity of notations. We divide the proof into two steps.

{\it Step 1.} For $1 \leq k \leq n-2$, noting $\bar{F}_{k+2-n}=0$, we have from $\eqref{1.14}$ that
\begin{equation}\label{3.4}
	\left\{\begin{aligned}
		&0=\int_{\mathbb{R}^{3}} v_{\sp} \cdot \nabla_{\sp} \bar{F}_{k} d v+\partial_{y} \int_{\mathbb{R}^{3}} v_{3} \bar{F}_{k+1} d v, \\
		&0=\int_{\mathbb{R}^{3}} v_{i} v_{\sp} \cdot \nabla_{\sp} \bar{F}_{k} d v+\partial_{y} \int_{\mathbb{R}^{3}} v_{3} v_{i} \bar{F}_{k+1} d v,\quad i=1,2,3,\\
		&0=\int_{\mathbb{R}^{3}}|v|^{2} v_{\sp} \cdot \nabla_{\sp} \bar{F}_{k} d v+\partial_{y} \int_{\mathbb{R}^{3}} v_{3}|v|^{2} \bar{F}_{k+1} d v.
	\end{aligned}\right.
\end{equation}
Due to the macro-micro decomposition, a direct calculation shows that
\begin{equation}\nonumber
	\left\{
	\begin{aligned}
		&\int_{\mathbb{R}^{3}} \bar{F}_{k} d v=\bar{\rho}_{k}, \quad \int_{\mathbb{R}^{3}}v_{i}\bar{F}_{k} d v=\bar{u}_{k,i}, \quad \int_{\mathbb{R}^{3}} v_{i} v_{j} \bar{F}_{k} d v=\delta_{i j}\left(\bar{\rho}_{k}+\bar{\theta}_{k}\right)+\left\langle\mathcal{A}_{i j},(\mathbf{I}-\mathbf{P}) \bar{f}_{k}\right\rangle,\\
		&\int_{\R^3}|v|^2\bar{F}_{k}dv=3(\bar{\rho}_{k}+\bar{\theta}_{k})+\sum\limits_{i=1}^{3}\langle \mathcal{A}_{ii}, (\mathbf{I}-\mathbf{P})\bar{f}_{k}\rangle,\\
		&\int_{\mathbb{R}^{3}} v_{i}|v|^{2} \bar{F}_{k} d v=5 \bar{u}_{k, i}+2\left\langle\mathcal{B}_{i},(\mathbf{I}-\mathbf{P})\bar{f}_{k}\right\rangle.
	\end{aligned}
\right.
\end{equation}
then \eqref{3.4} can be written as
\begin{equation}\label{3.6}
	\left\{\begin{aligned}
		&\operatorname{div}_{\sp} \bar{u}_{k, \sp}+\partial_{y} \bar{u}_{k+1,3}=0, \\
		&\partial_{x_{i}}\left(\bar{\rho}_{k}+\bar{\theta}_{k}\right)+\sum_{j=1}^{2} \partial_{x_{j}}\left\langle\mathcal{A}_{i j},(\mathbf{I}-\mathbf{P}) \bar{f}_{k}\right\rangle+\partial_{y}\left\langle\mathcal{A}_{3 i},(\mathbf{I}-\mathbf{P}) \bar{f}_{k+1}\right\rangle=0, \quad i=1,2, \\
		&\partial_{y}\left(\bar{\rho}_{k+1}+\bar{\theta}_{k+1}\right)+\sum_{j=1}^{2} \partial_{x_j}\left\langle\mathcal{A}_{3 j},(\mathbf{I}-\mathbf{P}) \bar{f}_{k}\right\rangle+ \partial_{y}\left\langle\mathcal{A}_{33},(\mathbf{I}-\mathbf{P}) \bar{f}_{k+1}\right\rangle=0, \\
		&5 \operatorname{div}_{\sp} \bar{u}_{k, \sp}+2\sum\limits_{j=1}^{2}\partial_{x_{j}}\la \mathcal{B}_{j},(\mathbf{I}-\mathbf{P}) \bar{f}_{k}\ra+5 \partial_{y} \bar{u}_{k+1,3}+2 \partial_{y}\left\langle\mathcal{B}_{3},(\mathbf{I}-\mathbf{P}) \bar{f}_{k+1}\right\rangle=0.
	\end{aligned}\right.
\end{equation}
If $1\leq k\leq n-3$, using \eqref{3.1}, one deduces that
\begin{equation}\label{3.7}
	\left\{\begin{aligned}
		&\operatorname{div}_{\sp} \bar{u}_{k, \sp}+\partial_{y} \bar{u}_{k+1,3}=0, \quad \partial_{x_i}\left(\bar{\rho}_{k}+\bar{\theta}_{k}\right)=0\quad i=1,2,\\
		&\partial_{y}\left(\bar{\rho}_{k+1}+\bar{\theta}_{k+1}\right)=0.
	\end{aligned}\right.
\end{equation}

Recalling the far-field condition $\eqref{3.3-3}$, we have from \eqref{3.7} that
\begin{equation}\label{3.7-1}
		\bar{u}_{k+1,3}=\int_{y}^{+\infty}\operatorname{div}_{\sp} \bar{u}_{k, \sp}(t,x_{\sp},s)ds,\quad
		\bar{\rho}_{k+1}+\bar{\theta}_{k+1}\equiv 0,
	\quad \text{for }1\leq k\leq n-3.
\end{equation}

Moreover, if $k=n-2$, it follows from $\eqref{3.1}_{2}$ and \cite[p. 649]{Guo2006} that
\begin{equation}\label{3.8}
	\begin{aligned}
		&(\mathbf{I}-\mathbf{P}) \bar{f}_{n-1} =\mathbf{L}^{-1}\{\Gamma\left(f_{0}^{-}, \bar{f}_{1}\right)+\Gamma(\bar{f}_{1}, f_{0}^{-})\}=(\mathbf{I}-\mathbf{P})(\frac{f_{0}^{-} \bar{f}_{1}}{\sqrt{\mu}}), \\
		&=\sum_{i, j=1}^{3} u_{0, i}^{-} \bar{u}_{1, j} \mathcal{A}_{i j}+\sum_{i=1}^{3} \bar{\theta}_{1} u_{0, i}^{-} \mathcal{B}_{i}+\sum_{i=1}^{3} \theta_{0}^{-} \bar{u}_{1, i} \mathcal{B}_{i} +\frac{1}{4}\t_{0}^{-}\bar{\t}_{1}(\mathbf{I-P})\{(|v|^2-5)^2\sqrt{\mu}\}.
	\end{aligned}
\end{equation}
Using \eqref{2.7-6}, we obtain that
\begin{align}\label{3.9}
		\left\langle\mathcal{A}_{l3}, (\mathbf{I}-\mathbf{P}) \bar{f}_{n-1}\right\rangle&=\sum_{i, j=1}^{3} u_{0, i}^{-} \bar{u}_{1, j}\left\langle\mathcal{A}_{l3}, \mathcal{A}_{i j}\right\rangle=u_{0,l}^{-} \bar{u}_{1,3}+u_{0,3}^{-} \bar{u}_{1,l}=0, \quad l=1,2,\nonumber\\
		\left\langle\mathcal{A}_{33}, (\mathbf{I}-\mathbf{P}) \bar{f}_{n-1}\right\rangle&=\sum_{i, j=1}^{3} u_{0, i}^{-} \bar{u}_{1, j}\left\langle\mathcal{A}_{33}, \mathcal{A}_{i j}\right\rangle=\frac{4}{3} u_{0,3}^{-} \bar{u}_{1,3}-\frac{2}{3}(u_{0,1}^{-} \bar{u}_{1,1}+u_{0,2}^{-} \bar{u}_{1,2})=-\frac{2}{3} u_{0,\sp}^{-}\cdot\bar{u}_{1,\sp}, \nonumber\\
		\left\langle\mathcal{B}_{3}, (\mathbf{I}-\mathbf{P}) \bar{f}_{n-1}\right\rangle&=\frac{5}{2}\bar{\t}_{1}u_{0,3}^{-}+\frac{5}{2} \theta_{0}^{-} \bar{u}_{1,3}=0,
\end{align}
where we have used $u_{0,3}^{-}=\bar{u}_{1,3}=0$. Substituting $\eqref{3.9}$ into $\eqref{3.6}$, one obtains
\begin{equation}\label{3.10}
		\bar{u}_{n-1,3}=\int_{y}^{\infty}\operatorname{div}_{\sp} \bar{u}_{n-2,\sp}(t,x_{\sp},s)ds, \quad
		\bar{\rho}_{n-1}+\bar{\theta}_{n-1}=\frac{2}{3} u_{0,\sp}^{-} \cdot \bar{u}_{1,\sp}.
\end{equation}

{\it Step 2.} We consider the case of $k\geq n-1$. Denoting $k=n-1+m$ with $m\geq 0$, one has from $\eqref{1.14}_5$ that
\begin{equation}\nonumber
	\left\{
	\begin{aligned}
		&\partial_{t}\int_{\R^3}\bar{F}_{m+1}dv+\operatorname{div}_{\sp}\int_{\R^3}v_{\sp} \bar{F}_{n-1+m}dv+\partial_{y}\int_{\R^3}v_{3}\bar{F}_{n+m}dv=0,\\
		&\partial_{t}\int_{\R^3}v_{i}\bar{F}_{m+1}dv+\operatorname{div}_{\sp}\int_{\R^3}v_{i}v_{\sp}\bar{F}_{n-1+m}dv+\partial_{y}\int_{\R^3}v_{i}v_{3}\bar{F}_{n+m}dv=0,\,\,i=1,2,3,\\
		&\partial_{t}\int_{\R^3}|v|^2\bar{F}_{m+1}dv+\operatorname{div}_{\sp}\int_{\R^3}|v|^2v_{\sp}\bar{F}_{n-1+m}dv+\partial_{y}\int_{\R^3}v_{3}|v|^2\bar{F}_{n+m}dv=0,
	\end{aligned}
	\right.
\end{equation}
which implies that
\begin{equation}\label{3.16}
	\left\{
	\begin{aligned}
		&\partial_{t}\bar{\rho}_{m+1}+\operatorname{div}_{\sp}\bar{u}_{n-1+m,\sp}+\partial_{y}\bar{u}_{n+m,3}=0,\\
		&\partial_{t}\bar{u}_{m+1,i}+\partial_{i}(\bar{\r}_{n-1+m}+\bar{\t}_{n-1+m})\\
		&\quad +\sum\limits_{j=1}^2\partial_{x_j}\la \mathcal{A}_{ij},(\mathbf{I-P})\bar{f}_{n-1+m}\ra +\partial_{y}\la \mathcal{A}_{3i},(\mathbf{I-P})\bar{f}_{n+m}\ra=0,\quad  (i=1,2),\\
		&\partial_{t}\bar{u}_{m+1,3}+\sum\limits_{j=1}^2\pa_{x_j}\la \mathcal{A}_{3j},(\mathbf{I-P})\bar{f}_{n-1+m}\ra\\
		&\quad+\partial_{y}(\bar{\r}_{n+m}+\bar{\t}_{n+m})+\partial_{y}\la \mathcal{A}_{33},(\mathbf{I-P})\bar{f}_{n+m}\ra =0,\\
		&3\partial_{t}(\bar{\r}_{m+1}+\bar{\t}_{m+1})+5\operatorname{div}_{\sp}\bar{u}_{n-1+m,\sp}
		\\&\quad +2\sum\limits_{i=1}^2\partial_{x_i}\la \mathcal{B}_{i},(\mathbf{I-P})\bar{f}_{n-1+m}\ra +5\partial_{y}\bar{u}_{n+m,3}+2\partial_{y}\la \mathcal{B}_{3},(\mathbf{I-P})\bar{f}_{n+m}\ra=0.
	\end{aligned}
	\right.
\end{equation}
Noting \eqref{3.1} and using similar arguments as in \eqref{2.7-3}-\eqref{2.7-4}, one has
\begin{equation}\label{3.16-1}
	\begin{aligned}
		&\begin{aligned}
			\la \mathcal{A}_{ij}, (\mathbf{I-P})\bar{f}_{n-1+m}\ra&=\la \mathcal{A}_{ij}, \mathbf{L}^{-1}\{\Gamma(f_{0}^{-},\mathbf{P}\bar{f}_{m+1})+\Gamma(\mathbf{P}\bar{f}_{m+1},f_{0}^{-})\}\ra+\la \mathcal{A}_{ij},W_{m}\ra\\
			&=\bar{u}_{m+1,j}u_{0,i}^{-}+u_{0,j}^{-}\bar{u}_{m+1,i}+\la \mathcal{A}_{ij},W_{m}\ra,\quad i\neq j,
		\end{aligned}\\
		&\begin{aligned}
			\la \mathcal{A}_{ii}, (\mathbf{I-P})\bar{f}_{n-1+m}\ra&=\la \mathcal{A}_{ii}, \mathbf{L}^{-1}\{\Gamma(f_{0}^{-},\mathbf{P}\bar{f}_{m+1})+\Gamma(\mathbf{P}\bar{f}_{m+1},f_{0}^{-})\}\ra+\la \mathcal{A}_{ii},W_{m}\ra\\
			&=\frac{4}{3}\bar{u}_{m+1,i}u_{0,i}^{-}-\frac{2}{3}\sum\limits_{j\neq i}u_{0,j}^{-}\bar{u}_{m+1,j}+\la \mathcal{A}_{ii},W_{m}\ra,\quad i=1,2,3,
		\end{aligned}\\
		&\begin{aligned}
			\la \mathcal{B}_{i}, (\mathbf{I-P})\bar{f}_{n-1+m}\ra&=\la \mathcal{B}_{i}, \mathbf{L}^{-1}\{\Gamma(f_{0}^{-},\mathbf{P}\bar{f}_{m+1})+\Gamma(\mathbf{P}\bar{f}_{m+1},f_{0}^{-})\}\ra+\la \mathcal{B}_{i},W_{m}\ra\\
			&=\frac{5}{2}\bar{u}_{m+1,i}\t_{0}^{-}+\frac{5}{2}\t_{m+1}u_{0,i}^{-}+\la \mathcal{B}_{i},W_{m}\ra,\quad i=1,2,3,
		\end{aligned}
	\end{aligned}
\end{equation}
where $W_{m}$ is the one defined in \eqref{3.15}.

For $(\mathbf{I}-\mathbf{P})\bar{f}_{n+m}$, we obtain from \eqref{3.1} that
\begin{align}\label{3.16-0}
		(\mathbf{I}-\mathbf{P})\bar{f}_{n+m}\equiv \mathbf{L}^{-1}\{&\Gamma(f_{0}^{-},\mathbf{P}\bar{f}_{m+2})+\Gamma(\mathbf{P}\bar{f}_{m+2},f_{0}^{-})\}+W_{m+1}\nonumber\\
		\equiv \mathbf{L}^{-1}\Big\{&-(\mathbf{I}-\mathbf{P})[ v_{3}\partial_y\mathbf{P}\bar{f}_{m+1}]+\big[\Gamma(f_{0}^{-},\mathbf{P}\bar{f}_{m+2})+\Gamma(\mathbf{P}\bar{f}_{m+2},f_{0}^{-})\big]\nonumber\\
		&+[\Gamma(\mathbf{P}f_{1}^{-},\mathbf{P}\bar{f}_{m+1})+\Gamma(\mathbf{P}\bar{f}_{m+1},\mathbf{P}f_{1}^{-})]\nonumber\\
		&+ y\cdot [\Gamma(\partial_{x_3}f_{0}^{-},\mathbf{P}\bar{f}_{m+1})+\Gamma(\mathbf{P}\bar{f}_{m+1},\partial_{x_3}f_{0}^{-})]\nonumber\\
		&+\sum\limits_{i+j=m+2\atop i=1\,\,\text{or}\,\,j=1}\frac{1}{2}[\Gamma(\bar{f}_{i},\mathbf{P}\bar{f}_{j})+\Gamma(\mathbf{P}\bar{f}_{j},\bar{f}_{i})]
		\Big\}+J_{m}\quad m\geq 0,
\end{align}
where $J_{m}$ is the one defined in \eqref{3.13}. Thus, similar as \eqref{3.16-1}, one has
\begin{equation}\label{3.16-2}
	\begin{aligned}
		&\la \mathcal{A}_{3i}, (\mathbf{I-P})\bar{f}_{n}\ra\\
		&=\left\la \mathcal{A}_{3i}, \mathbf{L}^{-1}\big\{-(\mathbf{I-P})v_3\pa_{y}\mathbf{P}\bar{f}_{1}\big\}\right\ra+\left\la\mathcal{A}_{3i},\mathbf{L}^{-1}\big\{\Gamma(f_{0}^{-},\mathbf{P}\bar{f}_{2})+\Gamma(\mathbf{P}\bar{f}_{2},f_{0}^{-})\big\}\right\ra\\
		&\quad+\left\la\mathcal{A}_{3i},\mathbf{L}^{-1}\big\{\Gamma(\mathbf{P}f_{1}^{-},\mathbf{P}\bar{f}_{1})+\Gamma(\mathbf{P}\bar{f}_{1},\mathbf{P}f_{1}^{-})\big\}\right\ra+ \left\la\mathcal{A}_{3i},\mathbf{L}^{-1}\big\{\Gamma(\bar{f}_{1},\bar{f}_{1})\big\}\right\ra\\
		&\quad +\left\la\mathcal{A}_{3i},\mathbf{L}^{-1}\big\{y\big[\Gamma(\partial_{x_3}f_{0}^{-},\mathbf{P}\bar{f}_{1})+\Gamma(\mathbf{P}\bar{f}_{1},\partial_{x_3}f_{0}^{-})\big]\big\}\right\ra+\la \mathcal{A}_{3i},J_{0}\ra\\
		&=-\lambda\pa_{y}\bar{u}_{1,i}+\bar{u}_{2,3}u_{0,i}^{-}+u_{0,3}^{-}\bar{u}_{2,i}+\bar{u}_{1,3}u_{1,i}^{-}+\bar{u}_{1,i}u_{1,3}^{-}\\
		&\quad +y\bar{u}_{1,3}\partial_{x_3}u_{0,i}^{-}+y\bar{u}_{1,i}\partial_{x_3}u_{0,3}^{-}+\bar{u}_{1,3}\bar{u}_{1,i}+\la\mathcal{A}_{3i},J_{0}\ra, \quad i=1,2,
	\end{aligned}
\end{equation}

\begin{equation}\label{3.16-5}
	\begin{aligned}
		&\la \mathcal{B}_{3}, (\mathbf{I-P})\bar{f}_{n}\ra\\
		&=\left\la \mathcal{B}_{3},\mathbf{L}^{-1}\big\{-(\mathbf{I-P})v_3\pa_{y}\mathbf{P}\bar{f}_{1}\big\}\right\ra+\left\la\mathcal{B}_{3},\mathbf{L}^{-1}\big\{\Gamma(f_{0}^{-},\mathbf{P}\bar{f}_{2})+\Gamma(\mathbf{P}\bar{f}_{2},f_{0}^{-})\big\}\right\ra\\
		&\quad+ \left\la\mathcal{B}_{3},\mathbf{L}^{-1}\big\{\Gamma(\mathbf{P}f_{1}^{-},\mathbf{P}\bar{f}_{1})+\Gamma(\mathbf{P}\bar{f}_{1},\mathbf{P}f_{1}^{-})\big\}\right\ra+\left\la\mathcal{B}_{3},\mathbf{L}^{-1}\big\{\Gamma(\bar{f}_{1},\bar{f}_{1})\big\}\right\ra\\
		&\quad +\left\la\mathcal{B}_{3},\mathbf{L}^{-1}\big\{y\big[\Gamma(\partial_{x_3}f_{0}^{-},\mathbf{P}\bar{f}_{1})+\Gamma(\mathbf{P}\bar{f}_{1},\partial_{x_3}f_{0}^{-})\big]\big\}\right\ra+\la \mathcal{B}_{3}, J_{0}\ra\\
		&=-\k\pa_{y}\bar{\t}_{1}+\frac{5}{2}\bar{u}_{2,3}\t_{0}^{-}+\frac{5}{2}u_{0,3}^{-}\bar{\t}_{2}+\frac{5}{2}\bar{u}_{1,3}\t_{1}^{-}+\frac{5}{2}\bar{\t}_{1}u_{1,3}^{-}\\
		&\quad +\frac{5}{2}y\bar{u}_{1,3}\partial_{x_3}\t_{0}^{-}+\frac{5}{2}y\bar{\t}_{1}\partial_{x_3}u_{0,3}^{-}+\frac{5}{2}\bar{u}_{1,3}\bar{\t}_{1}+\la\mathcal{B}_{3},J_{0}\ra,
	\end{aligned}
\end{equation}
\begin{equation}\label{3.16-4}
	\begin{aligned}
		&\la \mathcal{A}_{3i}, (\mathbf{I-P})\bar{f}_{n+m}\ra\\
		&=\left\la \mathcal{A}_{3i}, \mathbf{L}^{-1}\big\{-(\mathbf{I-P})v_3\pa_{y}\mathbf{P}\bar{f}_{m+1}\big\}\right\ra+\left\la\mathcal{A}_{3i},\mathbf{L}^{-1}\big\{\Gamma(f_{0}^{-},\mathbf{P}\bar{f}_{m+2})+\Gamma(\mathbf{P}\bar{f}_{m+2},f_{0}^{-})\big\}\right\ra\\
		&\quad+ \left\la\mathcal{A}_{3i},\mathbf{L}^{-1}\big\{\Gamma(\mathbf{P}f_{1}^{-},\mathbf{P}\bar{f}_{m+1})+\Gamma(\mathbf{P}\bar{f}_{m+1},\mathbf{P}f_{1}^{-})\big\}\right\ra\\
		&\quad+ \left\la\mathcal{A}_{3i},\mathbf{L}^{-1}\big\{\Gamma(\mathbf{P}\bar{f}_{1},\mathbf{P}\bar{f}_{m+1})+\Gamma(\mathbf{P}\bar{f}_{m+1},\mathbf{P}\bar{f}_{1})\big\}\right\ra\\
		&\quad +\left\la\mathcal{A}_{3i},\mathbf{L}^{-1}\big\{y\big[\Gamma(\partial_{x_3}f_{0}^{-},\mathbf{P}\bar{f}_{m+1})+\Gamma(\mathbf{P}\bar{f}_{m+1},\partial_{x_3}f_{0}^{-})\big]\big\}\right\ra+\la \mathcal{A}_{3i},J_{m}\ra\\
		&=-\lambda\pa_{y}\bar{u}_{m+1,i}+\bar{u}_{m+2,3}u_{0,i}^{-}+u_{0,3}^{-}\bar{u}_{m+2,i}+\bar{u}_{m+1,3}u_{1,i}^{-}+\bar{u}_{m+1,i}u_{1,3}^{-}+y\bar{u}_{m+1,3}\partial_{x_3}u_{0,i}^{-}\\
		&\quad +y\bar{u}_{m+1,i}\partial_{x_3}u_{0,3}^{-}+\bar{u}_{m+1,3}\bar{u}_{1,i}+\bar{u}_{1,3}\bar{u}_{m+1,i}+\la\mathcal{A}_{3i},J_{m}\ra, \quad m\geq 1,\,\,i=1,2,
	\end{aligned}
\end{equation}

\begin{equation}\label{3.16-3}
	\begin{aligned}
		&\la \mathcal{B}_{3}, (\mathbf{I-P})\bar{f}_{n+m}\ra\\
		&=\left\la \mathcal{B}_{3}, \mathbf{L}^{-1}\big\{-(\mathbf{I-P})v_3\pa_{y}\mathbf{P}\bar{f}_{m+1}\big\}\right\ra+\left\la\mathcal{B}_{3},\mathbf{L}^{-1}\big\{\Gamma(f_{0}^{-},\mathbf{P}\bar{f}_{m+2})+\Gamma(\mathbf{P}\bar{f}_{m+2},f_{0}^{-})\big\}\right\ra\\
		&\quad+ \left\la\mathcal{B}_{3},\mathbf{L}^{-1}\big\{\Gamma(\mathbf{P}f_{1}^{-},\mathbf{P}\bar{f}_{m+1})+\Gamma(\mathbf{P}\bar{f}_{m+1},\mathbf{P}f_{1}^{-})\big\}\right\ra\\
		&\quad+ \left\la\mathcal{B}_{3},\mathbf{L}^{-1}\big\{\Gamma(\mathbf{P}\bar{f}_{1},\mathbf{P}\bar{f}_{m+1})+\Gamma(\mathbf{P}\bar{f}_{m+1},\mathbf{P}\bar{f}_{1})\big\}\right\ra\\
		&\quad +\left\la\mathcal{B}_{3},\mathbf{L}^{-1}\big\{y\big[\Gamma(\partial_{x_3}f_{0}^{-},\mathbf{P}\bar{f}_{m+1})+\Gamma(\mathbf{P}\bar{f}_{m+1},\partial_{x_3}f_{0}^{-})\big]\big\}\right\ra+\la \mathcal{B}_{3},J_{m}\ra\\
		&=-\k\pa_{y}\bar{\t}_{m+1}+\frac{5}{2}\bar{u}_{m+2,3}\t_{0}^{-}+\frac{5}{2}u_{0,3}^{-}\bar{\t}_{m+2}+\frac{5}{2}\bar{u}_{m+1,3}\t_{1}^{-}+\frac{5}{2}\bar{\t}_{m+1}u_{1,3}^{-}\\
		&\quad +\frac{5}{2}y\bar{u}_{m+1,3}\partial_{x_3}\t_{0}^{-}+\frac{5}{2}y\bar{\t}_{m+1}\partial_{x_3}u_{0,3}^{-}+\frac{5}{2}\bar{u}_{m+1,3}\bar{\t}_{1}+\frac{5}{2}\bar{u}_{1,3}\bar{\t}_{m+1}+\la\mathcal{B}_{3i},J_{m}\ra, \quad m\geq 1,
	\end{aligned}
\end{equation}
and
\begin{equation}\label{3.16-6}
	\begin{aligned}
		&\la \mathcal{A}_{33}, (\mathbf{I-P})\bar{f}_{n+m}\ra=\la \mathcal{A}_{33},W_{m+1}\ra+\left\la\mathcal{A}_{33},\mathbf{L}^{-1}\big\{\Gamma(f_{0}^{-},\mathbf{P}\bar{f}_{m+2})+\Gamma(\mathbf{P}\bar{f}_{m+2},f_{0}^{-})\big\}\right\ra\\
		&=\frac{4}{3}\bar{u}_{m+2,3}u_{0,3}^{-}-\frac{2}{3}(u_{0,\sp}^{-}\cdot \bar{u}_{m+2,\sp})+\la \mathcal{A}_{33},W_{m+1}\ra,
	\end{aligned}
\end{equation}

Substituting \eqref{3.16-1}-\eqref{3.16-6} into \eqref{3.16} and using $u_{0,3}^{-}=\bar{u}_{1,3}=\operatorname{div}_{x}u_{0}=0$, after some simple but tedious calculations, one can deduce \eqref{3.12} with $W_{k-1}, J_{k-1}$ and $H_{k-1}$ defined in \eqref{3.15}-\eqref{3.17} respectively.

We shall point out that $\mathfrak{f}_{k-1,1}^{\pm},\mathfrak{f}_{k-1,2}^{\pm},\mathfrak{g}_{k-1}^{\pm}$ and $\mathfrak{s}_{k-1}^{\pm}$ in \eqref{3.12} depends on $f_{i}(0\leq i\leq k)$ and $\bar{f}_{i}(1\leq i\leq k-1)$. Indeed, it is clear that $\mathfrak{s}_{k-1}^{\pm}$ depends on $\bar{f}_{j}^{\pm}$ up to $(k-1)$-th order since $W_{k-1}^{\pm}$ and $H_{k-1}^{\pm}$ from \eqref{3.1} and $n\geq 3$. Using \eqref{3.21}, we see the terms involving $\bar{u}_{k,3}^{\pm}$ in $\mathfrak{f}_{k-1,1}^{\pm},\mathfrak{f}_{k-1,2}^{\pm}$ and $\mathfrak{g}_{k-1,}^{\pm}$ are determined by $\bar{f}_{k-1}^{\pm}$ and $\bar{f}_{k+1-n}^{\pm}$. Hence, we only have to show $\la \mathcal{A}_{31},J_{k-1}^{\pm}\ra$, $\la \mathcal{A}_{32},J_{k-1}^{\pm}\ra$ and $\la \mathcal{B}_{3},J_{k-1}^{\pm}\ra$ are determined by $\bar{f}_{j}^{\pm}(1\leq j\leq  k-1)$. For this,  we look at the highest order terms in those terms, that is,
$$
\begin{aligned}
	&I_{1,i}:=\la \mathcal{A}_{3i}, \mathbf{L}^{-1}\{\Gamma(f_{0}^{\pm},(\mathbf{I-P})\bar{f}_{k+1}^{\pm})+\Gamma((\mathbf{I-P})\bar{f}_{k+1}^{\pm},f_{0}^{\pm})\}\ra,\quad i=1,2,\\
	&I_{2,i}:=\la \mathcal{A}_{3i},\mathbf{L}^{-1}\{\Gamma((\mathbf{I-P})f_{1}^{\pm},\mathbf{P}\bar{f}_{k}^{\pm})+\Gamma(\mathbf{P}\bar{f}_{k}^{\pm},(\mathbf{I-P})f_{1}^{\pm})\}\ra,\quad i=1,2,\\
	&I_{3,1}:=\la \mathcal{B}_{3}, \mathbf{L}^{-1}\{\Gamma(f_{0}^{\pm},(\mathbf{I-P})\bar{f}_{k+1}^{\pm})+\Gamma((\mathbf{I-P})\bar{f}_{k+1}^{\pm},f_{0}^{\pm})\}\Big\ra,\\
	&I_{3,2}:=\la \mathcal{B}_{3}, \mathbf{L}^{-1}\{\Gamma((\mathbf{I-P})f_{1}^{\pm},\mathbf{P}\bar{f}_{k}^{\pm})+\Gamma(\mathbf{P}\bar{f}_{k}^{\pm},(\mathbf{I-P})f_{1}^{\pm})\}\ra.
\end{aligned}
$$
If $n\geq 4$, then it follows from $\eqref{2.1}_1$ that $(\mathbf{I-P})f_{1}=0$ and hence $I_{2,i}=I_{3,2}=0$. Moreover, noting \eqref{3.1}, one can deduce that $I_{1,i}, I_{3,1}$ are determined by $\bar{f}_{j}^{\pm}$ with $1\leq j\leq k-1$.

If $n=3$, it follows form $\eqref{2.1}_2$ that
$$
\begin{aligned}
	(\mathbf{I-P})f_{1}^{\pm}&
	=\frac{1}{2}(\mathbf{I-P})(\frac{f_{0}^{\pm}f_{0}^{\pm}}{\sqrt{\mu}})\\
	&=\frac{1}{2}\sum\limits_{l,m=1}^3u_{0,l}^{\pm}u_{0,m}^{\pm}\mathcal{A}_{lm}+\sum\limits_{l=1}^3\t_{0}^{\pm}u_{0,l}^{\pm}\mathcal{B}_{l}+\frac{(\t_{0}^{\pm})^2}{8}(\mathbf{I-P})\{(|v|^2-5)^2\sqrt{\mu}\}.
\end{aligned}
$$
Since $u_{0,3}^{\pm}=0$ and $\mathbf{L}^{-1}$ preserves the odevity of $v$, the nonvanishing parts in $I_{2,i}, I_{3,2}$ depend only on $\bar{u}_{k,3}^{\pm}$, which is then determined by $\bar{f}_{k-1}^{\pm}$ and $\bar{f}_{k+1-n}^{\pm}$ Using \eqref{3.1} and a similar argument, we know that $I_{1,i},I_{3,1}$ are also determined by $\bar{u}_{k,3}^{\pm}$ and hence $\bar{f}_{j}^{\pm}(1\leq j\leq k-1)$. This completes the proof of Proposition \ref{prop3.1}. $\hfill\square$

\

\noindent{\bf Acknowledgments.} Feimin Huang's research is partially supported by National Key R\&D Program of China No. 2021YFA1000800, and National Natural Sciences Foundation of China No. 12288201. Yong Wang's research is partially supported by National Natural Science Foundation of China No. 12022114, 12288201, CAS Project for Young Scientists in Basic Research, Grant No. YSBR-031, and Youth Innovation Promotion Association of CAS No. 2019002. Feng Xiao's research is partially supported by National Natural Science Foundation of China No. 12201209.

\end{document}